\documentclass[12pt]{jmlr}
\usepackage{isipta2023}
\usepackage[british]{babel}
\usepackage{booktabs} %
\usepackage{tikz} %

\usepackage{sepfootnotes}

\usepackage{adjustbox}

\newtheorem{assumption}[theorem]{Assumption}
\usepackage{enumitem}
\usepackage{bbm}
\usepackage{mathrsfs}
\usepackage{textgreek} %

\usepackage{comment}

\DeclarePairedDelimiter\abs{\lvert}{\rvert}%
\DeclarePairedDelimiter\norm{\lVert}{\rVert}%
\makeatletter
\let\oldabs\abs
\def\abs{\@ifstar{\oldabs}{\oldabs*}}
\let\oldnorm\norm
\def\norm{\@ifstar{\oldnorm}{\oldnorm*}}
\makeatother

\newcommand{\ceil}[1]{\left\lceil #1 \right\rceil}

\renewcommand{\epsilon}{\varepsilon}

\newcommand{\dcor}{d_{\textnormal{\textsc{COR}}}}
\newcommand{\dMult}{d_{\textnormal{\textsc{Mult}}}}

\newcommand{\dHam}{d_{\textsc{Ham}}}
\newcommand{\dSym}{d_{\triangle}}

\newcommand{\dLV}{d_{\textnormal{\textsc{LV}}}}
\newcommand{\dist}{d_{\textnormal{\textsc{dist}}}}

\newcommand{\IInt}[1]{\mathcal I \hspace{-0.1em}\left( #1 \right)}
\newcommand{\IIntP}[1]{\mathcal I_1 \left( #1 \right)}

\newcommand{\T}{\mathcal{T}}
\newcommand{\To}{\mathcal{T}^{\mathrm{o}}}
\newcommand{\F}{\mathscr{F}}
\newcommand{\X}{\mathcal{X}}
\newcommand{\supp}{\mathrm{supp}}

\newcommand{\ComX}{\textnormal{Comp}(\X)}

\DeclareMathOperator*{\esssup}{ess\,sup}
\DeclareMathOperator*{\essinf}{ess\,inf}

\newcommand{\Unif}{\mathrm{Unif}}

\newcommand{\Spairs}{\mathbb S}
\newcommand{\Puff}{$\epsilon$-$\mathcal P$uffer$\mathcal F$ish($\mathbb D$, $\Spairs$)}
\newcommand{\dPuff}{d_{\mathbb D, \Spairs}}
\newcommand{\GPuff}{G_{\mathbb D, \Spairs}}
\newcommand{\drm}{d_{\textnormal{\textsc{DR}}}}

\newcommand{\DPP}{M_G}

\newcommand\iid{\stackrel{iid}{\sim}} %

\let\originalleft\left
\let\originalright\right
\renewcommand{\left}{\mathopen{}\mathclose\bgroup\originalleft}
\renewcommand{\right}{\aftergroup\egroup\originalright}

\title[General Inferential Limits Under Differential and Pufferfish Privacy]{General Inferential Limits Under \\ Differential and Pufferfish Privacy}
\author{%
  \Name{James Bailie}\Email{jamesbailie@g.harvard.edu}\\ %
  \addr  Department of Statistics, Harvard University, USA
  \AND
  \Name{Ruobin Gong}\Email{ruobin.gong@rutgers.edu}\\
  \addr Department of Statistics, Rutgers University, USA
}

\begin{document}
\maketitle

\sepfootnotecontent{abstractFootnote}{This is an extended version of the conference paper \cite{bailieDifferentialPrivacyGeneral2023}. Along with some lesser additions, Sections~\ref{sec:pufferfish} and~\ref{sec:pufferfish-IP} on Pufferfish privacy are new; and Lemma~\ref{lemmaMarginalData}, Theorem~\ref{thmMarginalData} (these two results correspond to Theorem~6 in \cite{bailieDifferentialPrivacyGeneral2023}), Theorem~\ref{thmDPHypothesisTesting} (Theorem~8 in \cite{bailieDifferentialPrivacyGeneral2023}) and Theorem~\ref{thmPriorPredictive} (Theorem~10 in \cite{bailieDifferentialPrivacyGeneral2023}) all contain minor corrections.}

\begin{abstract}%
Differential privacy (DP) is a class of mathematical standards for assessing the privacy provided by a data-release mechanism. This work concerns two important flavors of DP that are related yet conceptually distinct: pure $\epsilon$-differential privacy ($\epsilon$-DP) and Pufferfish privacy. We restate $\epsilon$-DP and Pufferfish privacy as Lipschitz continuity conditions and provide their formulations in terms of an object from the imprecise probability literature: the interval of measures. We use these formulations to 
derive limits on key quantities in frequentist hypothesis testing and in Bayesian inference using data that are sanitised according to either of these two privacy standards. Under very mild conditions, the results in this work are valid for arbitrary parameters, priors and data generating models. These bounds are weaker than those attainable when analysing specific data generating models or data-release mechanisms. However, they provide generally applicable limits on the ability to learn from differentially private data -- even when the analyst's knowledge of the model or mechanism is limited. They also shed light on the semantic interpretations of the two DP flavors under examination, a subject of contention in the current literature.\sepfootnote{abstractFootnote}
\end{abstract}


\section{Introduction}\label{sec:intro}

The world today is witnessing an explosive growth of large-scale datasets containing personal information. Demographic and economic surveys, biomedical studies and massive online service platforms facilitate understanding of human biological functions and socio-behavioural environments. At the same time, they pose the risk of exposing confidential information about data contributors. Breaches of privacy can happen counter-intuitively and without malice. For example, \cite{homer2008resolving} demonstrated that even coarsely aggregated SNP (single-nucleotide polymorphisms \citep{kimSNPGenotypingTechnologies2007}) data from genome-wide association studies (GWAS) can still reliably reveal individual participants. This unsettling revelation led to the decision by the U.S. National Institute of Health to remove aggregate SNP data from open-access databases \citep{yu2013stability}. This incident, and similar occurrences across science, government and industry \cite{narayananRobustDeanonymizationLarge2008, dworkExposedSurveyAttacks2017, culnaneStopOpenData2019, sweeneySimpleDemographicsOften2000}, have attracted public attention and sparked debate about privacy-preserving data curation and dissemination. 

Commensurate with the increasing risk of privacy breaches, the recent decades have also seen rapid advances in formal approaches to statistical disclosure limitation (SDL).
These methodologies
supply a solid mathematical foundation for endeavors that enhance confidentiality protection without undue sacrifice to data quality. Notably, \emph{differential privacy} (DP) \citep{dwork2006calibrating,bunConcentratedDifferentialPrivacy2016,kifer2014pufferfish} puts forth a class of rigorous and practical standards for assessing the level of privacy provided by a data release. Many large IT companies, including Google \citep{erlingsson2014rappor}, Apple \citep{apple}, and Microsoft \citep{microsoft}, have been early adopters of DP. More recently, the U.S. Census Bureau deployed DP to protect the data publications of the 2020 Decennial Census \citep{abowd20222020}. The U.S. Internal Revenue Service is also exploring differentially private synthetic data methods for the publication of individual tax data \citep{bowen2022synthetic}. 
These decisions by statistical agencies and corporations showcase the growing popularity of DP among major data curators. 

Innovations in privacy protection methods have prompted quantitative researchers to confront a new reality, as existing modes, practices and expectations of data access are subject to renewal. We highlight two points of tension in this development. First, DP promises \emph{transparency}, in the sense that the design details about the protection method can be made public without degrading DP's mathematical assessment of the level of privacy protection. Transparency is one of the advantages of DP over traditional SDL\label{review2.1} methods since it supports valid statistical inference by providing the analyst with the ability to model the privacy noise. 
However, this promise often falls short in practice, leaving the statistician with tied hands \cite{gongTransparentPrivacyPrincipled2022}. Second, following the high-profile adoption of DP by the U.S. Census Bureau, a debate ensued concerning its interpretation, or its \emph{semantics}, as well as its reconciliation with other notions of statistical disclosure risk; see e.g. \cite{kenny2021use,hotz2022balancing,kifer2022bayesian, jarmin2023depth, boyd2022differential, francisNoteMisinterpretationUS2022, muralidharDatabaseReconstructionNot2023, garfinkelCommentMuralidharDomingoFerrer2023, sanchezExaminationAllegedPrivacy2023, kellerDatabaseReconstructionDoes2023}. These issues motivate theoretical investigations that may shed light on the pragmatic translation between rigorous privacy standards and usable statistical advice.

The current work takes multiple steps toward the resolution of these debates by examining DP via the lens of imprecise probabilities (IP). Our focus is restricted to two important flavors of DP: 1) the classic notion of pure $\epsilon$-differential privacy ($\epsilon$-DP) \cite{dwork2006calibrating}, and 2) Pufferfish privacy \cite{kifer2014pufferfish}, a conceptually-distinct variant of $\epsilon$-DP that is showing semantic promise (see e.g. \cite{jarmin2023depth}).  We  begin by describing $\epsilon$-DP as a Lipschitz continuity condition (\sectionref{sec:background}). This description enables the interpretation of $\epsilon$-DP as an \emph{interval of measures} \citep{derobertis1981bayesian} induced by the data-release mechanism (\sectionref{sec:DPIoM}). From here, we derive some implications of this interpretation on the problem of statistical inference using privacy-protected data releases. These results concern the probability model of the observable privatised data (\sectionref{sec:marginal}), as well as frequentist hypothesis testing (\sectionref{sec:frequentist}) and Bayesian posterior inference (\sectionref{sec:Bayesian}) using these data. Next we turn to address Pufferfish privacy (\sectionref{sec:pufferfish}) -- showing that it too can be described as a Lipschitz continuity condition -- and discuss its semantic interpretation as limits to frequentist and Bayesian inferences in an analogous manner (\sectionref{sec:pufferfish-IP}). Further, we link Pufferfish to another IP object: the density ratio neighbourhood (Theorem~\ref{thm:pufferfish-drn}). 
The results in Sections~\ref{sec:marginal}-\ref{sec:pufferfish-IP} establish bounds on key inferential objects while having general validity under mild assumptions about the data model, the privacy mechanism, and (when applicable) the analyst's prior. %
\sectionref{secTight} demonstrates that these results are optimal in the sense that the bounds we obtain are pointwise tight. \sectionref{sec:discussion} concludes the paper with a discussion.

Throughout this work, we demonstrate that various objects from the imprecise probability literature naturally arise when studying differential privacy. Specifically, definitions of DP often invoke \emph{distortion models}: neighbourhoods of precise probabilities defined as closed balls with respect to some metric -- or, more correctly, some \emph{distorting function} \cite{montesUnifyingNeighbourhoodDistortion2020,montesUnifyingNeighbourhoodDistortion2020a}. Moreover, the choice of the distortion model (partially) determines the flavor of the resulting privacy guarantee \cite{bailieRefreshmentStirredNot2024a}. In the following sections, we will outline the appropriate distortion models formulations as they arise.  

DP objects %
are naturally amenable to IP analysis. Indeed, the rich vocabulary of IP can help to articulate the %
properties of a DP object in a precise yet general manner. Within the current literature on DP, the study of data privacy protection using IP is a nascent endeavour. 
\cite{komarova2020identification} examines the issue of partial identification in inference from privacy-protected data, where in certain situations the identification set can be described with a belief function. 
In \cite{li2022local},\label{review1} the authors formulate local differential privacy definitions for belief functions, a proposal that amounts to a set-valued SDL mechanism whose probability distribution is given by the mass function associated with a belief function. \cite{liu2023two} examines constraints on DP mechanisms in terms of belief revision and updating. On the matter of using IP to explore mathematical formalisations of data privacy, we will return and remark on a few concrete potential future directions in the discussion (\sectionref{sec:discussion}).

\section{Pure \texorpdfstring{$\epsilon$}{ϵ}-Differential Privacy}\label{sec:background}

\sepfootnotecontent{metric}{Throughout this work, we allow metrics to have codomain $[0,\infty]$ rather than the more standard $[0,\infty)$. We precisely define a metric in Definition~\ref{defnMetric} of Appendix~\ref{appendixMetric}.}

Define the data universe $\mathcal X$ as the set of all theoretically-possible observable datasets. Let $d$ be a metric on $\mathcal X$.\sepfootnote{metric} Given confidential data $x \in \mathcal X$, consider releasing some (potentially randomised) summary statistic $T \in \mathcal T$ of $x$. To formalise this, equip the set $\mathcal T$ with a $\sigma$-algebra $\mathscr F$ and define a \emph{data-release mechanism} as a function $M : \mathcal X \times [0,1] \to \mathcal T$ whose inputs are the confidential data $x$ and a random \emph{seed} $U \in [0,1]$, and whose output is the summary statistic $T$. (We require that $M(x, \cdot)$ is $(\mathcal B[0,1], \mathscr F)$-measurable for each $x \in \mathcal X$, where $\mathcal B[0,1]$ denotes the Borel $\sigma$-algebra on $[0,1]$.) 

A distribution on the seed $U$ induces a probability on the summary statistic $T=M(x, U)$. Without loss of generality, we may take $U \sim \mathrm{Unif}[0,1]$. Denote by $P_x$ the probability measure of $M(x,U)$ induced by $U$, taking $x$ as fixed:
\begin{equation}\label{eqDefnPx}
    P_x(M(x,U) \in S) = \lambda \big( \{ u \in [0,1] : M(x,u) \in S\} \big),
\end{equation}
where $\lambda$ is the Lebesgue measure on $[0,1]$.

The realised value of the seed $U$ and the observed dataset $x$ are assumed to remain secret, while all other details of $M$ (including the distribution of $U$) may -- and should -- be made public \citep{gongTransparentPrivacyPrincipled2022}. An attacker is tasked with inferring $x$ (or some summary $h(x)$ of $x$) based on observing a draw $T = M(x,U) \sim P_x$. This set-up is analogous to fiducial inference \cite{hannig2016generalized}, with $x$ taking the role of the parameters, $T$ the data, and $M$ the data-generating process.

Pure $\epsilon$-DP is a Lipschitz condition on $M$:

\sepfootnotecontent{multDist}{\label{review2.3}In defining $\dMult$ we set $0/0 = \infty/\infty = 1$; and on the RHS of~\eqref{eqDPDefnLip}, we set $0 \times \infty = \infty$. On the space of probability measures, $\dMult$ is \emph{strongly equivalent} (Definition~\ref{defnMetricStrongEquivalence} in Appendix~\ref{appendixMetric}) to the \emph{density ratio metric} $\drm$ \citep{wasserman1992invariance} (Definition~\ref{defnDensityRatio} below). Namely,
	\begin{equation}\label{eqMultDRMEquivalent}
		\dMult(P,Q) \le \drm(P,Q) \le 2 \dMult(P,Q),
	\end{equation}
	for \emph{probability} measures $P,Q$ on $(\mathcal T, \mathscr F)$, so that $\epsilon$-DP can be defined with $\drm$ in place of $\dMult$, up to rescaling of $\epsilon$. Equation~\eqref{eqMultDRMEquivalent} is proven in Proposition~\ref{propInequalityMultDRM}.}

\begin{definition}\label{defnDP}
	Given a data universe $\mathcal X$ equipped with a metric $d$, a data-release mechanism $M : \mathcal X \times [0,1] \to \mathcal T$ satisfies \emph{(pure) $\epsilon$-differential privacy} if, for all $x, x' \in \mathcal X$,
	\begin{equation}\label{eqDPDefnLip}
		\dMult (P_x, P_{x'}) \le \epsilon d(x, x'),
	\end{equation}
	where \label{review2}
	\[\dMult(\mu,\nu) = \sup_{S \in \mathscr F} \abs{ \ln \frac{\mu(S)}{\nu(S)}},\]%
	is the \emph{multiplicative distance\sepfootnote{multDist}} between measures $\mu,\nu$ on $(\mathcal T, \mathscr F)$.
\end{definition}

The smallest Lipschitz constant $\epsilon \ge 0$ which satisfies~\eqref{eqDPDefnLip} is called the \emph{privacy loss} associated with releasing $T$. Larger $\epsilon$ intuitively corresponds to less privacy; smaller $\epsilon$ gives stronger privacy protection. A tenet of DP (in contrast with many other statistical disclosure risk frameworks) is that dependence of $M(x,U)$ on $x$ implies non-negligible privacy loss $\epsilon > 0$: Since $\dMult(\mu,\nu) = 0$ if and only if $\mu = \nu$, \label{statementCompletePrivacy} complete privacy\label{review3} ($\epsilon = 0$) is only possible by releasing pure noise -- or, more exactly, by releasing $T \sim P_x$ where $P_x$ is a function of $x$ only through its \emph{connected component} $[x] = \{x' \in \mathcal X \mid d(x, x') < \infty\}$ (see Definition~\ref{defnConnected} below). (This statement is formalised in Proposition~\ref{propCompletePrivacy}.)

In the ideal case, the data custodian decides upon a maximum value of $\epsilon$ which is acceptable when considering the sensitivity of the data $x$ and the privacy protection they deserve. The data custodian then designs a data-release mechanism which satisfies $\epsilon$-DP, for this chosen value of $\epsilon$. From this perspective, the maximum acceptable value of $\epsilon$ is called the \emph{privacy loss budget}.

Two common choices of the metric $d$ on $\mathcal X$ are: \label{review2.2}
\begin{enumerate}[label=\Alph*)]
    \item the Hamming distance
	\[\dHam(x,x') = \begin{cases}
		\sum_{i=1}^n 1_{x_i \ne x_i'} & \text{if } \abs{x} = \abs{x'} = n, \\
		\infty & \text{otherwise,}
	\end{cases}\]
where the data $x = (x_1,x_2,\ldots,x_n)$ are vectors and $\abs{x}$ is the size of $x$; and 
    \item the symmetric difference metric
	\[\dSym(x, x') = \abs{ x \setminus x'} + \abs{ x' \setminus x},\]
where the data $x, x' \in \mathcal X$ are multisets and $x \setminus x'$ is the (multi-)set difference.\footnote{\label{review4}We formally define a multiset $S$ to be a non-negative-integer-valued function $c_S$, where $c_S(a)$ is the number of times the element $a$ appears in $S$. The multiset difference $S \setminus S'$ is defined as the function $c_{S \setminus S'}(a) = \max\{0, c_S(a) - c_{S'}(a)\}$ and the multiset cardinality is defined as $\abs{S} = \sum_a c_S(a)$.%
}
\end{enumerate}

Equation \eqref{eqDPDefnLip} with $d$ the Hamming distance is referred to as \emph{bounded} DP and with the symmetric difference as \emph{unbounded} DP.

The intuition behind differential privacy considers each record $x_i$ in the data $x$ as representing a single distinct individual. A distance $d(x,x') = 1$ then implies that $x$ and $x'$ differ according to the change in behaviour of a single individual -- a change in the individual's response, for the Hamming distance; or a change in whether the individual responds or not, for the symmetric difference metric. Specifically, $\epsilon$-DP implies that a single individual can change the summary statistic $T=M(x,U)$ by at most $\epsilon$, where ``change'' is interpreted probabilistically in terms of the multiplicative distance. 

Under the mild assumption that $d$ is a graph distance with unit edges (Assumption~\ref{assumptionGraphDistance}, given below in Section~\ref{sec:DPIoM}), the converse implication also holds. That is, $\epsilon$-DP is equivalent to the Lipschitz condition \eqref{eqDPDefnLip} holding when $d(x,x') = 1$. (This follows by the triangle inequality; for details see the proof of Theorem~\ref{thm:dp-iom}.) From herein, we restrict our attention to the set of such metrics, which includes $\dHam$ and $\dSym$.

Since DP controls the change in $t$ due to perturbations in the data $x$, it can naturally be understood as a robustness property \citep{dwork2009differential_lei, avella-medinaRoleRobustStatistics2020, avella-medinaPrivacyPreservingParametricInference2021,  asiRobustnessPrivacyBack2023,  hopkinsRobustnessImpliesPrivacy2023}. Measuring the ``change in $t$'' by the multiplicative distance $\dMult$ -- in place of more-familiar metrics typically seen in the robustness literature, such as the Kolmogorov distance or the total variation distance -- is motivated by the strong notion of privacy as \emph{indistinguishability}. The formulations of both pure $\epsilon$-DP and Pufferfish as intervals of measure (which we describe in later sections) make this motivation clear; we therefore postpone further discussion of indistinguishability to Remarks~\ref{remarkIndistinguishability} and~\ref{remarkPuffIndistinguish}.

This link to robustness hints at the connection between DP and IP. As we will soon see in Section~\ref{sec:DPIoM}, the multiplicative distance $\dMult$ is a distorting function \cite{montesUnifyingNeighbourhoodDistortion2020}, and consequently pure $\epsilon$-DP can be characterised in terms of a distortion model, a point we expand on later in Remark~\ref{remarkDPCredalSet}. Indeed, $\dMult$ satisfies many of the common desiderata for distorting functions: it is positive definite, symmetric and quasi-convex, and it satisfies the triangle inequality. Additionally, $\dMult(\mu,\nu)$ is continuous with respect to the supremum norm if and only if $\T$ has finite cardinality and $\mu$ and $\nu$ have support $\T$. See Appendix~\ref{appDistort} for details on these desiderata.

\begin{example}[Laplace mechanism \citep{dwork2006calibrating}]\label{ex:laplace}
	Consider the problem of releasing a sanitised version of a deterministic summary statistic $q : \mathcal X \to \mathbb R^k$. (The terms `sanitised,' `privatised,' `privacy-protected,' and `privacy-preserving' are synonymous in the DP literature.) The Laplace mechanism adds noise with standard deviation proportional to the \emph{global $\ell_1$-sensitivity} of $q$:
    \begin{equation}\label{eqDefnGlobalSensitivity}
        \Delta(q) = \sup_{\substack{x, x' \in \mathcal X \\ d(x,x') = 1}} \norm{q(x) - q(x')}_1.
    \end{equation}
	Specifically, define $M(x,L) = q(x) + bL$, where $b =  \frac{\Delta(q)}{\epsilon}$ and $L$ (the seed) is a $k$-vector of i.i.d. Laplace random variables with density $f(z) = 0.5 \exp (- \abs{z})$. (See Figure~\ref{fig:Laplace_mech_ex}.) When $d(x, x') = 1$,
	\begin{align*}
		P_x(&S_1 \times \ldots \times S_k) = \prod_{i=1}^k \left[ 0.5 b^{-1} \int_{S_i} \exp \left( - \frac{\abs{z - q_i(x)}}{b} dz\right)\right] \\
		&\ge \exp \left( \frac{-\Delta(q)}{b} \right)  \prod_{i=1}^k \left[ 0.5 b^{-1} \int_{S_i} \exp \left( - \frac{\abs{z - q_i(x')}}{b} dz\right)\right] \\
		&= \exp (-\epsilon) P_{x'} (S_1 \times \ldots \times S_k).
	\end{align*}
	We will see in Theorem \ref{thm:dp-iom} that this suffices to prove that $M$ is $\epsilon$-DP.
\end{example}

\begin{figure}
	\floatconts{fig:Laplace_mech_ex}
	{\caption{\emph{An illustration of the Laplace mechanism (Example~\ref{ex:laplace})}. %
    Here $p_x(t)$ and $p_{x'}(t)$ are probability densities for the Laplace mechanism's output $T$. The datasets $x$ and $x'$ are chosen so that $d(x, x') = 1$ and $ q(x') - q(x) = \Delta(q) = 1$. The annotations $0.05$ and $0.05 \exp(\epsilon)$ are demonstrative of the property: $p_x(t) = e^\epsilon p_{x'}(t)$ for all $t \le q(x)$ and $p_x(t) = e^{-\epsilon} p_{x'}(t)$ for all $t \ge q(x')$. (In this illustration, $\epsilon = 2$ and $k=1$.)}}
	{\includegraphics[width=\textwidth]{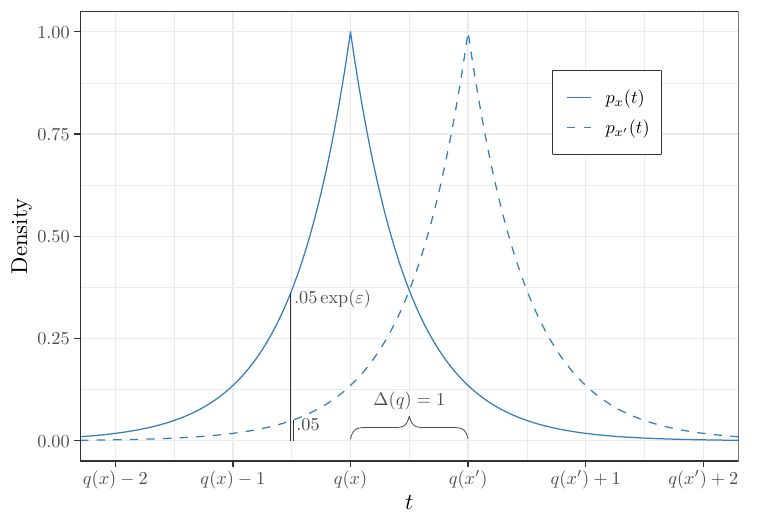}}
\end{figure}

\begin{example}[randomised response \citep{warnerRandomizedResponseSurvey1965}]\label{exampleRandomisedResponse}
    Taking $\mathcal X = \bigcup_{n \in \mathbb N} \{0,1\}^n$ as the data universe, the randomised response mechanism $M$ flips each bit $x_i$ with probability $p = (\exp \epsilon + 1)^{-1}$. That is, given a binary $n$-vector $x$ as input, $M$ outputs another binary $n$-vector $T$ with $i$-th component $T_i = x_i + B_i \mod 2$ where $B_1,B_2,\ldots \iid \mathrm{Bernoulli}(p)$. This mechanism is $\epsilon$-DP when $d=\dHam$. 
	
    Moreover, $M$ conforms with the \emph{local} model of DP \citep{kasiviswanathanWhatCanWe2011, duchi2018minimax}, since each data point can be independently infused with noise by the data respondents themselves. (For example, the $i$-th respondent can flip their own coin $B_i$ and report their noisy answer $T_i$.) Local DP models are typically used when data are collected by an untrusted entity (such as an IT company), since these models require that the privacy protection is applied to each record before data collection. In the non-interactive setting, this requirement implies that $P_x(T \in \cdot)$ must factor as a product measure $\prod_{i=1}^n P_{x_i}(T_i \in \cdot)$, where $n = \abs{x}$ is the number of records in $x$%
    . (In contrast, the local \emph{interactive} model of DP has the weaker condition that the distribution of user $i$'s response $T_i$ cannot depend on $x_j$ for $j \ne i$ (like in the non-interactive setting) but that this distribution can depend on the previous users' responses $T_j$ for $j < i$.) The local privacy model contrasts with the \emph{central} privacy model, under which the raw data $x$ can be aggregated by a central, trusted authority (such as a national statistical office) before privacy protection is applied. (And hence, the probability $P_x$ of a central privacy mechanism need not be factorizable.)
\end{example}

\section{Pure \texorpdfstring{$\epsilon$}{ϵ}-Differential Privacy as an Interval of Measures}\label{sec:DPIoM}

We introduce the definition of an \emph{interval of measures} (IoM), due to \citet{derobertis1981bayesian}:
\begin{definition}\label{def:iom}
	For measures $\mu$ and $\nu$ on the measurable space  $\left(\T, \mathscr{F}\right)$, write $\mu \le \nu$ to denote that $\mu(S) \le \nu(S)$ for all $S \in \mathscr F$. 
	
	Given measures $L, U$ on $\left(\T, \mathscr{F}\right)$ with $L\le U$, the convex set of measures
	\begin{equation*}
		\IInt{L,U} = \{ \mu \mathrm{\ a\ measure\ on\ } \left(\T, \mathscr{F}\right) \mid L \le \mu \le U \} 
	\end{equation*}
	is called an \emph{interval of measures}. $L$ and $U$ are called the \emph{lower} and \emph{upper measures} of $\IInt{L,U}$.
\end{definition}

\sepfootnotecontent{credal}{A credal set is simply a set of probabilities.}

Let $\Omega$ be the collection of all $\sigma$-finite measures on $\left(\T, \mathscr{F}\right)$ and $\Omega_1 = \{ P \in \Omega \mid P(\T) = 1\}$ be the collection of all probability measures on $(\T, \mathscr F)$. In the vast majority of cases we encounter (and indeed all the practically meaningful ones), the upper measure $U$ of an IoM $\IInt{L,U}$ is $\sigma$-finite and hence $\IInt{L,U} \subset \Omega$. The restriction $\IIntP{L,U} =  \IInt{L, U} \cap \Omega_1$ %
of an IoM $\IInt{L,U}$ to its probabilities forms a convex \emph{credal set} \cite{leviEnterpriseKnowledgeEssay1980}.\sepfootnote{credal} This set $\IIntP{L,U}$ -- which has previously been studied in the IP literature when $\abs{\T}$ is finite under the name \emph{probability interval} (PI) \cite{decamposPROBABILITYINTERVALSTOOL1994} -- is the fundamental object of analysis throughout this paper.

As a direct consequence of Definition~\ref{def:iom}, the odds $P\left(A\right) / P\left(B\right)$ -- for any $P \in \IInt{L,U} $ and any $A, B \in \mathscr{F}$ -- are bounded between $L\left(A\right)/U\left(B\right)$ and $U\left(A\right)/L\left(B\right)$, whenever these ratios are well-defined. An IoM can also be expressed as a \emph{density bounded class}, which is defined as follows: Fix some $\nu \in \Omega$ and pick $\nu$-densities $l \le u$. The density bounded class $\IInt{l,u}$ consists of $\nu$-densities $f$ satisfying $l \le f \le u$. (This is equivalent to Definition~\ref{def:iom} when $U \in \Omega$ since every $\mu \in \IInt{L,U}$ is absolutely continuous with respect to $U$ and so will always have a $\nu$-density when $\nu = U$. See Proposition~\ref{propDRCIoM} in Appendix~\ref{appendixAdditionalResults} for more details.) Density bounded class, or the closely-related density ratio classes, are often used as prior neighbourhoods in robust Bayesian analysis due to their attractive properties; see e.g. \cite{berger1990robust, lavineApproachRobustBayesian1991, wasserman1992invariance, seidenfeld1993dilation} and especially \cite{wassermanComputingBoundsExpectations1992}. Moreover, IoMs have also been used to represent neighbourhoods of sampling distributions \citep{lavine1991sensitivity}. When used in conjunction with prior neighbourhoods they augment Bayesian robustness beyond prior robustness without resorting to trivial posterior bounds. In fact, a neighbourhood of sampling distributions must have densities bounded away from zero and infinity -- as is the case for a density bounded class $\IInt{l,u}$ (with $0 < l$ and $u < \infty$), but is not so for other popular neighbourhood models -- to ensure that the resulting posterior neighbourhood have non-trivial extrema \citep[Example~1]{lavine1991sensitivity}, a point which is closely connected to the necessity of $\dMult$ for encoding the notion of ``privacy as indistinguishability'' (see Remark~\ref{remarkIndistinguishability}).

Here and elsewhere in this article, the term ``density'' is used in the broad sense of the Radon-Nikodym derivative $\frac{d\mu}{d\nu}$ of a measure $\mu \in \Omega$ with respect to a dominating measure $\nu \in \Omega$. Among other examples, this usage encompasses both probability density functions (PDFs) of continuous real-valued random variables (where the dominating measure is the Lebesgue measure) and probability mass functions (PMFs) of discrete random variables (where the dominating measure is the counting measure).

Theorem~\ref{thm:dp-iom} establishes an equivalence between the $\epsilon$-DP property of a data-release mechanism $M$ and the interval of measures $M$ induces.

\begin{definition}\label{defnConnected}
	Two datasets $x, x' \in \mathcal X$ are \emph{connected} -- or more precisely, $d$-connected -- if $d(x,x') < \infty$. In this case, we say that $x$ is a \emph{connection} of $x'$, and that the probability measures $P_x$ and $P_{x'}$ are \emph{connected}.  More generally, $S \subset \mathcal X$ is connected if all $x, x' \in S$ are. 
 
    The data universe $\mathcal X$ is partitioned into \emph{connected components} $[x] = \{x' \in \mathcal X \mid d(x, x') < \infty\}$. More generally, for $S \subset \X$, define 
    \[[S] = \{ x \in \X \mid \exists x' \in S \mathrm{\ s.t.\ } d(x, x') < \infty\}.\]
\end{definition}

Since the Lipschitz condition \eqref{eqDPDefnLip} is vacuous when $d(x, x') = \infty$, DP only constrains a mechanism $M$ to act similarly on connected datasets $x, x'$; it makes no (explicit) restrictions between outputs $M(x, U)$ and $M(x',U)$ for unconnected $x, x'$. That is, there is no privacy guarantee of indistinguishability between unconnected datasets (although restrictions between outputs for connected $x, x'$ may induce restrictions between outputs for unconnected $x, x'$). 

When $d = \dHam$, any dataset $x, x'$ of different dimension (i.e. $x, x'$ such that $\abs{x} \ne \abs{x'}$) are unconnected. Hence, $\epsilon$-DP with $d = \dHam$ does not protect against, for example, an attacker determining $\abs{x}$. Unconnected datasets also arise in the presence of \emph{invariants} \citep{gongCongenialDifferentialPrivacy2020, bailieRefreshmentStirredNot2024a}.

\begin{assumption}\label{assumptionGraphDistance}
	$d(x,x')$ is equal to the length of a shortest path between $x$ and $x'$ in a graph on $\mathcal X$ with unit-length edges.
\end{assumption}

When $d(x, x') > 1$, the Lipschitz condition \eqref{eqDPDefnLip} is called \emph{group privacy}. This terminology comes from the following intuition: When each $x_i$ represents an individual, condition \eqref{eqDPDefnLip} with $\dHam(x, x') > 1$ (or $\dSym(x, x') > 1$) is protecting multiple individuals' privacy simultaneously. We prove in Theorem~\ref{thm:dp-iom} that Assumption \ref{assumptionGraphDistance} and individual-only privacy (i.e. condition \eqref{eqDPDefnLip} for $x,x'$ with $d(x,x') = 1$) together imply group privacy.

The following lemma is useful for Theorem~\ref{thm:dp-iom} and subsequent discussions.

\begin{lemma}\label{lemmaMultIoM}
    For any $\mu, \nu \in \Omega$ and $\epsilon > 0$, we have $\nu \in \IInt{e^{-\epsilon} \mu, e^{\epsilon} \mu}$ if and only if 
        \[\dMult(\mu, \nu) \le \epsilon.\]
    Hence, for any $\mu, \nu \in \Omega$ and $0 < a \le 1 \le b < \infty$,
    \begin{enumerate}
        \item $\nu \in \IInt{a \mu, b \mu}$ implies $\dMult(\mu, \nu) \le \max (-\ln a, \ln b)$; and
        \item $\dMult(\mu, \nu) \le \min (-\ln a, \ln b)$ implies  $\nu \in \IInt{a \mu, b \mu}$.
    \end{enumerate}
\end{lemma}

The proof of Lemma~\ref{lemmaMultIoM} is given in Appendix~\ref{appendixProofs}, which also contains all other proofs which have been omitted from the main body of this paper.

\begin{theorem}\label{thm:dp-iom}
	Let $M : \mathcal X \times [0,1] \to \mathcal T$ be a data-release mechanism with the seed $U \sim \mathrm{Unif}[0,1]$ inducing a probability $P_x$ on $M(x,U)$ (where $x$ is taken as fixed).
	
	For $0 \le \epsilon < \infty$, the following statements are equivalent given Assumption \ref{assumptionGraphDistance}:
	\begin{enumerate}[label=\Roman*]
		\item\label{mainThmStatement1} $M$ is $\epsilon$-differentially private. 
		\item\label{mainThmStatement2} $P_{x'}(S) \le e^{\epsilon} P_{x} (S)$ for all $S \in \mathscr F$ and all $x, x' \in \mathcal X$ with $d(x,x') = 1$.
		\item\label{mainThmStatement3} For all $\delta \in \mathbb N$ and all $x, x' \in \mathcal X$ with $d(x, x') = \delta$,
		\begin{equation*}
			P_{x'} \in \IIntP{L_{x,\delta\epsilon},U_{x,\delta\epsilon}},
		\end{equation*} 
		where $L_{x,\delta\epsilon}=e^{-\delta\epsilon}P_{x}$ and $U_{x,\delta\epsilon}=e^{\delta\epsilon }P_{x}.$
		\item\label{mainThmStatement4} For all $x \in \mathcal X$ and all measures $\nu \in \Omega$, if $P_x$ has a density $p_x$ with respect to $\nu$, then every $d$-connected $x' \in [x]$ also has a density $p_{x'}$ (with respect to $\nu$) satisfying
		\begin{equation}\label{eqStatement4}
			p_{x'}(t) \in p_x(t) \exp \left[\pm \epsilon d(x,x') \right],
		\end{equation}
		for all $t \in \mathcal T$.
	\end{enumerate}
\end{theorem}

In \eqref{eqStatement4}, the notation $a \in \exp (\pm b)$ is shorthand for 
\[\exp (-b) \le a \le \exp (b).\]

\ref{mainThmStatement2} is the standard definition of pure $\epsilon$-DP \citep{dwork2006calibrating} and is listed here to justify our novel formulation given in Definition~\ref{defnDP}. Without Assumption \ref{assumptionGraphDistance}, group privacy is not implied by \ref{mainThmStatement2}. Hence Assumption \ref{assumptionGraphDistance} is needed only to extend \ref{mainThmStatement2} to provide group privacy; the equivalences between \ref{mainThmStatement1}, \ref{mainThmStatement3} and \ref{mainThmStatement4} are automatic. Without Assumption \ref{assumptionGraphDistance} (which almost always holds in practice, such as for $d = \dHam$ or $\dSym$), our definition of $\epsilon$-DP would be more stringent than the standard formulation.

\begin{proof}(sketch) %
    ``\ref{mainThmStatement1} $\Leftrightarrow$ \ref{mainThmStatement2}'': Since $d$ is a graph distance, there is a path $x = x_0,\ldots, x_n = x'$ such that $d(x,x') = n$ and $d(x_i, x_{i+1}) = 1$. By the triangle inequality, 
	\[\dMult(P_x, P_{x'}) \le \sum_{i=0}^{n-1} \dMult (P_{x_i}, P_{x_{i+1}}).\]
	Hence $\epsilon$-DP is equivalent to the Lipschitz condition \eqref{eqDPDefnLip} holding only when $d(x, x') = 1$. The equivalence between \ref{mainThmStatement1} and \ref{mainThmStatement2} then follows by an application of Lemma~\ref{lemmaMultIoM}: $e^{-\epsilon} P_x(S) \le P_{x'}(S) \le e^\epsilon P_{x}(S)$ for all $S \in \mathscr F$ if and only if $\dMult(P_x, P_{x'}) \le \epsilon$. ``\ref{mainThmStatement1} $\Leftrightarrow$ \ref{mainThmStatement3}'' is immediate by Lemma~\ref{lemmaMultIoM}.%
	
	``\ref{mainThmStatement3} $\Leftrightarrow$ \ref{mainThmStatement4}'': The direction $\Rightarrow$ is straightforward since the densities in an interval of measure $\IInt{L,U}$ are bounded by the densities of $L$ and $U$. In the other direction $\Leftarrow$, $P_x$ is always absolutely continuous with respect to itself, hence taking $P_x$ to be the dominating measure $\nu$, we have that \eqref{eqStatement4} implies $P_{x'} \in \IIntP{L_{x, \delta\epsilon},U_{x, \delta\epsilon}}$.
\end{proof}

\begin{remark}\label{remarkDPCredalSet}
    A \emph{distorting function} $\dist$ can be thought of as a generalised notion of distance between two probabilities or measures. (See Appendix~\ref{appDistort} for a precise definition.) %
    Given a distorting function $\dist$ and a \emph{distortion parameter} $r > 0$, the \emph{distortion model} on a probability $P$ is the closed ball $B_{\dist}^r(P) = \{Q \in \Omega_1 \mid \dist(Q, P) \le r\}$ centred at $P$ with radius $r$ \cite{montesUnifyingNeighbourhoodDistortion2020}.

    Given a probability $P$, the symmetric probability interval $\IIntP{e^{-\epsilon} P, e^{\epsilon} P}$ is a \emph{distortion model} because, by Lemma~\ref{lemmaMultIoM}, $\IInt{e^{-\epsilon} P, e^{\epsilon} P}$ is the closed $\epsilon$-multiplicative-distance-ball:
    \[\IInt{e^{-\epsilon} P, e^{\epsilon} P} = B_{\dMult}^\epsilon (P) = \{ \mu \in \Omega \mid \dMult(\mu, P) \le \epsilon \}.\]
    Since the multiplicative distance $\dMult$ is used in defining the ball $B_{\dMult}^\epsilon (P)$, it serves the role of the \emph{distorting function} of the distortion model $\IIntP{e^{-\epsilon} P, e^{\epsilon} P}$; and the radius $\epsilon$ is the \emph{distortion parameter}.%
    
    It is straightforward to verify that the credal set $\IIntP{L,U}$ is convex. In Appendix~\ref{appDistort}, we also prove that it is closed with respect to the supremum norm, when restricting to the setting of \cite{montesUnifyingNeighbourhoodDistortion2020}. That is, $\IIntP {L, U} \cap \Omega^*_1$ 
    is a closed subset of $\Omega^*_1 = \{ P \in \Omega_1 \mid P(\{t\}) > 0 \ \forall t \in \T \}$ (where closure is respect to the supremum norm), provided that the space $\T$ has finite cardinality.

    In general, a PI $\IIntP{L, U}$ does not satisfy the definition of a distortion model. %
    In fact, $\IIntP{L, U}$ is a distortion model only when $\IInt{L, U}$ is symmetric in the sense that the lower and upper measures are equidistant from a central (`nucleus') probability $P$. This is the case for the PI $\IIntP{e^{-\epsilon} P_x, e^{\epsilon} P_x}$ induced by an $\epsilon$-DP mechanism $M$. 
    
    As proven in Theorem~\ref{thm:dp-iom}, pure $\epsilon$-DP is the requirement that $P_x'$ lies in the neighbourhood $B_{\dMult}^{\delta\epsilon}(P_x)$, where $\delta = d(x, x')$. Many of the variants of $\epsilon$-DP -- such as $(\epsilon,\delta)$-DP \citep{dworkOurDataOurselves2006}, zero-concentrated DP \citep{dworkConcentratedDifferentialPrivacy2016, bunConcentratedDifferentialPrivacy2016a} and R\'enyi DP \citep{mironovRenyiDifferentialPrivacy2017} -- replace the multiplicative distance $\dMult$ with another distorting function $\dist$. Consequently, these variants can also be characterised as distortion models: each of them is the requirement that $P_{x'}$ lies in the neighbourhood $B_{\dist}^{\delta \epsilon} (P_x)$, for the appropriate choice of distorting function $\dist$. (See \cite{bailieRefreshmentStirredNot2024a} for the choices of $\dist$ corresponding to $(\epsilon,\delta)$-DP, zero-concentrated DP and R\'enyi DP, and see \cite{pejoGuideDifferentialPrivacy2022} for a catalogue of variants of DP.)    
\end{remark}

\begin{remark}\label{remarkIndistinguishability}
\ref{mainThmStatement4} of Theorem~\ref{thm:dp-iom} is a strong property. It provides an quantification of the ``indistinguishability'' between data $x, x' \in \mathcal X$: if $x, x'$ have densities $p_x, p_{x'}$ satisfying \eqref{eqStatement4}, then they are indistinguishable at the level $\epsilon$. (Equation~\eqref{eqStatement4} is termed $\epsilon$-indistinguishability in the literature, see e.g. \cite{dwork2006calibrating, dwork2014algorithmic, vadhanComplexityDifferentialPrivacy2017}).
More fundamentally, \ref{mainThmStatement4} provides a categorical notion of indistinguishability:
It implies that, for an $\epsilon$-DP mechanism, all connected $P_x$ are mutually absolutely continuous. Further, for all connected $x,x' \in \mathcal X$ and all $t \in \T$, either $p_x(t)$ and $p_{x'}(t)$ are both zero or both non-zero. In intuitive terms, this means that if any $x$ is plausible after observing $T = t$ (i.e. $p_x(t) > 0$) then all its connections $x' \in [x]$ are also plausible. This is a strong notion of privacy: regardless of the output $T = M(x, U)$, it's impossible for an attacker to distinguish between connected $x, x'$ with certainty. In other words, the fiducial distribution for $x$ is never degenerate (assuming that every $x$ has at least one connection).

This notion of privacy is the motivation for $\dMult$ in place of more standard concepts in the robustness literature such as total variation distance or $\epsilon$-contamination classes. Indeed, this categorical notion of indistinguishability requires that $p_x(t)/p_{x'}(t)$ is bounded away from zero and infinity, which is equivalent to $P_{x'} \in \IInt{a P_{x}, b P_{x}}$ for some $0 < a \le 1 \le b < \infty$. Yet Lemma~\ref{lemmaMultIoM} states that $P_{x'} \in \IInt{a P_{x}, b P_{x}}$ only if \label{review5} $\dMult(P_x, P_{x'}) \le \max (-\ln a, \ln b)$. Therefore, using the multiplicative distance $\dMult$ is necessary to encode the idea of privacy as indistinguishability between connected $x, x'$.  

This argument demonstrates that the Lipschitz condition \eqref{eqDPDefnLip} with another distorting function $\dist$ in place of $\dMult$ will not ensure indistinguishability (except in the trivial case where $\alpha \dist \ge \dMult$ for some constant $\alpha$). This is why common variants of pure $\epsilon$-DP -- such as $(\epsilon,\delta)$-DP \citep{dworkOurDataOurselves2006}, zero-concentrated DP \citep{dworkConcentratedDifferentialPrivacy2016, bunConcentratedDifferentialPrivacy2016a} and R\'enyi DP \citep{mironovRenyiDifferentialPrivacy2017} (which, as described in Remark~\ref{remarkDPCredalSet}, all replace $\dMult$ with another distorting function $\dist$)
-- do not guarantee this strong notion of privacy, even though they may be preferred over pure $\epsilon$-DP for data utility reasons.
\end{remark}

The observations of Theorem~\ref{thm:dp-iom}, specifically the equivalent characterization of $\epsilon$-DP via intervals of measures established by~\ref{mainThmStatement3} and~\ref{mainThmStatement4}, bear important consequences for statistical inference from privacy-protected data. Notably, they impose meaningful bounds on both the probability of the privatised query and on relevant quantities in the frequentist and Bayesian inference from the privatised queries. These bounds are valid under arbitrary statistical models for the unknown confidential database, assuming only mild conditions on the models' support. The next three sections explore these consequences in detail.

\section{Bounds on the Privatised Data Probability}\label{sec:marginal}

Consider the situation of statistical inference, where a data analyst supplies a parametric model $\mathcal P = \{P_\theta \mid \theta \in \Theta\}$ of data-generating distributions $P_\theta$. Nature generates data $X \sim P_{\theta}$ according to some unknown $\theta \in \Theta$. (We use capital $X$ to emphasise that the dataset is now random, whereas in the previous sections, it was considered fixed.) In the typical non-private setting, the data analyst observes $X$ directly. In the private setting, the data analyst only sees the summary statistic $T = M(X,U) \sim P_X$ outputted from a privacy-preserving data-release mechanism $M$. (We now require that the data universe $\mathcal X$ is equipped with a $\sigma$-algebra $\mathscr G$ and that every data-release mechanism $M$ is $(\mathscr G \otimes \mathcal B[0,1], \mathscr F)$-measurable, where $\mathcal B[0,1]$ is the Borel $\sigma$-algebra on $[0,1]$.)

The relevant vehicle for inference in the private setting is the marginal probability of the observed data $T$: 
\begin{align}
P( T \in S \mid \theta) %
&= \int_{\mathcal X} P_x(S) dP_{\theta}(x). \label{eq:marginal-data-probability}
\end{align}
We call $P( T \in S \mid \theta)$ the \emph{privatised data probability}. (Proposition~\ref{propDefnPrivDataProb} proves that it is well-defined.) Viewed as a function of $\theta$, $P( T \in S \mid \theta)$ is the \emph{marginal likelihood} of $\theta$. When the data observed by the analyst is privacy-protected, all frequentist procedures compliant with likelihood theory and all Bayesian inference hinge on this function \cite{bergerLikelihoodPrinciple1988}. The crucial role of~\eqref{eq:marginal-data-probability} for inference from privacy-protected data  was first recognized in the differential privacy literature by \citet{williams2010probabilistic}, and has since been utilized extensively to derive likelihood and Bayesian methodologies \cite[e.g.][]{Awan2018:Binomial,awan2020differentially,bernstein2018differentially,bernstein2019differentially,gong2022exact,ju2022data}.

When $M$ is $\epsilon$-DP and the support $\supp(P_\theta)$ of $P_\theta$ is $d$-connected, the existence of a density $p(t \mid \theta)$ for $P(T \in S \mid \theta)$ is implied by Theorem~\ref{thm:dp-iom}. The following result proves this density always exists -- as long as one restricts to a subspace of $\T$ and assumes that (informally) the support of ``$P(x \mid t_0, \theta)$'' is $d$-connected for some given $t_0 \in \T$. Other than this weak assumption, the following results hold for arbitrary data-generating models $\{P_\theta \mid \theta \in \Theta\}$ and $\epsilon$-DP mechanisms $M$.

To state this assumption more precisely, define $\supp(x \mid t, \theta)$ as the set of databases $x \in \mathcal X$ which could both generate $t$ and be generated by $P_\theta$. That is, $\supp(x \mid t, \theta)$ is informally the intersection of $\supp (P_\theta) \approx \{x \mid p_\theta(x) > 0\}$ and $\{x \mid p_x(t) > 0\} \approx \{x \mid t \in \supp (P_x)\}$. See Appendix \ref{appDefnSupport} for an exact definition.

\begin{theorem}\label{thmMarginalData}
    Let $M$ be an $\epsilon$-DP mechanism. Fix some $t_0 \in \T$ and suppose that $\supp(x \mid t_0, \theta)$ is $d$-connected. Define $\T_0 = \{t \in \T \mid \supp (x \mid t, \theta) \subset \supp (x \mid t_0, \theta)\}$. Then, there exist measures $L_{\theta, \epsilon}$ and $U_{\theta,\epsilon}$ on $(\T, \F)$ with densities $l_{\theta,\epsilon}$ and $u_{\theta,\epsilon}$ satisfying
    \begin{equation*}%
        l_{\theta,\epsilon}(t) = \esssup_{x_* \in \supp (x \mid t_0, \theta)} \exp \left( - \epsilon d_{*} \right) p_{x_*}(t), \quad \text{and}\quad u_{\theta,\epsilon}(t)=\essinf_{x_* \in \supp (x \mid t_0, \theta)} \exp \left( \epsilon d_{*} \right) p_{x_*}(t),
    \end{equation*}
    for all $t \in \T_0$, where $d_* = \sup_{x \in \supp (x \mid t_0, \theta)} d(x, x_*)$. 
    
    Furthermore, the privatised data probability $P(T \in \cdot \mid \theta)$ is bounded by $L_{\theta,\epsilon}$ and $U_{\theta,\epsilon}$ on $\T_0$:
    \begin{equation}\label{eqThmMarginalIoM}
        P(T \in \cdot \cap \T_0 \mid \theta) \in \IInt{L_{\theta,\epsilon}, U_{\theta,\epsilon}}.
    \end{equation}
\end{theorem}

\begin{proof}(sketch) 
    The existence of a density $p(t \mid \theta)$ for $P( T \in \cdot \cap \T_0 \mid \theta)$ follows from the fact that all $P_x$ with $x \in \supp (x \mid t_0, \theta)$ are mutually absolutely continuous by Theorem~\ref{thm:dp-iom}. For the upper bound of~\eqref{eqThmMarginalIoM}, first observe that
	\begin{align*}
		p(t \mid \theta) &= \int_{\supp (x \mid t_0, \theta)} p_x(t) dP_\theta(x) \\
		&\le \int_{\supp (x \mid t_0, \theta)} e^{\epsilon d(x, x_*)} p_{x_*}(t) d P_\theta (x) \\
		&\le e^{\epsilon d_{*}} p_{x_*}(t).
	\end{align*}
    Since the above inequalities hold for all $x_* \in \supp(x \mid t_0, \theta)$, we can take the essential infimum over $x_*$ to obtain the bound $p(t \mid \theta) \le u_{\theta,\epsilon}(t)$. %
    The lower bound of~\eqref{eqThmMarginalIoM} follows similarly. 
\end{proof}

It becomes apparent in the proof of Theorem~\ref{thmMarginalData} that this result can be generalised in the following way:
In defining $\T_0 = \{t \in \T \mid \supp (x \mid t, \theta) \subset \supp (x \mid t_0, \theta)\}$, one may replace $\supp (x \mid t_0, \theta)$ with any measurable $S$ satisfying
\begin{equation*}%
    \supp (x \mid t_0, \theta) \subset S \subset [\supp (x \mid t_0, \theta)].
\end{equation*}
(The notation $[\cdot]$ is defined in Definition~\ref{defnConnected}.) Theorem~\ref{thmMarginalData} holds with this new $\T_0$, provided that $\supp (x \mid t_0, \theta)$ is replaced by $S$ in the definitions of $l_{\theta,\epsilon}, u_{\theta,\epsilon}$ and $d_*$. This demonstrates that the density $p(t \mid \theta)$ exists on a larger $\T_0$, although the resulting bounds $l_{\theta,\epsilon}$ and $u_{\theta,\epsilon}$ on $p(t \mid \theta)$ may be wider.

Theorem~\ref{thmMarginalData} shows that the privatised data probability $P(T \in \cdot \mid \theta)$ is in a probability interval, and that this probability interval is bounded by $L_{\theta, \epsilon}$ and $U_{\theta, \epsilon}$ on $\T_0$. Broadly speaking, this theorem has two uses. Firstly, $t_0$ may be taken to be the realised value of $T$. Then Theorem~\ref{thmMarginalData} can be interpreted as bounds on the marginal likelihood $l(\theta \mid t_0)$ of $\theta$. (For this application, one must make the additional assumption that $\bigcup_{\theta \in \Theta} \supp(x \mid t_0, \theta)$ is $d$-connected, so that the densities $p(t_0 \mid \theta)$ of the privatised data probability, across the different values of $\theta$, share a common dominating measure. This ensures that the likelihood $l(\theta \mid t_0) = p(t_0 \mid \theta)$ is well-defined as a function of $\theta$.) Secondly, one may be interested in understanding the privatised data probability within some subspace $S \subset \T$. Although it may not always be possible, if one can find some $t_0 \in \T$ such that $S \subset \T_0$, then Theorem~\ref{thmMarginalData} provides information on what the privatised data probability looks like within the subspace of interest $S$. 

Surprisingly, the interval of measures $\IInt{L_{\theta, \epsilon}, U_{\theta, \epsilon}}$ in \eqref{eqThmMarginalIoM} depends on the data-generating distribution $P_\theta$ only through $\supp (x \mid t_0, \theta)$. When $\supp (P_\theta)$ is constant, $\IInt{L_{\theta, \epsilon}, U_{\theta, \epsilon}}$ is completely free of $\theta$. Alternatively, one may take the essential-infimum of $L_{\theta, \epsilon}$ over $\theta \in \Theta$ to obtain a bound on $P(T \in \cdot \cap \T_0 \mid \theta)$ which is completely free of $\theta$, although it is likely such a bound will be vacuous.  %

Theorem~\ref{thmMarginalData} is only practically meaningful when $d_* < \infty$. Typically $\sup_{x, x' \in \X} d(x, x') = \infty$, which one might presume would imply that $d_* = \infty$. But $\supp(P_\theta)$ can be much smaller than the data universe $\mathcal X$ when the analyst has prior knowledge of the data $X$. The analyst is free to restrict $\supp (P_\theta)$ to the set of datasets they deem plausible; the tighter this restriction, the stronger Theorem~\ref{thmMarginalData} is. For example, the analyst may have an upper bound $b$ on the number of records $\abs{X}$. (Provided that $\supp (x \mid t_0, \theta)$ is connected -- a weak assumption, as explained\label{review2.4} in Remark~\ref{remarkConnect} -- this would imply $d_{*} \le b$ for typical choices of $d$.) Moreover, $\supp (x \mid t_0, \theta)$ can be much smaller than $\mathcal X$ when $t_0$ restricts the possible values of $X$, such as in the presence of invariants \citep{gongCongenialDifferentialPrivacy2020, bailieRefreshmentStirredNot2024a}. For example in local DP (or bounded DP more generally), the number of records $\abs{t}$ is invariant; this restricts $\supp (x \mid t_0, \theta)$ to data $x$ satisfying $\abs{x} = \abs{t_0}$, which would typically imply $d_* \le \abs{t_0}$.

\begin{remark}\label{remarkConnect}
    Theorem~\ref{thmMarginalData} only relies on a single assumption which concerns the connectedness of $\supp (x \mid t_0, \theta)$. This assumption is weak. In fact, we can always augment the data-release mechanism $M$ so that this assumption is satisfied without increasing $M$'s privacy loss $\epsilon$. Specifically, the (deterministic) mechanism $x \mapsto [x]$ (which publishes the connected component $[x]$ of the observed data $x$) is trivially $\epsilon$-DP with $\epsilon = 0$. Publishing $[x]$ alongside $M(x,U)$ ensures that $\supp (x \mid t, \theta)$ is always connected, for all $t$ and $\theta$. (This argument is formalised in Proposition~\ref{propMakeConnected}.)
\end{remark}

We illustrate Theorem~\ref{thmMarginalData} with two examples.

\begin{example}[privatised binary sum]\label{ex:binary-sum}
Suppose the database $x \in \mathcal X = \bigcup_{n=1}^\infty \{0, 1\}^n$ consists of $n$ records of binary features, and its sum  $q(x) = \sum_{i=1}^{n} x_i$ is to be queried. Consider sanitising $q(x)$ using the Laplace mechanism defined in Example~\ref{ex:laplace}. For every privacy loss $\epsilon > 0$ and every database $x$,  
\begin{equation*}
     p_{x}(t) = \frac{\epsilon}{2\Delta(q){}} \exp\left(\frac{\epsilon\lvert t - q(x) \rvert}{\Delta(q)}\right),
\end{equation*}
where, in this case, the global $\ell_1$-sensitivity $\Delta(q)$ (defined in \eqref{eqDefnGlobalSensitivity}) is one. 

The data analyst posits an arbitrary statistical model $X \sim P_\theta$ for $\theta \in \Theta$ with $\supp (P_\theta) \subset \{ x \in \X \mid \abs{x} \le 10\}$, and considers the confidential and unknown database $x$ to be a realization from this model. Since $\supp (P_x) = \mathbb R$ for all $x \in \X$, the assumption of Theorem~\ref{thmMarginalData} simplifies to the requirement that $\supp (P_\theta)$ is $d$-connected. Moreover, $\T_0 = \T = \mathbb R$. (Both of these points hold regardless of the choice of $t_0$.)

Figure~\ref{fig:marginal} displays the lower and upper densities, $l_\epsilon = \essinf_{\theta \in \Theta} l_{\theta, \epsilon}$ and $u_\epsilon = \esssup_{\theta \in \Theta} u_{\theta, \epsilon}$, for the privatised data probability $p(t \mid \theta)$. The analyst upper bounds the number of records $\abs{x}$ by $10$, so that $d_* = 10$. The left and right panels display bounds under two different settings of $\epsilon$. The bounds are tighter and more informative when privacy protection is more stringent ($\epsilon = 0.1$),  and looser as the privacy loss increases ($\epsilon = 0.25$). Notice that these bounds for  $p(t \mid \theta)$ are functions of the value of the privatised query $t$. In particular, they do not depend on $\theta$ nor the form of the posited data model $P_\theta$.
\end{example}

\begin{figure}
    \floatconts{fig:marginal}
    {\caption{\emph{Upper and lower bounds for the density $p(t \mid \theta)$ of the privatised binary sum (Example~\ref{ex:binary-sum})}. The privacy loss is $\epsilon = 0.1$ (left) and $\epsilon = 0.25$ (right). These bounds depend on the assumed data model $P_\theta$ only through the analyst's belief that the number of records $|x|$ is bounded by ten, which means that $\supp(P_\theta) \subset \{x \in \X \mid |x| \le 10\}$. They are tighter and more informative when the privacy protection is more stringent (i.e. when $\epsilon$ is smaller).}}
    {
    \includegraphics[width=0.49\textwidth]{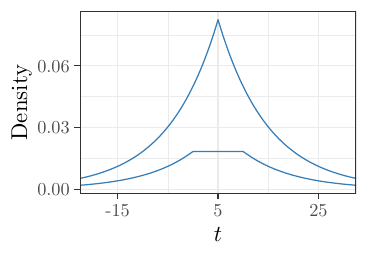}
    \includegraphics[width=0.49\textwidth]{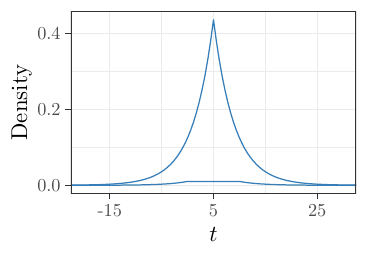}
	}
\end{figure}

\begin{example}[local $\epsilon$-DP]\label{ex:local-dp}
    Suppose the distribution $P_x$ of the published summary statistic $T$ factors as $\prod_{i=1}^n P_{x_i}$, where $n = \abs{x}$. (This always holds under the local, non-interactive model of DP, as we described in Example~\ref{exampleRandomisedResponse}.) Then $\abs{T} = \abs{X}$ and hence $\supp (P_x) \subset \{ t \in \T \mid \abs{t} = \abs{x}\}$. 
    
    Most local DP mechanisms satisfy the stricter assumption that $\supp (P_x) = \{ t \in \T \mid \abs{t} = \abs{x}\}$. Under this assumption, $\T_0 = \{t \in \T \mid \abs{t} = \abs{t_0}\}$ and, 
    if $d = \dHam$ (as is typical for local DP), then $d_* \le \abs{t_0}$ regardless of the choice of $x_*$. Hence, by Lemma~\ref{lemmaMarginalData} of Appendix~\ref{appendixProofs} (which is used in proving Theorem~\ref{thmMarginalData}), the density of an $\epsilon$-DP mechanism is bounded by
    \begin{equation*}
		p ( t \mid \theta) \in \prod_{i=1}^n p_{x_{i}}(t_i) \exp \left( \pm \epsilon n \right),
	\end{equation*}
    for any $x$ and any $t$ with $\abs{x} = \abs{t} = n$. Applying this result to the randomised response mechanism (Example~\ref{exampleRandomisedResponse}), $\min_{x_i} p_{x_i} (t_i) = (\exp \epsilon + 1)^{-1}$ and $\max_{x_i} p_{x_i} (t_i) = e^\epsilon (\exp \epsilon + 1)^{-1}$, so that 
    \begin{equation}\label{eq:marginal_rr}
    \frac{1}{(\exp \epsilon + 1)^{\abs{t}}} \le p(t \mid \theta) \le \frac{\exp (\abs{t} \epsilon)}{(\exp \epsilon + 1)^{\abs{t}}}.
    \end{equation}
    The bounds in~\eqref{eq:marginal_rr} depend on $t$ only through $\abs{t}$ (the number of records), regardless of the records' values. Figure~\ref{fig:marginal_rr} displays these bounds as a function of $\abs{t}$ for $\epsilon = 1$. As more records are released (larger $\abs{t}$), both bounds tend to zero with a narrowing gap.
\end{example}

\begin{figure}
	\floatconts{fig:marginal_rr}
	{\caption{\emph{Upper and lower density bounds for $p(t \mid \theta)$ under randomised response (Example~\ref{ex:local-dp})}. The privacy loss is $\epsilon = 1$. These bounds are a function of $t$ only through $\lvert t\rvert$ (the number of observed records).}}
	{\includegraphics[width=\textwidth]{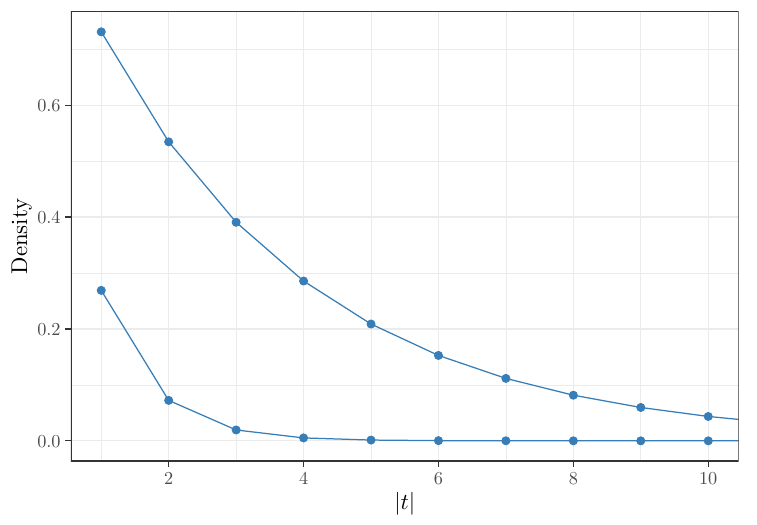}}
\end{figure}

\section{Frequentist Privacy-Protected Inference}
\label{sec:frequentist}

The interval of measures formulation of $\epsilon$-DP also shows that Neyman-Pearson hypothesis testing is restricted in the private setting, as demonstrated by the following theorem.

\begin{theorem}\label{thmDPHypothesisTesting}
    Consider testing $H_0 : \theta = \theta_0$ versus $H_1 : \theta = \theta_1$ for some $\theta_0 \ne \theta_1 \in \Theta$. Let $S_i = \supp(P_{\theta_i})$ and suppose that $S_0 \cup S_1$ is $d$-connected. In the private setting where the observed data $T$ is the output of an $\epsilon$-DP mechanism, the power of any level-$\alpha$ test is bounded above by $\alpha \exp (d_{**} \epsilon)$ where 
    \[d_{**} = \sup_{x \in S_0,x' \in S_1} d(x,x').\]
\end{theorem}

\begin{proof}(sketch)
	By \ref{mainThmStatement4} of Theorem~\ref{thm:dp-iom}, 
	\begin{align*}
		\frac{p(t \mid \theta_1)}{p(t \mid \theta_0)} &= \frac{\int_{S_1} p_x(t) dP_{\theta_1}(x)}{\int_{S_0} p_{x'}(t) dP_{\theta_0}(x')} \\
		&= \int_{S_1} \left[ \int_{S_0} \frac{p_{x'}(t)}{p_x(t)} d P_{\theta_0}(x') \right]^{-1} dP_{\theta_1}(x) \\
		&\in \exp \left( \pm \epsilon d_{**} \right).
	\end{align*}
	Let $R$ be the rejection region of a test with size $P(T \in R \mid \theta_0) \le \alpha$ and let $\nu$ be the dominating measure of the densities $p(t \mid \theta_0)$ and $p(t \mid \theta_1)$. Then
	\begin{align}
		P(T \in R \mid \theta_1) &= \int_R p(t \mid \theta_1) d\nu(t) \label{eqProofThmDPHypothesisTesting}\\
		&\le \exp ( d_{**} \epsilon) \int_R p(t \mid \theta_0) d\nu(t) \notag\\
		& \le \alpha \exp ( d_{**}\epsilon).\notag %
	\end{align}
\end{proof}

Compare Theorem~\ref{thmDPHypothesisTesting} to the hypothesis test $H_0 : x_{1:m} = y$ versus $H_1 : x_{1:m} = y'$ where $m \le \abs{x}$. (Here we assume that the datasets $x \in \mathcal X$ are vectors of length $\abs{x}$ and consist of records $x_1, x_2, \ldots, x_n$, where $n = \abs{x}$. For $1 \le i \le j \le \abs{x}$, the notation $x_{i:j}$ denotes the sub-vector $(x_i, x_{i+1}, \ldots, x_{j})$, consisting of the $i$-th through $j$-th records of $x$.) This test models an attacker trying to distinguish the first $m$ records of the database. \citet{wasserman2010statistical} showed that any level-$\alpha$ test of $x_{1:m}$ has power at most $\alpha \exp (\epsilon m )$ when the records $X_i$ are i.i.d. 

If the data analyst restricts $S_0$ and $S_1$ to datasets of length $m$, then typically $d_{**} = m$. Thus, any level-$\alpha$ test on the parameter $\theta$ has the same bound $\alpha \exp (\epsilon m)$ on its power (under an arbitrary data-generating model, not just i.i.d. $X_i$).

\sepfootnotecontent{Thm11Nuisance}{This ignores one minor technicality: the attacker may take some records as nuisance parameters, which they do not want to test. It is straightforward to generalise Theorem~\ref{thmDPHypothesisTesting} to this situation. Without loss of generality, suppose $x_{m+1:n}$ are nuisance parameters when testing $x_{1:m}$ against $x_{1:m}'$. By assigning a conditional probability on $x_{m+1:n}$ satisfying $\pi(x_{m+1:n} \mid x_{1:m}) = \pi(x_{m+1:n} \mid x_{1:m}')$, the nuisance parameters can be integrated out in \eqref{eqProofThmDPHypothesisTesting}.
This gives the same power bound $\alpha \exp (d_{**} \epsilon)$, except now with
\[d_{**} = \sup_{x_{m+1:n}} d \Big(\big[x_{1:m}, x_{m+1:n}\big], \left[x_{1:m}', x_{m+1:n}\right]\Big),\]
which is typically equal to $m$, as before.}

Theorem~\ref{thmDPHypothesisTesting} strictly generalises the result of \citet{wasserman2010statistical}. By taking $\Theta \subset \mathcal X$ and setting $P_\theta$ as degenerate point masses, we recover the set-up of an attacker's hypothesis test.\sepfootnote{Thm11Nuisance} Thus, Theorem~\ref{thmDPHypothesisTesting} is applicable to both the attacker testing $x$ (as in \cite{wasserman2010statistical}) and the analyst testing $\theta$ (with non-degenerate $P_\theta$). This highlights the fundamental tension between data privacy and data utility: bounding an attacker's power will bound the power of a legitimate analyst. However, another look at Theorem~\ref{thmDPHypothesisTesting} seems to suggest a possible way to partially resolve this tension under certain circumstances. The data custodian might have the liberty to choose a metric $d$ on $\X$ that ensures the connectedness assumption of Theorem~\ref{thmDPHypothesisTesting} holds for the hypothesis tests of the typical attacker but not for those of the legitimate analyst. If this happens, the hypothesis test of the attacker -- but not of the legitimate analyst -- will be constrained by Theorem~\ref{thmDPHypothesisTesting}. Such a choice for $d$ would therefore resolve the tension between privacy and utility as it appears in Theorem~\ref{thmDPHypothesisTesting}. (This is not to suggest that the analyst will be totally unaffected by such privacy protection -- any noise infusion can in general decrease the power of their test -- but, at least, they will not be affected to the extent suggested by Theorem~\ref{thmDPHypothesisTesting}.)%

\begin{corollary}\label{corDPHT}
    Under the set-up of Theorem~\ref{thmDPHypothesisTesting}, the power $1-\beta$ of any size-$\alpha$ test is bounded by the inequalities:
    \begin{equation}\label{eqDPHTCor}
        \max \big( \alpha e^{- d_{**} \epsilon}, 1 - e^{d_{**} \epsilon}[1-\alpha] \big) \le 1 - \beta \le \min \big( \alpha e^{d_{**} \epsilon}, 1 - e^{- d_{**} \epsilon} [1 - \alpha] \big).
    \end{equation}
\end{corollary}

\begin{proof}\cite[Section~6.1]{kifer2022bayesian} Let $R$ be the rejection region of a test with size $\alpha = P(T \in R \mid \theta_0)$ and power $1 - \beta = P(T \in R \mid \theta_1)$. In the proof of Theorem~\ref{thmDPHypothesisTesting}, we showed that
\[P(T \in R \mid \theta_1) \le \exp (d_{**} \epsilon) P(T \in R \mid \theta_0).\]
In the same way, one can show that 
\[P(T \in R \mid \theta_1) \ge \exp (-d_{**} \epsilon) P(T \in R \mid \theta_0),\]
and that
\[\exp (-d_{**} \epsilon) P(T \notin R \mid \theta_0) \le P(T \notin R \mid \theta_1) \le \exp (d_{**} \epsilon) P(T \notin R \mid \theta_0).\]
Combining these four inequalities gives~\eqref{eqDPHTCor}.
\end{proof}

\section{Bayesian Privacy-Protected Inference}\label{sec:Bayesian}

Following the set-up from the previous two sections, we further assume that the analyst is Bayesian and places a (proper) prior $\pi$ on $\Theta$. This setting can be seen as a Bayesian hierarchical model where the raw, confidential data $X$ acts as latent parameter in the Markov chain $\theta \to X \to T$. %

We make the following assumption throughout this section.

\begin{assumption}\label{assumptionConnected}
    Define  
    \[\supp(x \mid t) := \bigcup_{\theta \in \supp (\pi)} \supp (x \mid t, \theta).\]
    Fix some $t_0 \in \T$. Suppose that (A) $\supp(x \mid t_0)$ is $d$-connected. Further, assume that (B) the prior $\pi$ on $\theta$ is proper.
\end{assumption}

By the same reasoning as in Remark~\ref{remarkConnect}, Assumption~\ref{assumptionConnected}(A) is weak because it can always be satisfied by augmenting the data-release mechanism $M$ without additional privacy loss.

Theorem~\ref{thmPriorPredictive} establishes bounds on the analyst's prior predictive distribution $P(T \in S) = \iint P_x(S) dP_{\theta}(x) d\pi(\theta)$ for the privatised data $T$. 

\begin{theorem}\label{thmPriorPredictive} 
	Let $M$ be an $\epsilon$-DP mechanism. Define $\T_0 = \{t \in \T \mid \supp (x \mid t) \subset \supp (x \mid t_0)\}$. Then, there exist measures $L_\epsilon$ and $U_\epsilon$ on $(\T, \F)$ with densities $l_\epsilon$ and $u_\epsilon$ satisfying
     \begin{equation*}
        l_{\epsilon}(t) = \esssup_{x_* \in \supp (x \mid t_0)} \exp \left( - \epsilon d_{*} \right) p_{x_*}(t) \quad \text{and} \quad u_{\epsilon}(t) = \essinf_{x_* \in \supp (x \mid t_0)} \exp \left( \epsilon d_{*} \right) p_{x_*}(t),
	\end{equation*}
    for all $t \in \T_0$, where $d_* = \sup_{x \in \supp (x \mid t_0)} d(x, x_*)$. 
        
    Furthermore, the Bayesian analyst's prior predictive probability $P(T \in \cdot)$ is bounded by $L_\epsilon$ and $U_\epsilon$ on $\T_0$:
    \begin{equation}\label{eqthmPriorPredictive}
        P(T \in \cdot \cap \T_0) \in \IInt{L_\epsilon, U_{\epsilon}}.
    \end{equation}
\end{theorem}

\begin{proof}(sketch) Since $p(t) = \int_{\Theta} p(t \mid \theta) d\pi(\theta)$, Theorem~\ref{thmPriorPredictive} follows by showing $p(t \mid \theta)$ is bounded by $l_{\epsilon}\left(t\right)$ and $u_{\epsilon}\left(t\right)$ for almost all $t \in \T_0$. The proof of this is analogous to \eqref{eqThmMarginalIoM}.
\end{proof}

As for Theorem~\ref{thmMarginalData}, one can replace $\supp(x \mid t_0)$ in the definition of $\T_0$ and in the statement of Theorem~\ref{thmPriorPredictive} with any measurable $S \subset \X$ which satisfies
\[ \supp(x \mid t_0) \subset S \subset [\supp (x \mid t_0)].\]
In this way, one can obtain bounds $l_{\epsilon}(t)$ and $u_{\epsilon}(t)$ on the prior predictive density $p(t)$ which apply for a larger subspace $\T_0$, although these bounds will be wider.

The prior predictive distribution $p\left(t\right)$ plays an important role in Bayesian inference and model checking. Before observing the data, $p\left(t\right)$ captures the analyst's implied specification on the data-generation process. After observing the data, this quantity assessed at their value is called \emph{model evidence} where low $p(t)$ reveals potential \emph{conflict} between the data and the prior \cite{evans2006checking,walter2009imprecision}. In addition, it is also the normalizing constant for the posterior distribution $\pi(\theta\mid t)$ and hence is useful for computation.

Theorem~\ref{thmPriorPredictive} shows that the prior predictive distribution $P(T \in \cdot)$ is in a probability interval, and this probability interval is bounded by $L_\epsilon$ and $U_\epsilon$ on $\T_0$.

As an illustration, we can see from Figure~\ref{fig:marginal} of Example~\ref{ex:binary-sum} that when $\epsilon=0.1$, the prior predictive probability of the privatised query is lower-bounded at $\approx 0.02$ whenever $0\le t\le 10$, and can never exceed $\approx 0.08$ even when $t = 5$. On the other hand, when privacy protection is less stringent ($\epsilon=0.5$), the upper bound on the prior predictive probability increases to more than $0.4$.

An important observation on Theorem~\ref{thmPriorPredictive} is the following: While ${p}\left(t\right)$ is a function of both the data model $P_\theta$ and the prior $\pi$, the density bounds $l_{\epsilon}\left(t\right)$ and $u_{\epsilon}\left(t\right)$ are free of both. In this sense, these bounds provide a non-trivial yet almost assumption-free prior predictive model sensitivity analysis. Non-trivial bounds on $p(t)$ are not possible in general; in this case they are a consequence of the data $T$ being $\epsilon$-DP.

Theorem~\ref{thmPosterior} provides general bounds limiting the learning of a Bayesian analyst.

\begin{theorem}\label{thmPosterior}
    Suppose that an $\epsilon$-DP mechanism $M$ outputs the realisation $t_0$. The analyst's posterior probability given $t_0$ satisfies
	\begin{equation}\label{eq:posterior}
         \pi (\theta \in S \mid t_0) \in \pi (\theta \in S) \exp (\pm \epsilon d_{**}),
	\end{equation}
    for all $S \in \F$, where $d_{**} = \sup_{x, x' \in \supp (x \mid t_0)} d(x,x')$.
\end{theorem}

\begin{proof}(sketch) 
    As in the proof of Theorem~\ref{thmDPHypothesisTesting}, we can show that
    \[\frac{p(t_0 \mid \theta)}{p(t_0 \mid \theta')} \in \exp (\pm \epsilon d_{**} ),\]
    for all $\theta, \theta' \in \supp (\pi)$. Plugging this into $\pi (\theta \mid t_0) = \frac{ p(t_0 \mid \theta) \pi(\theta) }{\int_{\Theta} p(t_0 \mid \theta') d\pi(\theta') }$ gives the result.
\end{proof}

Theorem~\ref{thmPosterior} demonstrates that the posterior $\pi(\theta \in \cdot \mid t_0)$ is in a probability interval which is centred at the prior $\pi(\theta \in \cdot)$ and has radius $\exp (\epsilon d_{**})$:
\[ \pi (\theta \in \cdot \mid t_0) \in \IIntP{L, U},\]
where $L = \pi(\theta \in \cdot) \exp( - \epsilon d_{**})$ and $U = \pi(\theta \in \cdot) \exp ( \epsilon d_{**})$.

\begin{remark}\label{remarkUniqueness}
    By following the proof of Theorem~\ref{thmPosterior}, one can observe that $\dMult(P_x, P_{x'})$ being bounded away from infinity, for all $x, x' \in \supp (x \mid t_0)$, is a necessary condition for 
    \[ \dMult \big[ \pi(\theta \mid t_0), \pi(\theta) \big] < \infty.\]
	(Note that \eqref{eq:posterior} is equivalent to $\dMult \left[ \pi(\theta \mid t_0), \pi(\theta) \right] \le \epsilon d_{**}$.) Indeed this condition is required for the posterior to be in a non-vacuous probability interval centred at the prior -- i.e. for the posterior to be in an probability interval of the form $\pi( \cdot \mid t_0) \in \IIntP{a \pi, b \pi}$ where $0 < a \le 1 \le b < \infty$. Hence the use of $\dMult$ in the Lipschitz condition~\eqref{eqDPDefnLip} is the unique choice (modulo distorting functions $\dist$ satisfying $\alpha \dist \ge \dMult$ for some constant $\alpha$) that ensures a bound on the prior-to-posterior of the form \eqref{eq:posterior}. This is analogous to the fact that $\dMult$ is the unique choice of distorting function that encodes privacy as indistinguishability (see Remark~\ref{remarkIndistinguishability}). Furthermore, by similar logic $\dMult$ is also the unique choice of distorting function which enables bounds on hypothesis testing like those in Theorem~\ref{thmDPHypothesisTesting}.
    
    These uniqueness properties are mirrored in the results of \citet{wasserman1992invariance} and \citet{lavine1991sensitivity} on the uniqueness of intervals of measures in robust Bayesian inference.
\end{remark}

Theorem~\ref{thmPosterior} contributes to what is called the \emph{prior-to-posterior semantics} of differential privacy (see \cite{kasiviswanathan2014semantics, dwork2006calibrating, duncanDisclosurelimitedDataDissemination1986}), in the sense that~\eqref{eq:posterior} describes the extent to which a Bayesian agent's posterior about a parameter $\theta$ can depart from their prior when learning from an $\epsilon$-DP data product.\footnote{An alternative type of semantics for differential privacy is the \emph{posterior-to-posterior semantics} \cite{dinur2003revealing,kasiviswanathan2014semantics}, whose focus is on the extent to which a Bayesian agent's posterior may vary were it derived from privacy-protected queries based on different (counterfactual) confidential databases. Previous literature in differential privacy predominantly adopted posterior-to-posterior semantics; see e.g. \cite{kifer2022bayesian}. However, prior-to-posterior semantics have recently attracted increasing attention as they circumvent counterfactuals and are closely connected with the literature on statistical disclosure risk; see e.g. \cite{gongCongenialDifferentialPrivacy2020,hotz2022balancing}.} Analogous to the discussion on frequentist attackers at the end of Section~\ref{sec:frequentist}, Theorem~\ref{thmPosterior} demonstrates the trade-off between restricting a Bayesian attacker while allowing for legitimate Bayesian learning: By setting $\Theta \subset \mathcal X$ and $P_{\theta}$ as degenerate point masses, we strictly generalise the result of \citet{gongCongenialDifferentialPrivacy2020} which bounds an attacker's prior-to-posterior change in a single record $x_i$.\footnote{By fixing some $x \in \mathcal X$ and setting $\pi(X_{-i} = x_{-i}) = 1$, we get $d_{**} = 1$ and thereby rederive the result from \cite{gongCongenialDifferentialPrivacy2020}. (Here $x_{-i}$ denotes the dataset $x$ -- which is assumed to be a vector -- with the $i$-th record removed.)

Alternatively, one could set $\pi(\theta) = \pi(x_i \mid x_{-i})$, in which case Theorem~\ref{thmPosterior} implies that
\[\pi(x_i \mid t, x_{-i}) \in \pi(x_i \mid x_{-i}) \exp(\pm \epsilon d_{**}),\]
(with $d_{**} = 1$ when, for example, $d = \dHam$).} Hence, we see that Theorem~\ref{thmPosterior} applies to both the legitimate analyst who is inferring population-level characteristics and the illegitimate attacker who is inferring individual-level information. Restricting the attacker (by decreasing $\epsilon$) necessarily hurts the analyst; whilst furnishing the analyst (by increasing $\epsilon$) also assists the attacker. What makes this dilemma tractable is that $d_{**}$ is typically much larger for the analyst than for the attacker because the analyst is interested in population quantities while the attacker is interested in individual records. Hence, Theorem~\ref{thmPosterior}'s bounds on the inference of the analyst are wider than those of the attacker. However, this argument breaks down when the analyst is interested in small subpopulations (such as in small-area estimation) because in these situations there is little light between the attacker's and the analyst's interests, and as such the values of $d_{**}$ associated with the attacker and the analyst will be similar. (This commentary -- on why the tradeoff between analysts and attackers is tractable -- also applies to Theorem~\ref{thmDPHypothesisTesting}.)

Theorem~\ref{thmPosterior} is powerful because it holds for arbitrary specifications of the data model $P_\theta$ and is applicable to the agent's arbitrary (proper) prior $\pi(\theta)$. So long as $d_{**}$ is finite (see the discussion after Theorem~\ref{thmMarginalData} on why this is not unreasonable), the bounds in~\eqref{eq:posterior} are non-trivial.

With that said, whenever $d_{**}$ is large, the bounds provided by Theorem~\ref{thmPosterior} are wide, rendering the results weakly informative at best. Indeed, rather than a pair of wide posterior bounds, the agent would be better off with a precise Bayesian posterior, which is theoretically derivable via the simple relation
\begin{equation}\label{eq:posterior-precise}
    \pi(\theta \mid t ) \propto \pi(\theta)p(t\mid \theta),
\end{equation}
where $p(t\mid \theta)$ can in turn be derived from the convolution of the data model $P_{\theta}$ and the privacy mechanism $P_x$ according to~\eqref{eq:marginal-data-probability}. In practice, however, direct computation or sampling from~\eqref{eq:posterior-precise} is not always possible or feasible. Such difficulties arise in situations A) where the privacy mechanism $P_x$ is not fully transparent to the analyst due to its complex dependence on $x$, whether by design or by post-processing \cite{gongTransparentPrivacyPrincipled2022}; B) where the data model $P_{\theta}$ is intractable, such as if defined algorithmically or treated as a black-box; or C) where their convolution~\eqref{eq:marginal-data-probability}, typically an $n$-dimensional integral, is intractable. Under any of these situations, the analyst may still rely on Theorem~\ref{thmPosterior} to obtain bounds on their posterior. 

Despite their width, these bounds are optimal whenever $\epsilon$ is the smallest constant satisfying the Lipschitz condition \eqref{eqDPDefnLip}. Without adding further assumptions on $M, P_\theta$, or $\pi$, these bounds cannot be shrunk. (This also applies to the bounds from Sections~\ref{sec:marginal} and \ref{sec:frequentist}. We prove this in Section~\ref{secTight}.) Yet they are not necessarily tight at a given $\theta$. This deficiency is an inevitable consequence of our analysis, which replaced the average case, $\int p_x(t) dP_\theta(x)$, with the extremal case, $p_{x_*}(t) \exp \left( \epsilon d_* \right)$. Such an analysis is necessarily loose whenever there is any variation away from the extreme. But the analysis cannot be tightened without making assumptions about the nature of this variation -- i.e. by making further assumptions on $M, P_\theta$, or $\pi$.

We illustrate the posterior bounds of Theorem~\ref{thmPosterior} with an example of Bayesian inference for a privatised count. 

\begin{example}[privatised single count] \label{ex:private-count}
Suppose the database consists of a single count record $x \in \mathbb{N}$. We wish to  query  the value of $x$ after it has been \emph{clamped} to a pre-specified range $[a_0, a_1]$. That is, $q(x) = a_0$ if $x < a_0$,  $q(x) = a_1$ if $x > a_1$, and  $q(x) = x$ otherwise. In differentially private mechanism design, clamping can be a necessary procedure when the intended query has otherwise unbounded global sensitivity. Under clamping, the sensitivity is reduced to $\Delta(q) = a_1 - a_0$.

The analyst's Bayesian model is
\begin{align*}
   \theta	&\sim \mathrm{Gamma}\left(\alpha,\beta\right), \\
x\mid\theta	&\sim Pois(\theta), \\
t\mid x	&\sim Lap\left(q\left(x\right);\epsilon^{-1}\Delta(q) \right). 
\end{align*}

For illustration, set $a_0 = 0$, $a_1 = 6$, $\alpha = 3$, $\beta = 1$. Figure~\ref{fig:posterior} depicts in blue solid lines the upper and lower density bounds on the analyst's posterior distribution $p(\theta \mid t)$ as given by Theorem~\ref{thmPosterior}. With $\epsilon = 1$ and $d_{**} = 1$, they are equal to the $\mathrm{Gamma}(3,1)$ prior density (blue dashed line), scaled by $\exp\left({\pm 1}\right)$. Overlaid in grey are Monte Carlo posterior densities $p(\theta \mid t^{(k)})$, $k = 1,\ldots, 10$, produced via the exact sampling algorithm proposed by~\cite{gong2022exact}. Each $t^{(k)}$ is independently simulated from the prior predictive distribution of the above Bayesian model.
\end{example}

\begin{figure}
	\floatconts{fig:posterior}
	{\caption{\emph{Density bounds for the posterior $p(\theta \mid t)$ from a privatised single count (Example~\ref{ex:private-count}).} The dashed blue line is the density of the $\mathrm{Gamma}(3, 1)$ distribution, the analyst's prior for $\theta$. In grey are simulation-based posterior densities based on 10 realizations of $t$ from its prior predictive distribution under the Poisson data model \citep{gong2022exact}. Upper and lower density bounds for the posterior $p(\theta \mid t)$ are in solid blue. The clamping range is $[0,6]$ and the privacy loss is $\epsilon = 1$.}}
	{\includegraphics[width=\textwidth]{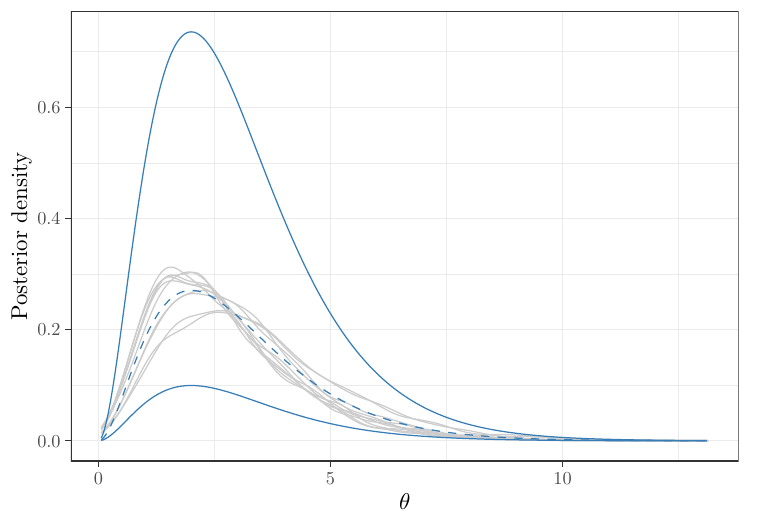}}
\end{figure}

Several aspects of Example~\ref{ex:private-count} are worth noting. First, the posterior density bounds (solid blue) are functions of the analyst's chosen prior $\pi\left(\theta\right)$ and the privacy mechanism parameters  $\epsilon$ and $d_{**}$ only. They are valid for any data model $P_\theta$ that the analyst wishes to employ, including (but not limited to) the Poisson data model that underlie the depicted precise posteriors densities $p(\theta \mid t^{(k)})$ in grey. On the other hand, while these precise posterior densities display moderate variations among each other, they do not depart much from the prior density (dashed blue). This is due to the heavy-handedness of the privacy mechanism employed for this analysis, resulting in poor statistical utility of the privatised count $t$. Indeed, the mechanism injects Laplace noise with standard deviation of $\sqrt{2}\epsilon^{-1}\Delta(q) = 8.48$ into a statistic clamped between $a_0 = 0$ and $a_1 = 6$. That $t$ cannot be highly informative for the inferential problem at hand is correctly identified by the full Bayesian analysis which precisely accounts for the uncertainty induced by the privacy mechanism (grey lines). Furthermore, these precise posterior distributions are generally far from the bounds implied by Theorem~\ref{thmPosterior}; this re-enforces the shallowness of these bounds due to their validity for very general classes of the data model $P_\theta$ and priors $\pi$.

\section{Pufferfish Privacy}\label{sec:pufferfish}

As the classic formulation of differential privacy, pure $\epsilon$-DP has inspired many variants, most of which closely resemble the original (for example by replacing $\dMult$ with some other distorting function $\dist$, see Remark~\ref{remarkDPCredalSet}). In contrast, Pufferfish privacy \citep{kifer2014pufferfish} is conceptually distinct from $\epsilon$-DP in two ways. Firstly, while $\epsilon$-DP conceptualises privacy as indistinguishability between pairs of comparable datasets $x, x' \in \mathcal X$ (Remark~\ref{remarkIndistinguishability}), Pufferfish reconceptualises privacy as indistinguishability between pairs of competing conjectures about the unobserved, confidential data $x$ (as we will see in Remark~\ref{remarkPuffIndistinguish}). Secondly -- and consequently -- $\epsilon$-DP  is concerned solely with the design of the data-release mechanism, while the object of Pufferfish's interest is the composition of the data-generating process and the data-release mechanism. We call this composite function the \emph{data-provision procedure}:
\begin{definition}
    Given a data-generating process $G(\theta, U_1)$ and a data-release mechanism $M(x, U_2)$, the \emph{data-provision procedure} $\DPP : \Theta \times [0,1]^2 \to \mathcal T$ is defined as 
    \[\DPP(\theta, U_1, U_2) = M(G(\theta, U_1), U_2),\]
    where $U_1, U_2$ are independent and (without loss of generality) identically distributed $\Unif[0,1]$.
\end{definition}

Here $U_1$ and $U_2$ are the random components (i.e. \emph{seeds}) of the data-generation $G$ and of the data-release $M$ respectively; %
$\theta$ is the data-generating model parameter; $X = G(\theta, U_1)$ is the (stochastic) dataset; and (as before) $T = M(X, U_2)$ is the released summary statistic. 

As in Section~\ref{sec:marginal}, we require that the data universe $\mathcal X$ is equipped with a $\sigma$-algebra $\mathscr G$; that the data-release mechanism $M$ is $(\mathscr G \otimes \mathcal B[0,1], \mathscr F)$-measurable; and, further, that $G(\theta, \cdot)$ is $(\mathcal B[0,1], \mathscr G)$-measurable for all $\theta \in \Theta$. (Recall that $\mathcal B[0,1]$ is the Borel $\sigma$-algebra on $[0,1]$ and $\mathscr F$ is the $\sigma$-algebra of the output space $\mathcal T$.) We now also assume that the data-generating model parameter set $\Theta$ is in a one-to-one correspondence with the set of probability measures on $(\mathcal X, \mathscr G)$. (So the data-generating model is not parametric in the usual sense of the term.)

Under these requirements, the seed $U_1$ induces a probability measure for the dataset $X$:
\begin{equation}\label{eqDefnProbPTheta}
    P_\theta(X \in E) = \lambda \big( \{u_1 \in [0,1] : G(\theta, u_1) \in E\} \big),
\end{equation}
(where $\lambda$ is the Lebesgue measure) and -- together with the seed $U_2$ -- also for the output $T$ of the data-provision procedure $\DPP$:
\begin{align}
    P(T \in S \mid \theta) &= \lambda \big( \{ u_1, u_2 \in [0,1] : M_G(\theta, u_1, u_2) \in S \} \big) \notag\\
    &= \int_{\mathcal X} P_x(S) dP_{\theta}(x),\label{eq:marginal-data-probability2}
\end{align}
(Equation~\eqref{eq:marginal-data-probability2} follows from the previous line by Fubini's theorem -- see Proposition~\ref{propDefnDataProvisionProb}.) Recall from Equation~\eqref{eq:marginal-data-probability} %
that $P(T \in \cdot \mid \theta)$ is the privatised data probability. Because Pufferfish is concerned with the data-provision procedure $\DPP$ -- and because $P(T \in \cdot \mid \theta)$ is the probability induced by $\DPP$ -- the privatised data probability plays a central role in Pufferfish. In fact, we will see (in Definition~\ref{defnPuffSat}) that Pufferfish can be conceived as a Lipschitz condition on the map $\theta \mapsto P(T \in \cdot \mid \theta)$, just as $\epsilon$-DP is a Lipschitz condition on the map $x \mapsto P_x$.

Pufferfish provides a framework for developing tailored privacy definitions. The data custodian constructs their custom Pufferfish privacy definition according to their judgement of:%
\begin{enumerate}[label=\Alph*.]
    \item The \emph{attackers}: Against what kinds of background knowledge, or beliefs about the data, should the data-release mechanism $M$ guard? (These knowledge and beliefs are modelled by probability distributions $\theta$ on the data $X$.)
    \item The \emph{attackers' conjectures on confidential information}: Which parts of the dataset require protection (i.e. what are the confidential information?), and what conjectures may an attacker make about these information? (Conjectures are modelled as events $E \in \mathscr G$ on the data universe $\mathcal X$.)
    \item The \emph{pairs of competing conjectures}: Which pairs of conjectures should remain indistinguishable to the attackers, even after observing (the realized value of) $T$? (Or, for a Bayesian attacker, between which pairs of conjectures should it be impossible for the attacker to significantly improve their ability to distinguish?) 
\end{enumerate}

(Note that %
B. is only used as a stepping-stone for answering C.; it does not have an independent role to play in Pufferfish.) Putting the above discussion more formally, the data custodian specifies their privacy definition with a \emph{Pufferfish instantiation}:

\begin{definition}\label{defnPuffInstant}
    A \emph{Pufferfish instantiation \Puff}\ is a tuple with three components:
    \begin{enumerate}
        \item A set $\mathbb D \subset \Theta$ of \emph{attackers};
        \item A set $\Spairs \subset \mathscr G \times \mathscr G$ of \emph{pairs of competing conjectures} about ``confidential'' information in the dataset $x$;\emph{\footnote{Note that elsewhere in the literature (for example in \cite{kifer2014pufferfish}), what we term the ``attackers' conjectures'' are referred to as ``secrets,'' and the ``pairs of competing conjectures'' are referred to as ``discriminative pairs.'' Moreover, in \cite{kifer2014pufferfish}, the set of discriminative pairs is denoted by $\mathbb S_{\mathrm{pairs}}$ and $\Spairs$ instead denotes the set of secrets/conjectures. We choose to omit the set of secrets/conjectures from a Pufferfish instantiation as it is superfluous, and instead use $\Spairs$ to denote the set of pairs of competing conjectures.}} and %
        \item A \emph{privacy loss budget} $0 \le \epsilon \le \infty$.
    \end{enumerate}
\end{definition}

The sets $\mathbb D$ and $\Spairs$ correspond to points A. and C. respectively in the above discussion. The privacy loss budget $\epsilon$ has the same role in Pufferfish as in $\epsilon$-DP: %
it describes the degree of continuity of the map $\theta \mapsto P(T \in \cdot \mid \theta)$ -- and hence, intuitively, the degree of privacy afforded to the data -- with smaller $\epsilon$ corresponding to more continuity/privacy.

Pufferfish privacy is a Lipschitz condition on $\DPP$:
\begin{definition}\label{defnPuffSat}
    Fix the data-generating process $G$. A data-release mechanism $M$ \emph{satisfies the Pufferfish instantiation \Puff}\ if the associated data-provision procedure $\DPP$ satisfies the inequality:
    \[ \dMult \Big( P(T \in \cdot \mid \theta), P(T \in \cdot \mid \theta') \Big) \le \epsilon \dPuff(\theta, \theta'),\]
    for all $\theta, \theta' \in \Theta$.
\end{definition}

\sepfootnotecontent{footnotePufferfishWellDefined}{\label{footnoteWellDefined}Pufferfish limits its consideration to $\theta \in \Theta$ and $E, E' \in \mathscr G$ for which $P_\theta (X \in \cdot \mid X \in E)$ and $P_{\theta} (X \in \cdot \mid X \in E')$ are well-defined, in order to ensure that the data-provision procedures $P(T \in \cdot \mid \theta, X \in E)$ and $P(T \in \cdot \mid \theta, X \in E')$ are themselves well-defined.

However, determining whether or not $P_\theta(X \in \cdot \mid X \in E)$ can be well-defined is beyond the scope of this paper. Answering this question to the necessary level of generality is difficult (see \cite{changConditioningDisintegration1997b} and references therein), but the majority of cases encountered in practice are covered by two approaches: When $P_\theta (X \in E)$ is non-zero, $P_\theta(X \in A \mid X \in E)$ is defined as $P_\theta(X \in A \cap E)/P_\theta(X \in E)$. And when $X = (Y, Z)$ has a canonical density $f(Y = y, Z = z)$ on a product measure $\mu = \mu_1 \times \mu_2$ with $\{X \in E\} = \{Z = z_E\}$ for some $z_E$, then $P_\theta (X \in \cdot \mid X \in E)$ is defined as the regular conditional probability $f_{Y \mid Z}(\cdot \mid z_E) d\mu_1(\cdot)$ where
\[f_{Y \mid Z}(y \mid z) = 
\begin{cases}
	\frac{ f(y,z) }{f(z)} & \text{if } f(z) > 0, \\
	\phi(y) & \text{otherwise,}
\end{cases}\]
with $f(z) = \int f(y,z) d\mu_1(y)$ and $\phi(y)$ an arbitrary density on $\mu_1$ \cite[Example~4.1.6]{durrettProbabilityTheoryExamples2019}.}

Here $\dPuff$ is a metric on $\Theta$ which is given by the Pufferfish instantiation \Puff. It is defined as follows: Firstly, as $\Theta$ is in a one-to-one correspondence with the set of probability measures on $(\mathcal X, \mathscr G)$, it is closed under conditioning. That is, for each $\theta \in \Theta$ and for each event $E \in \mathscr G$, there exists a unique $\theta' \in \Theta$ such that $P_{\theta'} (X \in \cdot) = P_\theta(X \in \cdot \mid X \in E)$ (provided that $P_\theta(X \in \cdot \mid X \in E)$ is well-defined\sepfootnote{footnotePufferfishWellDefined}). Denote this $\theta'$ by $\theta|_E$.

Then let the graph $\GPuff$ on $\Theta$ have edges $(\theta, \theta')$ if there exists some $\theta^* \in \mathbb D$ and some $(E, E') \in \Spairs$ such that $\theta = \theta^*|_E$ and $\theta' = \theta^*|_{E'}$ -- i.e. such that $P_\theta$ is equal to the conditional distribution $P_{\theta^*}(X \in \cdot \mid X \in E)$ and $P_{\theta'}$ is equal to $P_{\theta^*}(X \in \cdot \mid X \in E')$. Finally, define $\dPuff(\theta, \theta')$ as the length of a shortest path between $\theta$ and $\theta'$ in $\GPuff$.

Therefore, for all $\theta^* \in \mathbb D$ and all $(E, E') \in \Spairs$ (with $P_{\theta^*} (X \in \cdot \mid X \in E)$ and $P_{\theta^*} (X \in \cdot \mid X \in E')$ both well-defined\footref{footnoteWellDefined}), the data-generating probabilities $\theta^*|_E$ and $\theta^*|_{E'}$ are adjacent in the graph $\GPuff$ and hence 
\begin{equation}\label{eqPufferfishExplanation}
    \dMult \Big( P(T \in \cdot \mid \theta^*, X \in E), P(T \in \cdot \mid \theta^*, X \in E') \Big) \le \epsilon.
\end{equation}

The above discussion sheds light on the differences between Pufferfish and pure $\epsilon$-DP: As observed earlier, $\epsilon$-DP is concerned with indistinguishability of datasets $x, x' \in \mathcal X$. Hence, its starting point is the data universe $\mathcal X$ and it is a Lipschitz condition on %
the data-release mechanism $M : \mathcal X \times [0,1] \to \mathcal T$. On the other hand, Pufferfish is concerned with competing conjectures $\theta^*|_E$ and $\theta^*|_{E'}$ (for $\theta^* \in \mathbb D$ and $(E, E') \in \Spairs$). Its starting point is thus the data-generating parameter set $\Theta$ and it is a Lipschitz condition on %
the data-provision procedure $\DPP : \Theta \times [0,1]^2 \to \mathcal T$. %
Yet the data custodian only has partial control over $\DPP$. That is to say, while Pufferfish is a property of the data-provision procedure $\DPP$, %
the data custodian can achieve this property only through the design of $M$. In contrast, the data custodian often has full control of the object of $\epsilon$-DP's interest, the data-release mechanism $M$. 

While $\epsilon$-DP allows for the use of an arbitrary distance $d$, Pufferfish makes a very particular choice for the distance on its input space: $\dPuff$. (However, beyond the interpretation of $\dPuff$ in terms of attackers and competing conjectures, there is no reason in principle that Pufferfish cannot be generalised to allow for arbitrary distances $d$ on $\Theta$. In fact, all of the results below generalise immediately from $\dPuff$ to any metric $d$ on $\Theta$ which satisfies Assumption~\ref{assumptionGraphDistance}.)

Despite their differences, $\epsilon$-DP and Pufferfish share one characteristic which is very important for our purposes: They both use the multiplicative distance $\dMult$ to measure the change in output variations.
This means that Pufferfish is fundamentally linked to the concept of an interval of measures, just as $\epsilon$-DP is. This connection to intervals of measures and the resulting implications %
on the indistinguishability of important inferential quantities are the subject of the next section.

\section{An IP View of Pufferfish Privacy}\label{sec:pufferfish-IP}

As an analog to Theorem~\ref{thm:dp-iom} for pure $\epsilon$-differential privacy, Theorem~\ref{thm:pufferfish-iom}  below establishes the connection between Pufferfish and intervals of measures. Specifically, \ref{thm:pufferfish-iom}.\ref{puffThmStatement2} is the standard definition of Pufferfish as given in \cite{kifer2014pufferfish}. The equivalence \ref{thm:pufferfish-iom}.\ref{puffThmStatement1} $\Leftrightarrow$ \ref{thm:pufferfish-iom}.\ref{puffThmStatement2} justifies the formulation of Pufferfish as Lipschitz continuity. In addition, \ref{thm:pufferfish-iom}.\ref{puffThmStatement3} and \ref{thm:pufferfish-iom}.\ref{puffThmStatement4} give novel formulations of Pufferfish in terms of intervals of measures.

\sepfootnotecontent{footenotePuffThmStatement4}{Analogous to the concept of connected data $x, x' \in \mathcal X$ (Definition~\ref{defnConnected}), we say that $\theta, \theta' \in \Theta$ are \emph{$\dPuff$-connected} if $\dPuff(\theta, \theta') < \infty$ and we define $[\theta]$ to be $\theta$'s connected component: $[\theta] = \{\theta' \in \Theta \mid \dPuff(\theta, \theta') < \infty\}$.} 

\begin{theorem}\label{thm:pufferfish-iom}
    Fix the data-generating process $G(\theta, U_1)$ and the data-release mechanism $M(x, U_2)$. %
    For any Pufferfish instantiation \Puff\ with privacy loss budget $\epsilon < \infty$, 
    the following statements are equivalent:
    \begin{enumerate}[label=\Roman*]
        \item\label{puffThmStatement1} $M$ satisfies %
        \Puff.
        \item\label{puffThmStatement2} For all $S \in \mathscr F$, all competing conjectures $(E, E') \in \Spairs$ and all attackers $\theta^* \in \mathbb D$ (such that $P_{\theta^*}(X \in \cdot \mid X \in E)$ and $P_{\theta^*}(X \in \cdot \mid X \in E')$ are both well-defined\footref{footnoteWellDefined}), the following inequalities are satisfied:
        \begin{align*}
            P(T \in S \mid \theta^*, X \in E) &\le e^\epsilon P(T \in S \mid \theta^*, X \in E'), \\
            P(T \in S \mid \theta^*, X \in E' ) &\le e^{\epsilon} P(T \in S \mid \theta^*, X \in E).
        \end{align*}
        \item\label{puffThmStatement3} For all $\delta \in \mathbb N$ and all $\theta, \theta' \in \Theta $ with $\dPuff(\theta, \theta') = \delta$,
        \begin{equation*}
            P(T \in \cdot \mid \theta') \in \IIntP{L_{\theta,\delta\epsilon},U_{\theta,\delta\epsilon}},
        \end{equation*} 
        where $L_{\theta,\delta\epsilon} = e^{-\delta\epsilon}P(T \in \cdot \mid \theta)$ and $U_{\theta,\delta\epsilon} = e^{\delta\epsilon }P(T \in \cdot \mid \theta)$.
        \item\label{puffThmStatement4} For all $\theta \in \Theta$ and all measures $\nu \in \Omega$, if $P(T \in \cdot \mid \theta)$ has a density $p(t \mid \theta)$ with respect to $\nu$, then for every $\dPuff$-connected\sepfootnote{footenotePuffThmStatement4} $\theta' \in [\theta]$, $P(T \in \cdot \mid \theta')$ also has a density $p(t \mid \theta')$ (with respect to $\nu$) satisfying
        \begin{equation*}%
            p(t \mid \theta') \in p(t \mid \theta) \exp \left[\pm \epsilon \dPuff(\theta,\theta') \right],
        \end{equation*}
        for all $t \in \mathcal T$.
    \end{enumerate}
\end{theorem}

\begin{proof}
    First note that \eqref{eqPufferfishExplanation} is equivalent to both \ref{puffThmStatement1} and \ref{puffThmStatement2}. Specifically, that \eqref{eqPufferfishExplanation} implies \ref{puffThmStatement1} follows by applying the triangle inequality to a shortest path between $\theta$ and $\theta'$ in $\GPuff$, similar to the proof of \ref{thm:dp-iom}.\ref{mainThmStatement1} $\Leftrightarrow$ \ref{thm:dp-iom}.\ref{mainThmStatement2} in Theorem~\ref{thm:dp-iom}. %
    The remainder of the proof is analogous to that of Theorem~\ref{thm:dp-iom}. 
\end{proof}

\begin{remark}\label{remarkPuffIndistinguish}
Statement \ref{mainThmStatement4} of Theorem~\ref{thm:dp-iom} is the backbone for reasoning about indistinguishability between $d$-connected datasets $x, x' \in \mathcal X$ under $\epsilon$-DP (see Remark~\ref{remarkIndistinguishability}). In contrast, statement \ref{puffThmStatement4} of Theorem~\ref{thm:pufferfish-iom} 
provides the rationale for indistinguishability between $\dPuff$-connected distributions $\theta, \theta' \in \Theta$ under Pufferfish privacy. Specifically, for $\theta = \theta^*|_E$ and $\theta' = \theta^*|_{E'}$ (with $\theta^* \in \mathbb D$ and $(E, E') \in \Spairs$), an attacker cannot distinguish with certainty between $\theta$ and $\theta'$ because $p(t \mid \theta)/p(t \mid \theta')$ is bounded away from zero and infinity, regardless of the value of $t$. More generally, whenever $\theta$ is plausible (i.e. when $p(t \mid \theta) > 0$) then all $\dPuff$-connected $\theta' \in [\theta]$ are also plausible (regardless of the choice of dominating measure $\nu$).
\end{remark}

We now turn to discussing the impact of Pufferfish privacy on statistical inference in both Bayesian and frequentist paradigms. From the Bayesian view, Pufferfish limits the ability of an attacker $\theta^* \in \mathbb D$ to discern between two competing conjectures $(E, E') \in \Spairs$ relative to their prior (baseline) ability to do so:
\begin{equation}\label{eqPuffBayes}
    e^{-\epsilon} \le \left. \frac{P_{\theta^*} (X \in E \mid T = t)}{P_{\theta^*} (X \in E' \mid T = t)} \middle/ \frac{P_{\theta^*} (X \in E)}{P_{\theta^*}(X \in E')} \right. \le e^\epsilon,
\end{equation}
where the attacker's ``ability to discern between $(E, E')$'' is quantified as the odds of $E$ against $E'$, so that \eqref{eqPuffBayes} is a bound on the prior-to-posterior odds ratio. In fact, $M$ satisfies \Puff\ if and only if $M$ satisfies \eqref{eqPuffBayes} for all $\theta^* \in \mathbb D$, all\footref{footnoteWellDefined} $(E, E') \in \Spairs$ and almost all $t \in \mathcal T$. (This result was first described in \cite[p.~6]{kifer2014pufferfish} %
and follows from Statement \ref{thm:pufferfish-iom}.\ref{puffThmStatement4} by setting $\theta = \theta^*|_{E}$ and $\theta' = \theta^*|_{E'}$, and then applying Bayes rule. We formally state and prove this result in Proposition~\ref{propPuffOddsRatio} of Appendix~\ref{appendixAdditionalResults}.)

As previous sections contend, there is an important type of competing conjectures that privacy mechanisms aim to make indistinguishable. These conjectures concern the values of the records in the dataset $x$. To this end, Pufferfish provides a Bayesian semantic guarantee that conforms to the structure of a \emph{density ratio neighbourhood} \cite{wassermanComputingBoundsExpectations1992, wasserman1992invariance}, defined below. 

Recall (from Definition~\ref{def:iom}) that $\Omega$ is the set of $\sigma$-finite measures on the measurable space $(\mathcal T, \mathscr F)$.

\begin{definition}\label{defnDensityRatio} 
The \emph{density ratio neighbourhood} of $\mu \in \Omega$ with radius $r \ge 0$ is defined as
\[N_r(\mu) = \{ \nu \in \Omega : \drm (\mu,\nu) \le r\},\]
where $\drm$ is the \emph{density ratio metric}:
\begin{equation}\label{eqDefnDRM}
    \drm(\mu,\nu) = \begin{cases}
    0 &\mathrm{if\ } \mu = \nu = 0\\% \mathrm{\ are\ zero,}\footnotemark{}}\\
    \esssup_{t, t' \in \mathcal T^{\mathrm{o}}} \ln \left( \frac{f(t)}{f(t')} \middle/ \frac{g(t)}{g(t')} \right) &\mathrm{else\ if\ } \mu,\nu \mathrm{\ are\ mutually\ absolutely\ continuous},\\
    \infty &\mathrm{otherwise,}
\end{cases}
\end{equation}
with $f$ and $g$ densities of $\sigma$-finite measures $\mu$ and $\nu$ respectively, with respect to some common dominating measure $\tau \in \Omega$; $\mathcal T^{\mathrm{o}} = \{ t \in \mathcal T \mid 0 < f(t), g(t) < \infty \}$; and the essential supremum is with respect to $\tau$.\footnote{The property $f, g < \infty$ holds $\tau$-almost everywhere (because $\mu$ and $\nu$ are $\sigma$-finite -- see the proof of Lemma~\ref{lemmaEquivalenceBetweenSupEssSup}), and $f, g > 0$ holds $\mu$- and $\nu$-almost everywhere. Hence, practically one may take the essential supremum in equation~\eqref{eqDefnDRM} over $\mathcal T$; restricting to $\mathcal T^{\mathrm{o}}$ simply removes the complications of dividing by zero or infinity.}
\end{definition}

The definition of the density ratio metric $\drm$ is well-defined in the sense that $\drm(\mu, \nu)$ does not depend on the choice of $f,g$ and $\tau$ in~\eqref{eqDefnDRM}. (See Appendix~\ref{appendixDRM} for details.)

The following theorem characterises Pufferfish privacy (under a particular choice of $\Spairs$) as the requirement that an attacker $\theta$'s posterior on $X$ is in the $\epsilon$-density ratio neighbourhood of their prior:

\begin{theorem}\label{thm:pufferfish-drn} 
Fix some $\theta^* \in \mathbb D$. Let $\mathcal S_{\mathcal X}$ be a partition of $\mathcal X$ such that $P_{\theta^*} (X \in E)$, for $E \in \mathcal S_{\mathcal X}$, is given by a density $p_{\theta^*}(Z = z_{E})$ of some marginalisation $Z$ of $X$. Define $\Spairs = \mathcal S_{\mathcal X} \times \mathcal S_{\mathcal X}$. 

If $M$ satisfies \Puff, then 
\begin{equation}\label{eqthmPufferfishDrn}
    P_{\theta^*} (Z \in \cdot \mid T = t) \in N_{\epsilon} \big( P_{\theta^*} (Z \in \cdot) \big),
\end{equation}
for $P(T \in \cdot \mid {\theta^*})$-almost all $t \in \mathcal T$.

In the other direction, suppose that $P_\theta(X \in E)$, for $E \in \mathcal S_{\mathcal X}$, is given by a density $p_\theta(Z = z_E)$ for all $\theta \in \mathbb D$. Then \eqref{eqthmPufferfishDrn} holding for all $\theta \in \mathbb D$ and $P(T \in \cdot \mid \theta)$-almost all $t \in \mathcal T$ implies that $M$ satisfies \Puff.

\end{theorem}

The proof of Theorem~\ref{thm:pufferfish-drn} is immediate from~\eqref{eqPuffBayes}.

Two special cases of Theorem~\ref{thm:pufferfish-drn} are worth noting. Firstly, when the partition $\mathcal S_{\mathcal X} = \{ \{x\} : x \in \mathcal X\}$ consists of all the singleton subsets of $\mathcal X$, an \Puff\ mechanism is tasked with providing indistinguishability between competing conjectures $E = \{x\}$ and $E' = \{x'\}$. That is, a Pufferfish mechanism must protect against the conjecture $X = x$ versus $X = x'$, for any arbitrary choices of $x, x' \in \mathcal X$. %
This is a tall order, because the dataset may contain a large number of individual records, each with values $x_i, x_i'$ that are nothing alike. The resulting privacy guarantee is thus \label{review2.5}a stringent one: For any $\theta \in \mathbb D$,
\[P_{\theta \mid t} \in N_{\epsilon} \big( P_{\theta} \big),\]
where $P_\theta$ and $P_{\theta \mid t}$ are $P_{\theta} (X \in \cdot)$ and $P_{\theta}(X \in \cdot \mid T = t)$ respectively. 
In other words, the Bayesian attacker $\theta$'s ability to discern between two arbitrary datasets relative to their prior discernability (i.e. the prior-to-posterior odds ratio) is limited by a multiplicative factor of $e^\epsilon$. %
(This is closely related to the `no-free-lunch privacy' of \cite{kiferNoFreeLunch2011} and \cite[Section~3.2]{kifer2014pufferfish}.)

A second, and more pragmatic, special case arises when the partition  $\mathcal S_{\mathcal X}$ consists of the level sets given by fixing a small number of records in the dataset -- in particular, by fixing a single record. Assume that the datasets $x \in \mathcal X$ are vectors $(x_1, \ldots, x_n)$ of records and let $\mathcal R$ be the set of all possible values that a record $x_i$ can take, so that $\mathcal X \subset \bigcup_{n = 1}^\infty \mathcal R^n$. Define 
\begin{equation}\label{eqSXDefn}
    \mathcal S_{\mathcal X} = \{E(r,1) : r \in \mathcal R\},
\end{equation} 
where $E(r,i)$ is the level set which fixes the $i$th record to be the value $r$:
\begin{equation}\label{eqEriDefn}
    E(r,i) = \{ x \in \mathcal X : x_i = r\}.
\end{equation}
For this choice of $\mathcal S_{\X}$, the two densities $p_{\theta}(Z = z_{E(r,1)})$ and $p_{\theta}(Z = z_{E(r,1)} \mid T = t)$ are, respectively, the prior marginal density, $p_\theta(X_1 = r)$, and the posterior marginal density, $p_\theta(X_1 = r \mid T = t)$,  of the first record taking the value $r$, where the marginalisation is over all the other records with respect to the data-generating process $\theta$. Theorem~\ref{thm:pufferfish-drn} states that these prior and posterior marginal densities are restricted to the same density ratio neighbourhood of radius $\epsilon$.

\begin{remark}\label{remarkDPPuff}
    When the competing conjectures $\Spairs$ are given by the level sets, \Puff\ has a close connection to pure $\epsilon$-DP. Indeed, suppose that $\X = \mathcal R^n$ for a fixed $n$ and let $\Spairs = \bigcup_{i=1}^n \{E(r,i) : r \in \mathcal R\}^2$. If $\mathbb D$ is the collection of distributions on $\X$ which take the records $X_1,\ldots,X_n$ as mutually independent, then a mechanism $M$ satisfying \Puff\ is equivalent to $M$ satisfying $\epsilon$-DP with $d = \dHam$ \cite[Theorem~6.1]{kifer2014pufferfish}. This result follows from observing
    \[P(T \in S \mid \theta, X_i = x_i) = \int_{\mathcal R^{n-1}} P_x(S) dP_\theta(X_{-i} = x_{-i}),\]
    for $\theta \in \mathbb D$, which implies that Statements~\ref{thm:dp-iom}.\ref{mainThmStatement2} and \ref{thm:pufferfish-iom}.\ref{puffThmStatement2} are equivalent.
\end{remark}

Pufferfish also has a frequentist interpretation as a limit to the power of any attacker's level-$\alpha$ test between competing conjectures (c.f. Theorem~\ref{thmDPHypothesisTesting}): 

\begin{theorem}\label{thmPufferfishHypothesis}
    A data-release mechanism $M$ satisfies \Puff\ if and only if, for all $\theta_0 \ne \theta_1 \in \Theta$, the power $1-\beta$ of all size-$\alpha$ tests of 
    \[H_0 : \theta = \theta_0 \mathrm{\ versus\ } H_1 : \theta = \theta_1,\]
    is bounded by the inequalities:
    \begin{equation}\label{eqThmPufferfishHypothesis}
        \max \big( \alpha/\phi, 1 - [1-\alpha]\phi \big) \le 1 - \beta \le \min \big( \alpha \phi, 1 - [1 - \alpha]/\phi \big),
    \end{equation}
    where $\phi = \exp [ \epsilon \dPuff(\theta_0, \theta_1) ]$.
\end{theorem}

\begin{proof}
    The result is immediate from \ref{puffThmStatement4} of Theorem~\ref{thm:pufferfish-iom} and the Neyman-Pearson lemma. (The derivation is analogous to the second half of the proof of Theorem~\ref{thmDPHypothesisTesting} and the proof of Corollary~\ref{corDPHT}.)
\end{proof}

As an application of Theorem~\ref{thmPufferfishHypothesis}, suppose an attacker $\theta^* \in \mathbb D$ is interested in testing the null hypothesis $X \in E$ against the alternative $X \in E'$, for some $(E, E') \in \Spairs$. This is equivalent to setting $\theta_0 = \theta^*|_{E}$ and $\theta_1 = \theta^*|_{E'}$ in Theorem~\ref{thmPufferfishHypothesis} (assuming that $\theta^*|_{E}$ and $\theta^*|_{E'}$ are well-defined\footref{footnoteWellDefined}). Hence, any such test with level 0.05 will have power at most $0.05\exp (\epsilon)$ under \Puff.

\begin{remark}\label{remarkDRNCOR}
    Because the density ratio neighbourhood $N_r(\mu)$ is the closed ball $B_{\drm}^r(\mu)$, it is a distortion model (Remark~\ref{remarkDPCredalSet}). Moreover, it is closely related to the \emph{constant odds ratio model} \cite{walley1991statistical}, which is the distortion model associated with the distorting function $\dcor$ \cite{montesUnifyingNeighbourhoodDistortion2020}. Here, $\dcor$ is the constant odds ratio metric:
    \[\dcor (\mu, \nu) = \begin{cases}
        0 &\mathrm{if\ } \mu = \nu = 0, \\
        1 - \inf_{S, S' \in \F^*} \frac{\mu(S) \nu(S')}{\mu(S') \nu(S)} &\mathrm{else\ if\ } \mu, \nu \mathrm{\ are\ mutually\ absolutely\ continuous,}\\ 
        1 &\mathrm{otherwise,}
    \end{cases}\]
    for finite $\mu, \nu \in \Omega$, where $\F^* = \{ S \in \F : \mu(S) > 0\}$.

    Corollary~\ref{corDRMDCOR} in Appendix~\ref{appendixAdditionalResults} proves that, for finite $\mu, \nu \in \Omega$,
    \[\drm(\mu, \nu) = - \ln \left[ 1 - \dcor (\mu, \nu) \right].\]
    Hence, when restricting to finite measures, the density ratio neighbourhood $N_r$ is equal to the constant odds ratio model with distortion parameter $\delta = 1 - \exp( - r ).$
\end{remark}

\section{Optimality of This Paper's Results}\label{secTight}

The bounds presented in this paper cannot be improved without additional assumptions on the data-release mechanism $M$, the data-generating model $P_\theta$ or the prior $\pi$. In this section, we provide examples which demonstrate the optimality of these bounds. However, it is important to reiterate that these bounds are only tight pointwise. Indeed it would be impossible for it to be otherwise, since the bounds are on probability measures, yet the bounds themselves are not probability measures.

Throughout this section, we rely on the Laplace mechanism $M$ for the count query $q(x) = \sum_i x_i$ (Example~\ref{ex:laplace}). The density of $T \sim M(x, U)$ is $p_x(t) = \frac{\epsilon}{2} \exp \left(-\epsilon \abs{t - q(x)} \right)$ when $\mathcal X = \{0,1\}^n$ and the metric $d$ on $\mathcal X$ is the Hamming distance $\dHam$. We assume that $n$ is fixed, so that $\X$ is $\dHam$-connected and $\sup_{x, x' \in \X} \dHam(x,x') = n < \infty$.

We begin with Theorem~\ref{thmMarginalData} which states that the privatised data probability $P(T \in \cdot \mid \theta)$ is bounded by $L_{\theta,\epsilon}$ and $U_{\theta, \epsilon}$ on $\T_0$. Because $p(t \mid \theta) = \int p_x(t) dP_\theta(x)$ for a.e. $t \in \T_0$, the lower bound $p(t \mid \theta) \ge l_{\theta, \epsilon}(t)$ is tight if $p_x(t) = l_{\theta, \epsilon}(t)$ for $P_\theta$-a.e. $x \in \supp (x \mid t_0, \theta)$. Consider the Laplace mechanism under the setting given above. In this case, $\T_0 = \mathbb R$ for any $t_0$ and any $\theta$. Further, the essential-supremum in $l_{\theta, \epsilon}(t)$ is achieved by $x_* = (1,\ldots, 1)$  when $t \le 0$. Hence $l_{\theta, \epsilon}(t) = p_{x_0}(t)$ where $x_0 = (0, \ldots, 0)$. Therefore, $p(t \mid \theta)$ can be arbitrarily close to $l_{\theta, \epsilon}(t)$ for $t \le 0$ as $P_{\theta}$ concentrates on $x_0$. This implies $P(T \in S \mid \theta)$ can be arbitrarily close to $L_{\theta, \epsilon}(S)$ for a bounded, measurable set $S \subset \mathbb R^{\le 0}$. The upper bound $P(T \in \cdot \mid \theta) \le U_{\theta, \epsilon}$ follows similarly, by considering $t \ge n, x_* = (0,\ldots,0)$ and $x_0 = (1,\ldots,1)$.

We now move to Theorem~\ref{thmDPHypothesisTesting} which concerns the power of hypothesis tests in the private setting. To see that this result is tight, consider the model $\mathcal P = \big\{ P_{\theta} \mid \theta \in \{0,1\}^n \big\}$ where $P_\theta(X \in \cdot)$ is the point mass on $x = \theta$. Set $\theta_0 = (0,\ldots,0)$ and $\theta_1 = (1,\ldots, 1)$. By examining the density $p_x(t)$ of the Laplace mechanism for $x = (0,\ldots,0)$ and for $x = (1,\ldots,1)$, one can conclude that the Neyman-Pearson (NP) test must have a rejection region of the form $R = \{t > t_1\}$ for some $t_1$. Moreover, for small enough $\epsilon$, $t_1$ must be at least $n = d_{**}$ (assuming $\alpha < 0.5$). Then, $p(t \mid \theta_1) = \exp (\epsilon n) p(t \mid \theta_0)$ for all $t \in R$, which means the NP test has power exactly $\alpha \exp (\epsilon n)$. 

Theorem~\ref{thmPriorPredictive} provides bounds on a Bayesian analyst's prior predictive probability. If one sets the prior $\pi$ to be a point mass on a single $\theta_0$, then the prior predictive probability $P(T \in S) = \int_{\Theta} P(T \in S \mid \theta) d\pi(\theta)$ reduces to the privatised data probability $P(T \in S \mid \theta_0)$. In this case, proving optimality of Theorem~\ref{thmPriorPredictive} is analogous to proving that the bounds $L_{\theta_0, \epsilon} \le P(T \in \cdot \mid \theta_0) \le U_{\theta_0, \epsilon}$ are tight. Hence, the argument outline above for Theorem~\ref{thmMarginalData} can also be used to show optimality of Theorem~\ref{thmPriorPredictive}.

For Theorem~\ref{thmPosterior} -- which demonstrates that a Bayesian's posterior is within a probability interval of the prior -- take $\Theta = [0,1]$ with the prior $\pi = \mathrm{Unif}[0,1]$. Let $P_{\theta}(x)$ be the point mass on $(1,\ldots,1)$ if $\theta = 1$ and the point mass on $(0,\ldots,0)$ otherwise. For $t > n$, we have $\pi(\theta = 1 \mid t) = \pi (\theta = 1) \exp (\epsilon n)$. Thus, the bound in Theorem~\ref{thmPosterior} is achieved since $d_{**} = n$. 

Finally, we prove that our results on the inferential limits induced by Pufferfish privacy (Theorems~\ref{thm:pufferfish-drn} and~\ref{thmPufferfishHypothesis}) are tight. For this, consider the setting described in Remark~\ref{remarkDPPuff}. In this case, the Laplace mechanism $M$ with $\X = \{0,1\}^n$ satisfies \Puff. Theorem~\ref{thm:pufferfish-drn} states that, for any $i \in \{1,\ldots,n\}$ and any $\theta^* \in \mathbb D$, the posterior $P_{\theta^*}(X_i \in \cdot \mid T = t)$ is in the density ratio neighbourhood of radius $\epsilon$ that is centred at the prior $P_{\theta^*}(X_i \in \cdot)$. Let $\theta^*$ be such that
\[P_{\theta^*}(X = x_0) = P_{\theta^*}(X = x_1) = 0.5,\]
where $x_0 = (0,\ldots,0)$ and $x_1 = (1,0,\ldots,0)$. Then
\[ \left. \frac{P_{\theta^*} (X_1 = 0 \mid T = t)}{P_{\theta^*} (X_1 = 1 \mid T = t)} \middle/ \frac{P_{\theta^*} (X_1 = 0)}{P_{\theta^*}(X_1 = 1)} \right. = \frac{p_{x_0}(t)}{p_{x_1}(t)},\]
by Bayes rule%
, where $p_x(t)$ is the density of the Laplace mechanism. Yet $p_{x_0}(t)/p_{x_1}(t) = \exp(\epsilon)$ for $t \le 0$. Hence,
\[ \drm \bigg( P_{\theta^*}(X_1 \in \cdot \mid T = t), P_{\theta^*}(X_1 \in \cdot) \bigg) = \epsilon,\]
and thus the bound~\eqref{eqthmPufferfishDrn} of Theorem~\ref{thm:pufferfish-drn} is tight. 

Now we consider Theorem~\ref{thmPufferfishHypothesis}, which provides a bound on the power of any size-$\alpha$ test. As before, we use the Pufferfish instantiation given in Remark~\ref{remarkDPPuff} and the Laplace mechanism $M$. Let $\theta_0$ and $\theta_1$ be such that $P_{\theta_i}(X \in \cdot)$ is the point mass on $x_i$, where $x_0 = (0, \ldots, 0)$ and $x_1 = (1, \ldots, 1)$ with $\abs{x_0} = \abs{x_1} = n$. Then $\dPuff(\theta_0, \theta_1) = n$ and, furthermore, the test $H_0 : \theta = \theta_0$ versus $H_1 : \theta = \theta_1$ is exactly the test we examined when proving the optimality of Theorem~\ref{thmDPHypothesisTesting}. In that proof, we demonstrated that the Neyman-Pearson test has power $\alpha \exp (\epsilon n)$ (assuming that $\epsilon$ is small and $\alpha < 0.5$). This implies the bounds~\eqref{eqThmPufferfishHypothesis} in Theorem~\ref{thmPufferfishHypothesis} are tight.

\section{Discussion}\label{sec:discussion}

The results we obtain in this paper make novel contributions to the differential privacy literature in the following ways. Firstly, the bounds we obtain in Theorems~\ref{thmMarginalData},~\ref{thmDPHypothesisTesting},~\ref{thmPriorPredictive},~\ref{thmPosterior},~\ref{thm:pufferfish-drn} and~\ref{thmPufferfishHypothesis}  are non-trivial, due to the validity of these results across a broad range of data models, privacy mechanisms and prior distributions. When the analyst has little knowledge or is only willing to make minimal assumptions about their model, these bounds are useful representations of the limits of statistical learning under privacy constraints. This draws a contrast with the existing DP literature, which has largely focused on asymptotic lower bounds or on constructing (asymptotically-)optimal data-release mechanisms for specific data use cases \cite{smithPrivacypreservingStatisticalEstimation2011, caiCostPrivacyOptimal2021, chhorRobustEstimationDiscrete2023, duchi2018minimax, bassilyPrivateEmpiricalRisk2014, talwarNearlyOptimalPrivate2015, dworkAnalyzeGaussOptimal2014, wasserman2010statistical, awan2020differentially}. This literature aligns with \emph{query-based access} \cite{hotz2022balancing} where the user can choose what statistics are released. Our results, on the other hand, are finite-sample and apply to the \emph{dissemination mode} of data release where the mechanism is not tailored for the analyst's use case. This setting is typical of official statistics (e.g. censuses and surveys) and, more generally, data products with multiple users, and is more common in the research community than query-based access \cite{hotz2022balancing}.

Secondly, the generality of our bounds implies that they are inherent consequences of the privacy standards themselves, be it pure $\epsilon$-DP or Pufferfish.
Specifically, 
these bounds stem only from the requirement that the mechanism $M$ is $\epsilon$-DP -- or \Puff\ -- and not on any particularities of $P_\theta, M$ or $\pi$. That these bounds are typically wide in practice -- as can be seen from Examples~\ref{ex:binary-sum} and~\ref{ex:private-count} -- is in part due to the near-total lack of assumption under which they are derived. While these bounds can approach vacuity as the data size $n$ grows, in practical examples that need not be the case if, for example, the data analyst has probabilistic knowledge about the privacy mechanism (see e.g. Example~\ref{ex:private-count}) or the data space $\mathcal X$. For a given choice of $P_\theta, M$ and $\pi$, we may obtain tighter bounds than those in this paper. In addition to the asymptotic results in the aforementioned papers, sharp bounds for specific $P_\theta, M$ and $\pi$ may be derivable from the existing literature on measurement errors (\emph{errors-in-variables}) in statistics and econometrics, particularly in the case of point identification problems (see e.g. 
\citep{carroll1988optimal, horowitzIdentificationRobustnessContaminated1995}).

Through the lens of Theorems~\ref{thmDPHypothesisTesting},~\ref{thmPosterior},~\ref{thm:pufferfish-drn} and~\ref{thmPufferfishHypothesis}, we obtain valuable insights in both frequentist and Bayesian paradigms on the \emph{privacy-utility trade-off} which is fundamental to differential privacy as a quantitative privacy standard \cite{hotz2022balancing}. (For details, see the discussion accompanying these theorems.) Qualitatively speaking, there exists an inherent tension between protecting private information and deriving scientific knowledge. To date, quantitative approaches to this trade-off predominantly rely on the privacy loss budget as the sole metric to balance this trade-off \cite{abowd2019economic, hsuDifferentialPrivacyEconomic2014, heffetzWhatWillIt2022}. However, from the suite of IP analyses presented here, we see that other building blocks  -- notably the metric structure $(\mathcal X, d)$ of the data universe, the associated database distances (such as $d_*$ and $d_{**}$), and the clamping parameters on $\X$ -- are all relevant factors that, together with the privacy loss budget $\epsilon$, collectively determine the limits to statistical learning for attackers and scientists alike. Therefore, $\epsilon$ is not the only parameter of concern -- and perhaps not even the central concern -- when assessing and trading-off privacy and utility \cite{bailieRefreshmentStirredNot2024a}.

While this paper qualifies the tradeoff of $\epsilon$-DP and Pufferfish privacy in concrete terms, it narrowly conceives privacy and utility in terms of the extent of statistical estimation attainable under either frequentist or Bayesian paradigms. There are, of course, many other aspects of utility that are worth examining -- such as the ease of analysis, use of computational resources, facial validity and logical consistency \cite{boyd2022differential, hotz2022chronicle, rugglesDifferentialPrivacyCensus2019} -- and other paradigms (in particular decision theory) with which the notions of privacy and utility can be quantified. In fact, both notions are multi-faceted and context-specific and, as one of the reviewers of \cite{bailieDifferentialPrivacyGeneral2023} pointed out, a judicious conceptualisation of privacy and utility may improve their tradeoff's efficiency frontier.  %
Acknowledging the complex makeup of this tradeoff, we advocate for future design and analysis of data-release mechanisms to treat the conceptualisations of privacy and utility -- and the roles that $\mathcal X, d$ and the distorting function $\dist$ play in these conceptualisations -- with scrutiny, given their scarcity in the current literature \citep{bailieRefreshmentStirredNot2024b}.

Tools from the IP literature harbour potential in aiding future endeavours to study statistical data privacy in a rigorous yet general manner. This work has examined some examples in which a DP definition can be formulated as the requirement that a mechanism's probability $P_x$ under one input $x \in \X$ is in a certain distortion model of its probability $P_{x'}$ under another input $x'$, whenever those two inputs $x, x'$ are connected. In this case, the choice of %
distorting function $\dist$ (partially) determines the flavor of the privacy guarantee. One direction for future research is thus to explore IP characterisations of other common flavors of DP, in particular $(\epsilon,\delta)$-DP \citep{dworkOurDataOurselves2006, machanavajjhala2008privacy, kifer2022bayesian}, zero-concentrated DP (zCDP) \citep{dworkConcentratedDifferentialPrivacy2016, bunConcentratedDifferentialPrivacy2016a}, R\'enyi DP  \citep{mironovRenyiDifferentialPrivacy2017} and Gaussian DP \citep{dong2022gaussian}, which are popular in practice due to flexible privacy mechanism design, better privacy budget accounting and increased statistical efficiency. Since these variants, and others such as subspace DP \citep{gao2022subspace}, stem from changing $\mathcal X, d$ or the distorting function $\dist$ (Remark~\ref{remarkDPCredalSet}), three key questions are 1) how the distorting function $\dist$ corresponding to a DP variant can be characterised as an IP object; 2) what IP properties does $\dist$ have; and 3) what are the consequences on statistical inference from using $\dist$ to constrain $P_x$ and $P_{x'}$ (for connected $x, x' \in \X$). %
For example, our preliminary analysis shows that $(\epsilon,\delta)$-DP cannot be described by an interval of measures, at least not alone. 
Moreover, we have demonstrated the necessity of using $\dMult$ as the distorting function for ensuring the types of bounds on frequentist and Bayesian inference found in Sections~\ref{sec:frequentist} and~\ref{sec:Bayesian} (see Remark~\ref{remarkUniqueness}). This implies other variants of DP which replace $\dMult$ with some other distorting function $\dist$ cannot satisfy the Bayesian and frequentist semantics described in this article.

A second direction for future research is suggested by the characterisation of Pufferfish as a Lipschitz condition with the metric $\dMult$ on $\Omega$ and the metric $\dPuff$ on $\Theta$. Replacing $\dMult$ with some other distorting function $\dist$ -- for example the distorting functions corresponding to $(\epsilon, \delta)$-DP or to zCDP (see \cite{bailieRefreshmentStirredNot2024a} for the definitions of these distorting functions) -- would generate new variants of Pufferfish. (Some such variants have already been proposed in \cite{zhangAttributePrivacyFramework2022, dingApproximationPufferfishPrivacy2024}.) Similarly, replacing $\dPuff$ with another metric on $\Theta$ would generalise Pufferfish beyond the interpretation of attackers and competing conjectures, and would provide a more-flexible framework for expressing privacy as $\epsilon$-indistinguishability (Remarks~\ref{remarkIndistinguishability} and~\ref{remarkPuffIndistinguish}) between distributions $\theta, \theta' \in \Theta$: By specifying their own metric on $\Theta$, the data custodian can choose which $\theta$ and $\theta'$ should be $\epsilon$-indistinguishable according to their knowledge and expertise, rather than being restricted to those $(\theta, \theta')$ pairs which correspond to $(\theta^*|_E, \theta^*|_{E'})$ for some $\theta^* \in \mathbb D$ and some $(E, E') \in \Spairs$.

Thirdly, a conceptually distinct IP approach to data privacy protection is the employment of SDL mechanisms that produce \emph{set-valued} outputs. This approach has not been explicitly considered in this paper, although one can take the elements $t$ of the output space $\T$ to themselves be sets, and as such, the results from this paper may still be applied. SDL mechanisms which produce set-valued outputs -- such as the ``leaky'' variant of Warner's randomised response considered in~\cite{li2022local} and the randomised censoring mechanisms considered in~\cite{dingIntervalPrivacyFramework2022} -- can be understood as an intentional \emph{coarsening} of the data product~\cite{heitjan1990inference}. A mechanism that produces set-valued outputs has a precise probability distribution that is given by the mass function associated with a \emph{belief function} -- a special type of coherent lower prevision \cite{shafer1976mathematical}. As such, one can view set-valued mechanisms as inducing imprecise probabilities on $\T$ -- rather than precise probabilities $P_x$ -- where this imprecise probability is a belief function. Compared to more general forms of IP, including the distortion models considered in this work, the mass function formulation of belief functions lends a computational advantage (particularly for Markov chain Monte Carlo -- see e.g.~\cite{jacob2021gibbs}). On the other hand, it is less clear that set-valued outputs are practically acceptable for many of the real-world use-cases of SDL. Data users may anticipate point-valued data in most situations, and may not be prepared to conduct further statistical processing of set-valued outputs. In sum, the utility of the set-valued approach to SDL remains open to formulation and assessment in future research.  %

\appendix

\section{Definition of \texorpdfstring{$\supp (x \mid t, \theta)$}{the Support of x}}\label{appDefnSupport}

We assume that $\X$ is equipped with a topology $\tau_{\X}$ and that $\mathscr G$ is the Borel $\sigma$-algebra induced by $\tau_{\X}$. Denote the support of $P_\theta$ by
\begin{equation}\label{eqDefnSupp}
	\supp(P_\theta) = \bigcap \{S \subset \X \mathrm{\ closed\ and\ Borel\ measurable} \mid P_\theta(S) = 1\}.
\end{equation}
Here we mean `closed' with respect to $\tau_{\X}$, not necessarily the topology induced by the metric $d$. (Generally we should not use the topology induced by $d$: Since $d$ is typically discrete, this topology results in $\supp (P_\theta) = \emptyset$ whenever $\mathcal X$ is uncountable\footnote{If $\X$ is uncountable, then $P_\theta(\{x\}) = 0$. Yet $\{x\}^c$ is closed under the discrete topology, and hence $x \notin \supp (P_\theta)$. This argument applies for all $x \in \X$; thus $\supp (P_\theta) = \emptyset$.} and then Theorems~\ref{thmMarginalData} and \ref{thmPosterior} would be vacuous.) A standard example would be $\mathcal X = \bigcup_{n \in \mathbb N} \mathbb R^n$ (i.e. the universe of datasets with a finite number of real-valued records) with the topology induced by the map 
\begin{align*}
    \mathcal X &\to \mathbb R^{\mathbb N} \\
    x &\mapsto (x, 0, 0, \ldots),
\end{align*}
when $\mathbb R^{\mathbb N}$ is equipped with the product Euclidean topology. (We assume that $\mathscr G$ is Borel to ensure that $\supp(P_\theta) \in \mathscr G$.)

Similarly, assume that $\T$ is equipped with a topology $\tau_{\T}$ and that $\F$ is the Borel $\sigma$-algebra induced by $\tau_{\T}$. Analogously to \eqref{eqDefnSupp}, define $\supp (P_x) \subset \mathcal T$ for each $x \in \X$. Write $\supp_0 (x \mid t) = \{ x \mid t \in \supp (P_x)\}$ and finally define
\[\supp(x \mid t, \theta) = \supp (P_\theta) \cap \supp_0 (x \mid t).\]
We assume that $(\X, \tau_{\X})$ and $(\T, \tau_{\T})$ are second countable (that is, there exist countable bases for $\tau_{\X}$ and $\tau_{\T}$) to ensure that 
\begin{equation}\label{eqZeroOutsideSupport}
    P_\theta(E_1) = P_x(E_2) = 0,
\end{equation}
for any measurable $E_1 \subset \supp (P_\theta)^c$ and $E_2 \subset \supp (P_x)^c$.

\section{The Density Ratio Metric Is Well-Defined}\label{appendixDRM}

\begin{proposition}\label{propDRMWell}
    The density ratio metric $\drm$ is well-defined. That is, $\drm(\mu,\nu)$ does not depend on the choice of $f,g$ and $\tau$ in~\eqref{eqDefnDRM}.
\end{proposition}

For $\mu, \nu \in \Omega$, write $\mu \ll \nu$ to denote that $\mu$ is absolutely continuous with respect to $\nu$.

We need the following result, which can be found in a standard probability-theory textbook (such as \cite[Exercise~32.6]{billingsleyProbabilityMeasure2012}): 

\begin{lemma}\label{lemmaLikelihoodRatioRN}
    Let $\mu, \nu, \tau \in \Omega$ such that $\mu \ll \nu$ and $\nu \ll \tau$. Then $\mu \ll \tau$ and the Radon-Nikodym derivative $\frac{d\mu}{d\nu}$ satisfies
    \begin{equation}\label{eqLemmaLikelihoodRatioRN}
        \frac{d\mu}{d\nu} = \frac{f}{g}, \quad \mu\text{-}\mathrm{a.e.},
    \end{equation}
    where $f, g$ are the $\tau$-densities of $\mu$ and $\nu$ respectively and, on the RHS of~\eqref{eqLemmaLikelihoodRatioRN}, $0/0 = 0$.
\end{lemma}

\begin{proof}of Proposition~\ref{propDRMWell}:
    Let $\tau_1, \tau_2 \in \Omega$ and suppose $\mu$ and $\nu$ are both non-zero and absolutely continuous with respect to both $\tau_1$ and $\tau_2$. Let $f_1$ and $f_2$ be the densities of $\mu$ with respect to $\tau_1$ and $\tau_2$ respectively. Similarly, define $g_1$ and $g_2$ as the densities of $\nu$ with respect to $\tau_1$ and $\tau_2$. 

    Define $\tau = \tau_1 + \tau_2$. Then $\tau_1 \ll \tau$ and $\tau_2 \ll \tau$. By Lemma~\ref{lemmaLikelihoodRatioRN},
    \[\frac{f_1}{g_1} = \frac{ \left. \frac{d\mu}{d\tau} \middle/ \frac{d\tau_1}{d\tau} \right.}{\left. \frac{d\nu}{d\tau} \middle/ \frac{d\tau_1}{d\tau} \right.} = \frac{\frac{d\mu}{d\tau}}{\frac{d\nu}{d\tau}},\quad \tau_1\text{-a.e.}\]
    Hence
    \begin{equation}\label{eqProofProp231}
        \tau_1\text{-}\esssup_{t, t' \in \mathcal T^{\mathrm{o}}} \frac{f_1(t)}{g_1(t)} \frac{g_1(t')}{f_1(t')} = \tau_1\text{-}\esssup_{t, t' \in \mathcal T^{\mathrm{o}}} \frac{\frac{d\mu}{d\tau}(t)}{\frac{d\nu}{d\tau}(t)} \frac{\frac{d\nu}{d\tau}(t')}{\frac{d\mu}{d\tau}(t')} \le \tau\text{-}\esssup_{t, t' \in \mathcal T^{\mathrm{o}}} \frac{\frac{d\mu}{d\tau}(t)}{\frac{d\nu}{d\tau}(t)} \frac{\frac{d\nu}{d\tau}(t')}{\frac{d\mu}{d\tau}(t')}
    \end{equation}
    (Note we use the notation $\tau\text{-}\esssup$ to refer to the essential supremum with respect to the measure $\tau$.) 
    
    Now we prove the reverse inequality of~\eqref{eqProofProp231}. For any $E \in \mathscr F$ with $\tau_1(E) = 0$, we have that
    \begin{equation}\label{eqProofProp232}
        \tau \left( \left\{ \frac{d\mu}{d\tau} > 0 \right\} \cap E \right) = 0.
    \end{equation}
    (Otherwise $\mu(E) = \int_E \frac{d\mu}{d\tau} d\tau > 0$ and hence $\mu \not\ll \tau_1$.) By symmetry, \eqref{eqProofProp232} also holds with $\frac{d\nu}{d\tau}$ in place of $\frac{d\mu}{d\tau}$. This implies 
    \begin{equation}\label{eqProofProp233}
        \tau_1\text{-}\esssup_{t, t' \in \mathcal T^{\mathrm{o}}} \frac{f_1(t)}{g_1(t)} \frac{g_1(t')}{f_1(t')} \ge \tau\text{-}\esssup_{t, t' \in \mathcal T^{\mathrm{o}}} \frac{\frac{d\mu}{d\tau}(t)}{\frac{d\nu}{d\tau}(t)} \frac{\frac{d\nu}{d\tau}(t')}{\frac{d\mu}{d\tau}(t')}.
    \end{equation}
    By combining~\eqref{eqProofProp232} and~\eqref{eqProofProp233}, we get an equality between these two essential suprema. By exactly the same reasoning, we have that
    \[\tau_2\text{-}\esssup_{t, t' \in \mathcal T^{\mathrm{o}}} \frac{f_2(t)}{g_2(t)} \frac{g_2(t')}{f_2(t')} = \tau\text{-}\esssup_{t, t' \in \mathcal T^{\mathrm{o}}} \frac{\frac{d\mu}{d\tau}(t)}{\frac{d\nu}{d\tau}(t)} \frac{\frac{d\nu}{d\tau}(t')}{\frac{d\mu}{d\tau}(t')},\]
    and hence
    \[\tau_1\text{-}\esssup_{t, t' \in \mathcal T^{\mathrm{o}}} \frac{f_1(t)}{g_1(t)} \frac{g_1(t')}{f_1(t')} = \tau_2\text{-}\esssup_{t, t' \in \mathcal T^{\mathrm{o}}} \frac{f_2(t)}{g_2(t)} \frac{g_2(t')}{f_2(t')}.\]
    Since logarithms are continuous, they are interchangeable with essential suprema:
    \begin{align*}
        \tau_1\text{-}\esssup_{t, t' \in \mathcal T^{\mathrm{o}}} \ln \left( \frac{f_1(t)}{g_1(t)} \frac{g_1(t')}{f_1(t')} \right) &= \ln \left( \tau_1\text{-}\esssup_{t, t' \in \mathcal T^{\mathrm{o}}} \frac{f_1(t)}{g_1(t)} \frac{g_1(t')}{f_1(t')} \right). %
    \end{align*}
    This proves that the value of $\drm(\mu, \nu)$ is the same when computed using $f_1, g_1$ and $\tau_1$, as when computed using $f_2, g_2$ and $\tau_2$.
\end{proof}

That the density ratio metric $\drm$ is well-defined (Proposition~\ref{propDRMWell}) is also an easy corollary of Proposition~\ref{propDRMDefn} below.

\section{Metric Spaces}\label{appendixMetric}

\begin{definition}\label{defnMetric}
    A \emph{metric} $d$ on a set $S$ is a function $S \times S \to [0,\infty]$ that satisfies the following properties for all $x,y,z \in S$:
    \begin{enumerate}
        \item Positive definiteness: $d(x,y) = 0$ if and only if $x = y$;\label{defnMetricProp1}
        \item Symmetry: $d(x,y) = d(y,x)$; and\label{defnMetricProp2}
        \item Triangle inequality: $d(x,z) \le d(x,y) + d(y,z)$.\label{defnMetricProp3}
    \end{enumerate}
    A \emph{premetric} $d$ on a set $S$ is a function $S \times S \to [0,\infty]$ that satisfies $d(x,x) = 0$ for every $x \in S$.
\end{definition}

Note that the co-domain of a metric is the extended, non-negative real numbers; we allow a metric to take the value of positive infinity. Metrics of this kind are sometimes referred to as \emph{extended}-metrics, or $\infty$-metrics, to distinguish them from those metrics with co-domain $[0,\infty)$.

Premetrics naturally arise in the context of differential privacy: The distorting function associated with a DP flavour's distortion model is a premetric \cite{bailieRefreshmentStirredNot2024a}. %
It is also common in many DP flavours that the `metric' $d$ on $\X$ is not in fact a metric, but only a premetric.

\begin{definition}\label{defnMetricStrongEquivalence}
    Two metrics $d_1$ and $d_2$ (which are defined on the same set $S$) are \emph{strongly equivalent} \citep[p.~121]{carothersRealAnalysis2000} if there exist constants $0 < a \le b < \infty$ such that
    \[a d_1(x, y) \le d_2(x, y) \le b d_1(x, y),\]
    for all $x, y \in S$.
\end{definition}

Proposition~\ref{propInequalityMultDRM} below proves that the multiplicative distance $\dMult$ and the density ratio metric $\drm$ are strongly equivalent on the space of probability measures (but not on the space $\Omega$ of $\sigma$-finite measures).

\begin{lemma}\label{lemmaMultMetric}
    The multiplicative distance $\dMult$ is a metric on the collection of measures on $(\mathcal T, \mathscr F)$.
\end{lemma}

\begin{proof}
    Suppose $\mu, \nu$ and $\tau$ are measures on $(\mathcal T, \mathscr F)$. To prove property \ref{defnMetricProp1}. of a metric, note that $\dMult(\mu, \nu) = 0$ if and only if $\mu(S) = \nu(S)$ for all $S \in \mathscr F$. Yet this holds if and only if $\mu = \nu$. Property \ref{defnMetricProp2}. follows by observing
    \[\dMult(\mu, \nu) = \sup_{S \in \mathscr F} \abs{\ln \frac{\mu(S)}{\nu(S)}} = \sup_{S \in \mathscr F} \abs{-\ln \frac{\mu(S)}{\nu(S)}} = \sup_{S \in \mathscr F} \abs{\ln \frac{\nu(S)}{\mu(S)}} = \dMult(\nu, \mu).\]
    Property \ref{defnMetricProp3}. is implied by
    \begin{align*}
        \dMult(\mu, \nu) &= \sup_{S \in \mathscr F} \abs{\ln \mu(S) - \ln \nu(S)} \\
        &= \sup_{S \in \mathscr F} \abs{\ln \mu(S) - \ln \tau(S) + \ln \tau (S) - \ln \nu(S)} \\
        &\le \sup_{S \in \mathscr F} \abs{\ln \mu(S) - \ln \tau(S)} + \sup_{S \in \mathscr F} \abs{\ln \tau (S) - \ln \nu(S)} \\
        & = \dMult(\mu, \nu) + \dMult(\tau, \nu).
    \end{align*}%
\end{proof}

\begin{lemma}\label{lemmaDRMMetric}
    The density ratio metric $\drm$ is a metric on the space $\Omega_1$ of probability measures on $(\mathcal T, \mathscr F)$.
\end{lemma}

\begin{proof}
    Suppose $P, Q$ and $R$ are probability measures on $(\mathcal T, \mathscr F)$. Proof of property \ref{defnMetricProp1}.: Suppose that $P = Q$. Then $P$ and $Q$ are mutually absolutely continuous and have $P$-densities $f$ and $g$ respectively, with 
    \[f(t) = g(t) = 1,\]
    for all $t \in \mathcal T$. Thus $\drm(P, Q) = 0$. Now suppose that $P \ne Q$. If $P$ and $Q$ are not mutually absolutely continuous, then $\drm(P, Q) = \infty > 0$. If they are mutually absolutely continuous, then they have $P$-densities $f$ and $g$ respectively. Moreover, because $P \ne Q$, %
    there exists $S \in \mathscr F$ such that $P(S) > Q(S)$. Hence $P(f > g) > 0$. This implies 
    \[ 0 < \esssup_{t \in \mathcal T^{\mathrm{o}}} \ln \frac{f(t)}{g(t)}.\]
    Also, $P(S^c) < Q(S^c)$ (here $S^c$ is the complement of $S$ in $\mathcal T$), so by exactly the same reasoning
    \[ 0 < \esssup_{t' \in \mathcal T^{\mathrm{o}}} \ln \frac{g(t')}{f(t')}.\]
    Combining these two results gives $\drm(P, Q) > 0$.

    Property \ref{defnMetricProp2}. follows by the symmetry of
    \[\left. \frac{f(t)}{f(t')} \middle/ \frac{g(t)}{g(t')}\right. = \left. \frac{g(t')}{g(t)} \middle/ \frac{f(t')}{f(t)}\right.,\]
    and the fact that $t$ and $t'$ are interchangeable in the definition of the density ratio metric.

    Finally, we prove property \ref{defnMetricProp3}. Suppose that $P$ and $R$ are not mutually absolutely continuous. Then $Q$ is either mutually absolutely continuous with $P$, or with $R$, but not both. Hence $\drm(P, Q) + \drm(Q, R) = \infty$ and thus 
    \[\drm(P, R) \le \drm(P, Q) + \drm(Q, R),\]
    holds vacuously. Now suppose that $P$ and $R$ are mutually absolutely continuous. We may also suppose that $Q$ is mutually absolutely continuous with respect to both $P$ and $R$. (When this is not the case, property \ref{defnMetricProp3}. again holds vacuously.) Let $f, g$ and $h$ be $P$-densities of $P, Q$ and $R$ respectively. Define 
    \begin{align*}
        \mathcal T_{f,g}^{\mathrm{o}} &= \{ t \in \mathcal T \mid 0 < f(t), g(t) < \infty \}, \\
        \mathcal T_{f,h}^{\mathrm{o}} &= \{ t \in \mathcal T \mid 0 < f(t), h(t) < \infty \}, \\
        \mathcal T_{g,h}^{\mathrm{o}} &= \{ t \in \mathcal T \mid 0 < g(t), h(t) < \infty \}.
    \end{align*}
    Then
    \begin{align*}
        \drm(P, R) &= \esssup_{t, t' \in \mathcal T_{f,h}^{\mathrm{o}}} \ln \left( \frac{f(t)}{f(t')}\right) - \ln \left( \frac{h(t)}{h(t')} \right) \\
        &= \esssup_{t, t' \in \mathcal T_{f,h}^{\mathrm{o}}} \ln \left( \frac{f(t)}{f(t')}\right) - \ln \left( \frac{g(t)}{g(t')} \right) + \ln \left( \frac{g(t)}{g(t')} \right) - \ln \left( \frac{h(t)}{h(t')} \right) \\
        &\le \esssup_{t, t' \in \mathcal T_{f,g}^{\mathrm{o}}} \ln \left( \frac{f(t)}{f(t')}\right) - \ln \left( \frac{g(t)}{g(t')} \right) + \esssup_{t, t' \in \mathcal T_{g,h}^{\mathrm{o}}} \ln \left( \frac{g(t)}{g(t')} \right) - \ln \left( \frac{h(t)}{h(t')} \right) \\
        &= \drm(P, Q) + \drm(Q, R),
    \end{align*}
    where all the essential suprema are with respect to $P$. We can exchange $\mathcal T_{f,g}^{\mathrm{o}}, \mathcal T_{f,h}^{\mathrm{o}}$ and $\mathcal T_{g,h}^{\mathrm{o}}$ in the above computations because 
    \[P\left(\mathcal T_{f,g}^{\mathrm{o}}\right) = P\left(\mathcal T_{f,h}^{\mathrm{o}}\right) = P\left(\mathcal T_{g,h}^{\mathrm{o}}\right) = 1.\]
\end{proof}

The density ratio metric $\drm$ is not a metric on $\Omega$. While it is easy to verify that properties \ref{defnMetricProp2}. and \ref{defnMetricProp3}. hold (follow the same reasoning as in the above proof), $\drm$ does not satisfy property \ref{defnMetricProp1}: There exists $\mu \ne \nu \in \Omega$ such that $\drm(\mu, \nu) = 0$. For example, suppose that 
\begin{equation}\label{eqDRMCounterexample}
    \nu(S) = 2 \mu(S),    
\end{equation}
for all $S \in \mathscr F$. Then $f(t) = 1$ and $g(t) = 2$ are $\mu$-densities of $\mu$ and $\nu$ respectively. Yet this implies $\drm(\mu, \nu) = 0$.

\begin{proposition}\label{propDRMPseudoMetric}
    The density ratio metric $\drm$ is a pseudo-metric on $\Omega$. That is, $\drm$ is a function from $\Omega \times \Omega$ to the extended real line which satisfies the following properties for all $\mu, \nu, \tau \in \Omega$:
    \begin{enumerate}
        \item $\drm(\mu, \mu) = 0$;\label{defnPseudoMetricProp1}
        \item $\drm(\mu, \nu) \ge 0$;\label{defnPseudoMetricProp2}
        \item $\drm(\mu, \nu) = \drm(\nu, \mu)$; and\label{defnPseudoMetricProp3}
        \item $\drm(\mu, \tau) \le \drm(\mu, \nu) + \drm (\nu, \tau)$.\label{defnPseudoMetricProp4}
    \end{enumerate}
\end{proposition}

\begin{proof}
    Properties \ref{defnPseudoMetricProp1}., \ref{defnPseudoMetricProp3}. and \ref{defnPseudoMetricProp4}. follow by the same reasoning as in the proof of Lemma~\ref{lemmaDRMMetric} with minor adjustments to account for the fact that $\mu, \nu$ or $\tau$ may be zero. %
    Property \ref{defnPseudoMetricProp2}. follows by applying Properties \ref{defnPseudoMetricProp1}., \ref{defnPseudoMetricProp4}. and \ref{defnPseudoMetricProp3}. sequentially:
    \begin{equation*}
        0 = \drm (\mu, \mu) \le  \drm(\mu, \nu) + \drm(\nu, \mu) = 2\drm(\mu, \nu).%
    \end{equation*}
\end{proof}

\section{Distorting Functions and Distortion Models}\label{appDistort}

\begin{definition}
    \hspace{-0.2em}{\normalfont \textbf{\cite{montesUnifyingNeighbourhoodDistortion2020}}} A \emph{distorting function} $\dist$ is a function $S \times S \to [0, \infty]$ where $S$ is one of the following spaces:
    \begin{itemize}
        \item the set of all measures on a measurable space $(\T, \mathscr F)$;
        \item the set $\Omega$ of all $\sigma$-finite measures on $(\T, \mathscr F)$;
        \item the set of all finite measures on $(\T, \mathscr F)$;
        \item the set $\Omega_1$ of all probability measures on $(\T, \mathscr F)$; or
        \item the set $\Omega_1^* = \{ P \in \Omega_1 \mid P(\{t\}) > 0\ \forall t \in \T\}$ of probabilities on $(\T, 2^{\T})$ with non-zero mass on every event $E \subset \T$ (often with the assumption that $\T$ has finite cardinality).
    \end{itemize}
\end{definition}

\begin{definition} 
    \hspace{-0.2em}{\normalfont \textbf{\cite{montesUnifyingNeighbourhoodDistortion2020}}} Given a distorting function $\dist$, a positive constant $r > 0$ (termed the \emph{distortion parameter}), and a probability $P_0 \in \Omega_1$ (the \emph{nucleus}), the distortion model on $P_0$ associated with $\dist$ and $r$ is the closed $\dist$-ball centred at $P_0$ with radius $r$:
    \[B_{\dist}^r(P_0) = \{ P \in \Omega_1 \mid \dist(P, P_0) \le r \}.\]
\end{definition}

One may use $\Omega_1^*$ in place of $\Omega_1$ in the above definition, as in \cite{montesUnifyingNeighbourhoodDistortion2020}.

A distortion model is one example of the more general notion of neighbourhood models, found throughout the literature on IP and robustness. In general, a neighbourhood model on $P_0$ is simply a set of probabilities which contains $P_0$. An advantage of distortion models is that they are characterised by a small number of parameters (three). Beyond the symmetric IoM $\IIntP{e^{-\epsilon}P_0, e^{\epsilon} P_0}$ and the density ratio neighbourhood $N_r(P_0) \cap \Omega_1$, other examples of distortion models can be found in \cite{montesUnifyingNeighbourhoodDistortion2020, montesUnifyingNeighbourhoodDistortion2020a, montesNeighbourhoodModelsInduced2023, mirandaEvaluatingUncertaintyVertical2024, desterckeProcessingDistortionModels2022, pelessoniInferenceNearlyLinearUncertainty2021, walley1991statistical}.

\sepfootnotecontent{footnoteContinuity1}{Note that our definition of continuity is strictly stronger than that given in \cite{montesUnifyingNeighbourhoodDistortion2020} which allows $\delta$ to depend on $\nu_3$.}

\sepfootnotecontent{footnoteContinuity2}{Given a $\sigma$-finite measure space $(\T, \F, \mu)$, the supremum norm $\norm{\cdot}_\infty^\mu$ on the set $\Omega_{\nu}$ is defined as:
\[\norm{\nu}_\infty^\mu = \esssup \abs{f(t)},\]
where the essential supremum is with respect to $\mu$ and $f$ is a $\mu$-density of $\nu$.}

Invariance under marginalisation, invariance under updating and immunity to dilation are three desiderata for distortion models. For symmetric IoMs and density ratio neighbourhoods, these desiderata are studied in \cite{wassermanComputingBoundsExpectations1992, wasserman1992invariance, seidenfeld1993dilation}. %
\cite{decamposPROBABILITYINTERVALSTOOL1994} showed that, under certain restrictions, $\IIntP{e^{-\epsilon} P_0, e^{\epsilon} P_0}$ is 2-monotone. To the best of our knowledge, we are not aware of studies investigating the geometry of the symmetric IoM or the density ratio neighbourhood -- for example, are they polytopes and, if so, how many extremal points do they have?

In Appendix~\ref{appendixMetric}, we showed that the distorting functions $\dMult$ and $\drm$ are metrics -- that is, they are positive definite, symmetric and satisfy the triangle inequality on the space $\Omega_1$ of probability measures. Two other desirable properties of a distorting function $\dist$ are: 
\begin{enumerate}
    \item \emph{Quasi-Convexity} \citep{montesUnifyingNeighbourhoodDistortion2020, mirandaEvaluatingUncertaintyVertical2024}: Given a real- or complex-vector space $V$, a function $d : V \times V \to [0,\infty]$ is quasi-convex (in its second argument) if 
    \begin{equation}\label{eqOurDefnConvex}
        d\left( v_1, \alpha v_2 + [1-\alpha] v_3 \right) \le \max \{ d(v_1, v_2), d(v_1, v_3) \},
    \end{equation}
    for all $\alpha \in [0,1]$ and all $v_1, v_2, v_3 \in V$.%
    \item \emph{Continuity} \citep{montesUnifyingNeighbourhoodDistortion2020}:\sepfootnote{footnoteContinuity1} For $\mu \in \Omega$, let $\Omega_{\nu} = \{ \nu \in \Omega \mid \nu \ll \mu\}$ be the set of $\sigma$-finite measures $\nu$ which are absolutely continuous with respect to $\mu$. Given $S \subset \Omega$, a function $d : S \times S \to [0,\infty]$ is continuous (in its second argument) with respect to the supremum norm if, for all $\mu \in \Omega$, all $\nu_1, \nu_2 \in S \cap \Omega_{\mu}$ and all $\epsilon > 0$, there exists $\delta > 0$ such that, for all $\nu_3 \in S \cap \Omega_{\mu}$,
    \[\norm{\nu_2 - \nu_3}_\infty^\mu < \delta \Rightarrow \abs{d (\nu_1, \nu_2) - d(\nu_1, \nu_3)} < \epsilon,\]
    where $\norm{\cdot}_{\infty}^\mu$ denotes the supremum norm.\sepfootnote{footnoteContinuity2} (The continuity is uniform if $\delta$ does not depend on $\nu_1$ or $\nu_2$.)
\end{enumerate}

\begin{remark}\label{remarkConvex}
    The astute reader may consider \eqref{eqOurDefnConvex} to be a strange definition of convexity. A more standard definition of convexity (in the second argument) would be the requirement:
    \begin{equation}\label{eqStdDefnConvex}
        d\left( v_1, \alpha v_2 + [1-\alpha] v_3 \right) \le \alpha d(v_1, v_2) + (1-\alpha) d(v_1, v_3),
    \end{equation}
    for all $\alpha \in [0,1]$ and all $v_1, v_2, v_3 \in V$. This requirement is strictly stronger than %
    quasi-convexity (equation~\eqref{eqOurDefnConvex}). (It is straightforward to prove \eqref{eqStdDefnConvex} is stronger than \eqref{eqOurDefnConvex} and we provide some examples later in this remark to prove strictness.) In fact, \eqref{eqStdDefnConvex} is often too strong a requirement, for three reasons.
    
    Firstly, if a distorting function $\dist$ on $\Omega^*_1$ satisfies \eqref{eqOurDefnConvex} and continuity, then the ball $\{P \in \Omega^*_1 \mid \dist (P, P_0) \le r\}$ centred at $P_0 \in \Omega^*_1$ is equal to the credal set induced by the lower envelope of this ball \citep[Proposition~3.1]{montesUnifyingNeighbourhoodDistortion2020}. Further, under \eqref{eqOurDefnConvex} and continuity, there exist simple necessary and sufficient conditions for this lower envelope to be a probability interval \citep[Propositions~3.3, 3.4]{montesUnifyingNeighbourhoodDistortion2020}. Hence, \eqref{eqOurDefnConvex} is a useful requirement for a distorting function $\dist$ as it ensures the distortion model associated with $\dist$ is well-behaved. It would be unnecessary to require the stricter condition~\eqref{eqStdDefnConvex} solely to ensure this nice behaviour.
    
    Secondly, some of the common distorting functions found in the IP literature do not satisfy \eqref{eqStdDefnConvex} but do satisfy \eqref{eqOurDefnConvex}. One example is the linear vacuous distorting function:
    \begin{align*}
        \dLV : \Omega_1^* \times \Omega_1^* &\to [0, \infty), \\
        (P, Q) & \mapsto \max_{ \emptyset \ne S \subset \T} \frac{Q(S) - P(S)}{Q(S)}.
    \end{align*}
    \citet[Proposition~5.1]{montesUnifyingNeighbourhoodDistortion2020} show that $\dLV$ satisfies \eqref{eqOurDefnConvex}. However, the following counterexample demonstrates that $\dLV$ does not satisfy \eqref{eqStdDefnConvex}: Let $\alpha = 0.5$ and $\T = \{1,2\}$. Let $P_1$ and $P_2$ be the uniform probability on $(\T, 2^\T)$ and define $P_3 \in \Omega_1^*$ by $P_3(\{1\}) = 0.7$. Then
    \[ \dLV (P_1, \alpha P_2 + (1-\alpha) P_3) = \tfrac{1}{6} > \tfrac{1}{7} = \alpha \dLV(P_1, P_2) + (1-\alpha) \dLV (P_1, P_3).\]
    The multiplicative distance $\dMult$ also does not satisfy \eqref{eqStdDefnConvex} (although we show below that it satisfies \eqref{eqOurDefnConvex}), even when restricting to $\Omega_1^*$ and to $\T$ with finite cardinality. To see this, let $\alpha = 0.5$ and let $P_1$ be uniform on $\T = \{1,\ldots,10\}$; define $P_2$ by $P_2(\{1\}) = 1/5$ and $P_2(\{t\}) = 8/90$ for all $t \ne 1$; and define $P_3$ by $P_3(\{1\}) = 3/10$ and  $P_3(\{t\}) = 7/90$ for all $t \ne 1$. Then 
    \[\dMult(P_1, \alpha P_2 + (1-\alpha) P_3) = \ln 2.5 > \tfrac{1}{2} (\ln 2 + \ln 3) = \alpha \dMult(P_1, P_2) + (1-\alpha) \dMult(P_1, P_3).\]

    Thirdly and finally, \eqref{eqOurDefnConvex} aligns with the notion of convexity which arises independently in the DP literature \cite{kifer2022bayesian}: Under a mild assumption, a DP flavor \cite{bailieRefreshmentStirredNot2024a} with distorting function $\dist$ is convex in the sense of \cite{kifer2022bayesian} if and only if $\dist$ is quasi-convex in the sense of %
    \eqref{eqOurDefnConvex}. Because most of the commonly used DP flavors are convex, it follows that most of the distorting functions used in DP are quasi-convex.
\end{remark}

\begin{proof}that $\dMult$ satisfies quasi-convexity: 
    Let $\mu, \nu$ and $\tau$ be measures on $(\T, \F)$. For $E \in \F$, we have
    \begin{align*}
        \ln \frac{\mu(E)}{\alpha \nu (E) + (1-\alpha) \tau (E)} &\le \ln \frac{\mu(E)}{\min \{ \nu(E), \tau(E)\} } \\
        &\le \sup_{S \in \F} \left\{ \max \left\{ \ln \frac{\mu(S)}{\nu(S)}, \ln \frac{\mu(S)}{\tau(S)} \right\} \right\} \\
        &\le \max \{ \dMult(\mu, \nu), \dMult(\mu, \tau) \}.
    \end{align*} 
    Similarly,
    \begin{align*}
        \ln \frac{\alpha \nu (E) + (1-\alpha) \tau(E)}{\mu(E)} &\le \ln \frac{\max \{ \nu(E), \tau(E) \}}{\mu(E)} \\
        &\le \sup_{S \in \F} \left\{ \max \left\{ \ln \frac{\nu(S)}{\mu(S)}, \ln \frac{\tau(S)}{\mu(S)} \right\} \right\} \\
        &\le \max \{ \dMult(\mu, \nu), \dMult(\mu, \tau) \}.
    \end{align*}%
\end{proof}

\begin{lemma}\label{lemmaCtsSuffCond}
    Let $d : S \times S \to [0,\infty]$ be a metric (where $S \subset \Omega$). The following is a sufficient condition for $d$ to be continuous: For all $\mu \in \Omega$, all $\nu_1 \in S \cap \Omega_{\mu}$ and all $\epsilon > 0$, there exists $\delta > 0$ such that, for all $\nu_2 \in S \cap \Omega_{\mu}$,
    \[\norm{\nu_1 - \nu_2}_\infty^\mu < \delta \Rightarrow d (\nu_1, \nu_2) < \epsilon.\]
\end{lemma}

\begin{proof}
    This result follows by the triangle inequality and symmetry: $\abs{ d(\nu_1, \nu_2) - d(\nu_1, \nu_3)} \le d(\nu_2, \nu_3)$.
\end{proof}

The following proposition proves that $\dMult$ is continuous, in the setting considered in \cite{montesUnifyingNeighbourhoodDistortion2020}.

\begin{proposition}\label{propDMultCont}
    Suppose that $\T$ has finite cardinality. Let $\mathscr F = 2^{\T}$ and $\Omega^* = \{ \nu \in \Omega \mid \nu(\{t\}) > 0\ \forall t \in \T\}$. Then $\dMult$ is continuous when restricted to the domain $\Omega^* \times \Omega^*$.
\end{proposition}

Proposition~\ref{propDMultCont} implies that the interval of measures $\IInt {L, U} \cap \Omega^*$ is closed with respect to the topology induced by the supremum norm, as long as $\T$ has finite cardinality.

\begin{proof}
    Lemma~\ref{lemmaCtsSuffCond} describes a sufficient condition for continuity. We will prove that this sufficient condition holds under the assumption that $\T$ has finite cardinality and under the restriction to $\sigma$-finite measures $\nu$ with $\nu(\{t\}) > 0$. 
    
    Without loss of generality, we assume that the dominating measure $\mu$ is the counting measure on $\T$. Then $\Omega^* \subset \Omega_{\mu}$. Let $\nu_1 \in \Omega^*$ and $\epsilon > 0$. Let $f$ be the $\mu$-density of $\nu_1$. Choose some $\delta$ which satisfies
    \[ 0 < \delta < \min_{t \in \T} f(t) (1 - e^{-\epsilon}).\]
    (Such a $\delta$ exists because $f(t) > 0$ for all $t \in \T$ and because $\T$ has finite cardinality.) Fix some $\nu_2 \in \Omega^*$ with $\norm{ \nu_1 - \nu_2}_\infty^{\mu} < \delta$. Then, for any $S \ne \emptyset$,
    \[\frac{\nu_2(S)}{\nu_1(S)} \le \frac{\nu_1(S) + \delta}{\nu_1(S)} \le \max_{t \in \T} \frac{f(t) + \delta}{f(t)} \le \max_{t \in \T} \frac{f(t)}{f(t)-\delta} < e^\epsilon.\]
    (The last inequality follows because $\delta < f(t) - f(t) e^{-\epsilon}$ for all $t$.) Similarly,
    \[\frac{\nu_1(S)}{\nu_2(S)} \le \frac{\nu_1(S)}{\nu_1(S) - \delta} \le \max_{t \in \T} \frac{f(t)}{f(t)-\delta} < e^\epsilon.\]
    Hence $\dMult(\nu_1, \nu_2) < \epsilon$.
\end{proof}

If we remove either of the two restriction in Proposition~\ref{propDMultCont} (to finite $\T$ and to $\mu, \nu \in \Omega^*$), then $\dMult$ is no longer continuous. We will demonstrate that both of these restrictions are necessary with two counterexamples: In the first case, suppose $\T = \{1,2, \ldots\}$ with $\mathscr F = 2^{\T}$ and let $\mu$ be the counting measure. Define $\nu_1(\{t\}) = t^{-2}$. Let $\epsilon = \ln 2$ and fix $0 < \delta < 1$. Define
\[\nu_2(\{t\}) = \begin{cases}
    1-\delta/4 & \mathrm{if\ } t = 1,\\
    t^{-2} + \delta/4 &\mathrm{if\ } t = t_0, \\
    t^{-2} &\mathrm{otherwise,}
\end{cases}\]
where $t_0 = \ceil{\sqrt{4 \delta^{-1}} + 1}$. Then $\norm{\nu_1 - \nu_2}_{\infty}^\mu < \delta$ but 
\[\dMult(\nu_1, \nu_2) \ge \ln \frac{t_0^{-2} + \delta/4}{t_0^{-2}} > \ln 2.\]
In the second case, let $\T$ be finite, $\F = 2^\T$ and $\mu$ be the counting measure. For some $t_0 \in \T$, define $\nu_1(\{t_0\}) = \delta/2$ and $\nu_2(\{t_0\}) = 0$, and suppose that $\nu_1(\{t\}) = \nu_2(\{t\})$ for all $t \ne t_0$. Then $\norm{\nu_1 - \nu_2}_{\infty}^\mu < \delta$ but $\dMult(\nu_1, \nu_2) = \infty$.

These two counterexamples can easily be modified so that $\nu_1$ and $\nu_2$ are probability measures. Hence $\dMult$ is also not continuous on the space of probability measures, except when $\T$ is finite and we restrict to probabilities with support $\T$.

\section{Supplementary Results}\label{appendixAdditionalResults}

\subsection{An Equivalent Definition of the Multiplicative Distance\texorpdfstring{ $\dMult$}{}}

\begin{lemma}\label{lemmaDensityFiniteSigma}
    Suppose that $\mu \in \Omega$. Let $f$ be a $\mu$-density of $\nu$. Then $\nu$ is $\sigma$-finite if and only if $f$ is finite $\mu$-almost everywhere.
\end{lemma}

\begin{proof}
    ``$\Rightarrow$'' by the contrapositive: (This direction does not require that $\mu$ is $\sigma$-finite.) Suppose that $\mu ( f = \infty) > 0$. Let $\{E_n : n \in \mathbb N\} \subset \F$ be any countable partition of $\mathcal T$. We will show that necessarily there exists some $E_n$ with $\nu(E_n) = \infty$. Since
    \[0 < \mu(f = \infty) = \sum_{n=1}^{\infty} \mu(E_n \cap \{ f = \infty\}),\]
    there exists some $E_n$ with $\mu(E_n \cap \{ f = \infty \} ) > 0$. Then
    \[\nu(E_n) \ge \nu (E_n \cap \{f = \infty\}) = \int_{E_n \cap \{f = \infty\}} f d\mu = \infty.\]

``$\Leftarrow$'': (This direction requires that $\mu$ is $\sigma$-finite -- the case that $\mu = \nu$ with $f = 1$ serves as a counterexample.) Let $A = \{f = \infty\} \in \F$. Then $\nu(A) = \mu(A) = 0$ because $\nu$ is absolutely continuous with respect to $\mu$ by assumption. Define $E_n = \{ t \in \T \setminus A : n-1 \le f(t) < n\}$ for all $n \in \mathbb N$. Let $\{S_m : m \in \mathbb N\}$ be a partition of $\T$ such that $\mu(S_m)<\infty$ for all $m \in \mathbb N$. Then $\nu(E_n \cap S_m) < n \mu(S_m) < \infty$. Hence $\{A, E_n \cap S_m : n, m \in \mathbb N\}$ is a countable partition of $\T$ with each component having finite $\nu$-measure.
\end{proof}

The following proposition gives an alternative definition of the multiplicative distance $\dMult$.

\begin{proposition}\label{propMultDistanceDefn}
        Let $\mu, \nu \in \Omega$. Then
        \[\dMult(\mu, \nu) = \begin{cases}
            0 & \mathrm{if\ } \mu = \nu = 0,\\
            \esssup_{t \in \mathcal T^{\mathrm{o}}} \abs{ \ln \frac{f(t)}{g(t)} } &\mathrm{else\ if\ } \mu,\nu \mathrm{\ are\ mutually\ absolutely\ continuous},\\
            \infty &\mathrm{otherwise},
        \end{cases}\]
        where $f$ and $g$ are any densities of $\mu$ and $\nu$ respectively, with respect to any common dominating measure $\tau \in \Omega$; $\mathcal T^{\mathrm{o}} = \{ t \in \mathcal T \mid 0 < f(t), g(t) < \infty \}$; and the essential supremum is with respect to $\tau$.
        
\end{proposition}

To prove the above proposition, we need two lemmata.

\begin{lemma}\label{lemmaEquivalenceBetweenSupEssSup}
    Suppose that $\mu, \nu \in \Omega$ are non-zero and mutually absolutely continuous and $\tau \in \Omega$ is a dominating measure. Let $f$ and $g$ be $\tau$-densities of $\mu$ and $\nu$ respectively. 

    Then 
    \[\sup_{S \in \mathscr F} \ln \frac{\mu(S)}{\nu(S)} = \esssup_{t \in \To} \ln \frac{f(t)}{g(t)},\]
    where the essential supremum is over $\tau$.
\end{lemma}

\begin{proof}
    We first need that (A) $f, g < \infty$ holds $\tau$-almost everywhere and that (B) $f, g > 0$ holds $\mu$- and $\nu$-almost everywhere. (A) is a direct result of Lemma~\ref{lemmaDensityFiniteSigma}.
    To prove (B),  observe that
    \[\mu(f = 0) = \int_{\{f = 0\}} f d\tau = 0,\]
    and then $\nu(f = 0) = 0$ follows by absolute continuity. That $g > 0$ holds $\mu$- and $\nu$-almost everywhere has a similar proof.

    For $a \in \mathbb R$, define $E_a = \{ t \in \mathcal T^{\mathrm{o}} \mid f(t) > \exp(a) g(t) \}$.
    We need to prove the following result holds for all $a \in \mathbb R$: There exists $S \in \mathscr F$ such that $\mu(S) > \exp(a) \nu(S)$ if and only if $\tau (E_a) > 0$. Denote this result by ($*$).
    
    Suppose that there exists $S \in \mathscr F$ such that $\mu(S) > \exp(a) \nu(S)$. By (A) and (B), this implies 
    \[\mu(S \cap \To) > \exp (a) \nu (S \cap \To),\]
    and hence
    \[\int_{S \cap \To} (f - e^a g) d\tau > 0.\]
    Thus, there exists some $E \subset S \cap \To$ such that $\tau(E) > 0$ and $f(t) - \exp (a) g(t) > 0$ for all $t \in E$. Hence $\tau(E_a) \ge \tau(E) > 0$. 

    In the other direction, suppose $\tau(E_a) > 0$. Then 
    \[\int_{E_a} (f - e^a g) d\tau > 0,\]
    which implies $\mu(E_a) > \exp (a) \nu(E_a)$. This proves ($*$).

    Finally, we have
    \begin{align*}
        \sup_{S \in \mathscr F} \ln \frac{\mu(S)}{\nu(S)} &= \sup \{ a \in \mathbb R \mid \text{there exists } S \in \mathscr F \text{ s.t. } \mu(S) > \exp (a) \nu (S) \} \\
        &= \sup \{ a \in \mathbb R \mid \tau ( E_a ) > 0 \} \\
        &= \esssup_{t \in \To} \ln \frac{f(t)}{g(t)},
    \end{align*}
    where the first line follows by continuity of $\exp(\cdot)$; the second by ($*$); and the third by the definition of the essential supremum.
\end{proof}

\begin{lemma}\label{lemmaNotAbsCts}
    Let $\mu, \nu \in \Omega$ be non-zero and not mutually absolutely continuous. Then
    \[\sup_{S \in \mathscr F} \ln \frac{\mu(S)}{\nu(S)} = \infty \quad\mathrm{or}\quad \sup_{S \in \mathscr F} \ln \frac{\nu(S)}{\mu(S)} = \infty.\]
\end{lemma}

\begin{proof}
    Without loss of generality, there exists some $E \in \mathscr F$ such that $\mu(E) > 0$ but $\nu(E) = 0$. Then 
    \[\sup_{S \in \mathscr F} \ln \frac{\mu(S)}{\nu(S)} \ge \ln \frac{\mu(E)}{\nu(E)} = \infty.\]
\end{proof}

\begin{proof}of Proposition~\ref{propMultDistanceDefn}:
    Suppose that $\mu = \nu = 0$. Then, for all $S \in \mathscr F$, we have that $\ln \frac{\mu(S)}{\nu(S)} = \ln \frac{0}{0} = 0$ where we define $0/0 = 1$ as in the definition of the multiplicative distance (Definition~\ref{defnDP}). Hence $\dMult(\mu, \nu) = 0$.

    Now suppose that $\mu$ and $\nu$ are both non-zero and not mutually absolutely continuous. By Lemma~\ref{lemmaNotAbsCts}, $\dMult(\mu, \nu) = \infty$.%

    Finally, suppose that $\mu$ and $\nu$ are both non-zero and mutually absolutely continuous. Fix a dominating measure $\tau$ and $\tau$-densities $f$ and $g$. 
    By symmetry, it suffices to show that 
    \[ \esssup_{t \in \mathcal T^{\mathrm{o}}} \ln \frac{f(t)}{g(t)} = \sup_{S \in \mathscr F} \ln \frac{\mu(S)}{\nu(S)}.\]
    This equation is given by Lemma~\ref{lemmaEquivalenceBetweenSupEssSup}.
\end{proof}

\subsection{The Multiplicative Distance \texorpdfstring{$\dMult$ }{} Is Composable}

\begin{proposition}\label{propMultDistComp}
    Let $(\T_1, \F_1)$ and $(\T_2, \F_2)$ be measurable spaces and let $\Omega^{(i)}$ be the collection of $\sigma$-finite measures on $(\T_i, \F_i)$. Then, for all $\mu_1, \nu_1 \in \Omega^{(1)}$ and all $\mu_2, \nu_2 \in \Omega^{(2)}$,
    \begin{equation}\label{eqPropMultDistComp}
        \dMult \big( \mu_1 \times \mu_2, \nu_1 \times \nu_2 \big) \le \dMult (\mu_1, \nu_1) + \dMult(\mu_2, \nu_2),
    \end{equation}
    where $\mu_1 \times \mu_2$ and $\nu_1 \times \nu_2$ are product measures.
\end{proposition}

\begin{proof}
    We will use the definition of the multiplicative distance $\dMult$ given in Proposition~\ref{propMultDistanceDefn}. 

    Suppose that $\mu_1$ and $\nu_1$ are zero. Then $\mu_1 \times \mu_2$ and $\nu_1 \times \nu_2$ are also zero. Hence \eqref{eqPropMultDistComp} simplifies to $0 \le \dMult(\mu_2, \nu_2)$, which always holds because $\dMult$ is a metric.

    Suppose that $\mu_1$ and $\nu_2$ are non-zero but not absolutely continuous. Then the RHS of \eqref{eqPropMultDistComp} is infinite, and so \eqref{eqPropMultDistComp} holds vacuously.

    Suppose that $\mu_1$ and $\nu_1$ are non-zero and mutually absolutely continuous. Then there exists a common dominating measure $\tau \in \Omega_1$. Let $f$ and $g$ denote $\tau$-densities of $\mu_1$ and $\nu_1$ respectively. Then, Lemma~\ref{lemmaDensityFiniteSigma} and Proposition~\ref{propMultDistanceDefn} together imply that 
    \begin{equation}\label{eqProofpropMultDistComp}
        f(t) \le g(t) \exp \left[ \dMult (\mu_1, \nu_1) \right],
    \end{equation}
    for $\tau$-a.e. $t \in \T_1$. Thus, for $S \in \F_1 \otimes \F_2$,
    \begin{align*}
        (\mu_1 \times \mu_2)(S) &= \int_{\T_1} \mu_2(S^t) d\mu_1(t) \\
        &= \int_{\T_1} \mu_2(S^t) f(t) d\tau (t) \\
        &\le \int_{\T_1} \mu_2(S^t) g(t) \exp \left[ \dMult(\mu_1, \nu_1) \right] d\tau (t) \\
        &\le \int_{\T_1} \exp \left[ \dMult(\mu_2, \nu_2) \right] \nu_2(S^t) g(t) \exp \left[ \dMult(\mu_1, \nu_1) \right] d\tau (t) \\
        &= \exp \left[ \dMult(\mu_1, \nu_1) + \dMult(\mu_2, \nu_2) \right] \int_{\T_1} \nu_2(S^t) g(t) d\tau(t) \\
        &= \exp \left[ \dMult(\mu_1, \nu_1) + \dMult(\mu_2, \nu_2) \right]\ (\nu_1 \times \nu_2) (S),
    \end{align*}
    where the third line follows by~\eqref{eqProofpropMultDistComp} and the fourth by the definition of $\dMult$ (Definition~\ref{defnDP}). Hence
    \[ \ln \frac{(\mu_1 \times \mu_2)(S)}{(\nu_1 \times \nu_2) (S)} \le \dMult(\mu_1, \nu_1) + \dMult(\mu_2, \nu_2),\]
    for all $S \in \F_1 \otimes \F_2$. The proof of the bound 
    \[ \ln \frac{(\nu_1 \times \nu_2)(S)}{(\mu_1 \times \mu_2) (S)} \le \dMult(\mu_1, \nu_1) + \dMult(\mu_2, \nu_2),\]
    follows analogously.
\end{proof}

\subsection{An Equivalent Definition of the Density Ratio Metric\texorpdfstring{ $\drm$}{}}

The following proposition gives an alternative definition of the density ratio metric $\drm$. 

\begin{proposition}\label{propDRMDefn}
    For any non-zero $\mu, \nu \in \Omega$,
    \begin{equation}\label{eqPropDRMDefn}
        \drm(\mu, \nu) = \sup_{S \in \mathscr F} \ln \frac{\mu(S)}{\nu(S)} + \sup_{S \in \mathscr F} \ln \frac{\nu(S)}{\mu(S)},
    \end{equation}
    where $0/0 = \infty/\infty = 1$. Also, $\drm(\mu, \nu) = \infty$ if exactly one of $\mu$ or $\nu$ is zero; and $\drm(\mu, \nu) = 0$ if $\mu$ and $\nu$ are both zero.
\end{proposition}

Compare the two definitions of the multiplicative distance $\dMult$ -- Definition~\ref{defnDP} and Proposition~\ref{propMultDistanceDefn} -- with the two definitions of the density ratio metric $\drm$ -- Definition~\ref{defnDensityRatio} and Proposition~\ref{propDRMDefn}. Each of these two distances have definitions in terms of densities and in terms of measures.

\begin{proof}of Proposition~\ref{propDRMDefn}:
    Firstly, we must verify that the RHS of~\eqref{eqPropDRMDefn} is well defined -- i.e. that the RHS cannot take on the form of $\infty-\infty$. We need to prove that ($*$) both $\sup_{S \in \mathscr F} \ln \frac{\mu(S)}{\nu(S)}$ and $\sup_{S \in \mathscr F} \ln \frac{\nu(S)}{\mu(S)}$ are bounded away from negative infinity.

    Let $E_1, E_2, \ldots$ be a partition of $\T$ such that $\mu(E_n) < \infty$. Since $\nu(\mathcal T) > 0$ there must exist some $E_n$ with $\nu(E_n) > 0$. Then
    \[\sup_{S \in \mathscr F} \ln \frac{\nu(S)}{\mu(S)} \ge \ln \frac{\nu(E_n)}{\mu(E_n)} > -\infty.\]
    The proof of $\sup_{S \in \mathscr F} \ln \frac{\mu(S)}{\nu(S)} > -\infty$ is analogous. 

    Suppose that $\mu = \nu = 0$. Then 
    \[\ln \frac{\mu(S)}{\nu(S)} = \ln \frac{\nu(S)}{\mu(S)} = 0,\]
    for all $S \in \mathscr F$. Hence~\eqref{eqPropDRMDefn} holds.
    
    Suppose that $\mu$ and $\nu$ are non-zero and not mutually absolutely continuous. Then by Lemma~\ref{lemmaNotAbsCts} and ($*$), 
    \[\sup_{S \in \mathscr F} \ln \frac{\mu(S)}{\nu(S)} + \sup_{S \in \mathscr F} \ln \frac{\nu(S)}{\mu(S)} = \infty.\]
    Hence~\eqref{eqPropDRMDefn} holds once again.

    Finally, suppose that $\mu$ and $\nu$ are non-zero and mutually absolutely continuous. Let $\tau$ be a dominating measure of $\mu$ and $\nu$ and suppose that $f$ and $g$ are $\tau$-densities of $\mu$ and $\nu$ respectively. By Lemma~\ref{lemmaEquivalenceBetweenSupEssSup},
    \begin{align*}
        \sup_{S \in \mathscr F} \ln \frac{\mu(S)}{\nu(S)} + \sup_{S \in \mathscr F} \ln \frac{\nu(S)}{\mu(S)} &= \esssup_{t \in \To} \ln \frac{f(t)}{g(t)} + \esssup_{t' \in \To} \ln \frac{g(t')}{f(t')} = \drm(\mu, \nu),
    \end{align*}
    where the essential supremum is over $\tau$.
\end{proof}

\begin{corollary}\label{corDRMCombinedTerms}
    For finite $\mu, \nu \in \Omega$,
    \begin{equation*}
        \drm(\mu, \nu) = \begin{cases} 0 &\mathrm{if\ } \mu = \nu = 0,\\
            \sup_{S, S' \in \F^*} \ln \frac{\mu(S) \nu(S')}{\mu(S') \nu(S)} &\mathrm{else\ if\ } \mu, \nu \mathrm{\ are\ mutually\ absolutely\ continuous,}\\
            \infty &\mathrm{otherwise,}
        \end{cases}
    \end{equation*}
    where $\F^* = \{ S \in \F : \mu(S) > 0\}$.
\end{corollary}

\begin{proof}
    Proving the result when $\mu = \nu = 0$ or when $\mu, \nu$ not mutually absolutely continuous is straightforward. We may thus assume that $\mu$ and $\nu$ are non-zero and mutually absolutely continuous. Then $\F^* = \{ S \in \F : \nu (S) > 0\}$ and moreover
    \[\ln \frac{\mu(S)}{\nu(S)} = 0,\]
    for all $S \notin \F^*$. Hence
    \[\drm(\mu, \nu) = \sup_{S \in \mathscr F^*} \ln \frac{\mu(S)}{\nu(S)} + \sup_{S \in \mathscr F^*} \ln \frac{\nu(S)}{\mu(S)}.\]
    The result then follows because $0 < \mu(S), \nu(S) < \infty$ for all $S \in \F^*$ implies we may combine the log-terms without introducing an undefined operation. That is,
    \[\sup_{S \in \mathscr F^*} \ln \frac{\mu(S)}{\nu(S)} + \sup_{S \in \mathscr F^*} \ln \frac{\nu(S)}{\mu(S)} = \sup_{S, S' \in \F^*} \ln \frac{\mu(S) \nu(S')}{\mu(S') \nu(S)},\]
    because $0 < \mu(S), \nu(S) < \infty$ for all $S \in \F^*$.
\end{proof}

\begin{corollary}\label{corDRMDCOR}
    For finite $\mu, \nu \in \Omega$,
    \[\drm(\mu, \nu) = - \ln \left[ 1 - \dcor (\mu, \nu) \right],\]
    where $\dcor$ is the constant odds metric (defined in Remark~\ref{remarkDRNCOR}).
\end{corollary}

\begin{proof}
    The cases where $\mu$ and $\nu$ are zero or not mutually absolutely continuous are straightforward to verify. When $\mu$ and $\nu$ are non-zero and mutually absolutely continuous,
    \[\frac{1}{-\dcor(\mu, \nu)+1} = \exp( \drm(\mu, \nu) ),\]
    by Corollary~\ref{corDRMCombinedTerms}.
\end{proof}

\subsection{The Multiplicative Distance \texorpdfstring{$\dMult$ }{}and the Density Ratio Metric \texorpdfstring{$\drm$ }{} Are Strongly Equivalent for Probabilities}

The following proposition proves that the multiplicative distance $\dMult$ and the density ratio metric $\drm$ are strongly equivalent (Definition~\ref{defnMetricStrongEquivalence}) -- but not equal -- on the set of probability measures. This proposition also shows that $\dMult$ and $\drm$ are not strongly equivalent on the set $\Omega$ of $\sigma$-finite measures.

See also Theorem~1 of \citet{wassermanComputingBoundsExpectations1992}, which describes the relationships between intervals of measures (or density bounded classes) and density ratio neighbourhoods.

\begin{proposition}\label{propInequalityMultDRM}
    Let $P, Q \in \Omega_1$ be probability measures on $(\T, \F)$. Then \begin{equation}\label{eqMultDRMEquivalentv2}
        \dMult(P,Q) \le \drm(P,Q) \le 2 \dMult(P,Q).
    \end{equation}
    Moreover, for each of the two inequalities in~\eqref{eqMultDRMEquivalentv2}:
    \begin{itemize}
        \item there exist $P, Q \in \Omega_1$ such that the inequality is strict; and
        \item there exist probabilities $P, Q \in \Omega_1$ (which can even be mutually absolutely continuous and not equal) such that the inequality is an equality.
    \end{itemize}
    Finally, when the probabilities $P$ and $Q$ are replaced by measures $\mu, \nu \in \Omega$, the first inequality of~\eqref{eqMultDRMEquivalentv2} does not hold (even up to a non-zero multiplicative constant), but the second does.
\end{proposition}

\begin{proof}
    Let $P$ and $Q$ be probability measures on $(\T, \F)$. Proving~\eqref{eqMultDRMEquivalentv2} when $P = Q$ is trivial by the metric properties of $\dMult$ and $\drm$. Hence, we may assume that $P \ne Q$. Then there exists $S \in \F$ such that $P(S) > Q(S)$ and consequently,
    \begin{equation}\label{eqLogLikelihoodPos}
        \sup_{S \in \F} \ln \frac{P(S)}{Q(S)} > 0.
    \end{equation}
    Symmetrically,
    \begin{equation}\label{eqLogLikelihoodPos2}
        \sup_{S \in \F} \ln \frac{Q(S)}{P(S)} > 0.
    \end{equation}
    Thus,
    \begin{align}
        \dMult (P, Q) &= \max \left( \sup_{S \in \F} \ln \frac{P(S)}{Q(S)}, \sup_{S \in \F} \ln \frac{Q(S)}{P(S)} \right) \notag\\
        &\le  \sup_{S \in \F} \ln \frac{P(S)}{Q(S)} + \sup_{S \in \F} \ln \frac{Q(S)}{P(S)} \label{eqWillBeStrict}\\
        &= \drm(P,Q),\notag
    \end{align}
    where the second line follows by~\eqref{eqLogLikelihoodPos} and~\eqref{eqLogLikelihoodPos2}, and the third by Proposition~\ref{propDRMDefn}.

    Also,
    \[\drm(P,Q) \le \sup_{S \in \F} \abs{\ln \frac{P(S)}{Q(S)}} + \sup_{S \in \F} \abs{\ln \frac{Q(S)}{P(S)}} = 2 \dMult(P,Q),\]
    by Proposition~\ref{propDRMDefn}.

Note that the first inequality, $\dMult(P,Q) \le \drm(P,Q)$, does not hold when $P, Q$ are replaced by (non-probability) measures. For example, for $\mu$ and $\nu$ defined in \eqref{eqDRMCounterexample}, $\drm(\mu, \nu) = 0$ while $\dMult(\mu, \nu) = \ln 2$. (It is easy to modify this example to construct $\mu$ and $\nu$ such that $0 < \drm(\mu, \nu) < \dMult(\mu, \nu)$. Simply change $\nu$ on some $E \ne \T$ so that $0 < \mu(E) < \nu(E) < 2 \mu(E)$ while still maintaining $\nu(E') = 2 \mu(E') > 0$ on some other $E'$.) However, the second inequality. $\drm(\mu,\nu) \le \dMult(\mu,\nu)$, does hold for any $\mu, \nu \in \Omega$, by the same proof as given above. 

The inequalities in \eqref{eqMultDRMEquivalentv2} are tight, even when $P$ and $Q$ are mutually absolutely continuous and $P \ne Q$: To see that the first inequality is tight, define the probability $P$ on $[-1,1]$ %
by the density
\[f(x) = \begin{cases}
    \frac{2}{3} & \mathrm{if\ } x \in [-1,0], \\
    \frac{2}{3}(1-x) &\mathrm{if\ } x \in (0,1],
\end{cases}\]
and $Q$ by
\[g(x) = \begin{cases}
    \frac{2}{3}(1+x) & \mathrm{if\ } x \in [-1,0], \\
    \frac{2}{3} &\mathrm{if\ } x \in (0,1],
\end{cases}.\]
Then $\dMult(P,Q) = \drm(P,Q) = \infty$. However, if $0 < \dMult(P,Q) < \infty$, then the inequality in~\eqref{eqWillBeStrict} will be strict and thus $\dMult(P,Q) < \drm(P,Q)$. 

To see that the second inequality in \eqref{eqMultDRMEquivalentv2} is tight, define $P,Q$ on $\{-1,1\}$ by $P(-1) = Q(1) = 1/3$ and $P(1) = Q(-1) = 2/3$. Then $\drm(P,Q) = 2 \ln 2 = 2 \dMult(P,Q)$. However, it is also possible that $\drm(P,Q) < 2 \dMult(P,Q)$: Define $P,Q$ on $\{1,2, \ldots, 9\}$ by
\begin{itemize}
    \item $P(1) = 9/10$;
    \item $P(x) = 1/80$ for $x = 2, 3, \ldots, 9$; 
    \item $Q(1) = 1/5$; and
    \item $Q(x) = 1/10$; for $x = 1, 2, \ldots, 10$.
\end{itemize}
Then $\drm(P,Q) = \ln 36 < 2 \ln 8 = 2\dMult(P,Q)$.
\end{proof}

\subsection{An Equivalence Between Intervals of Measures and Density Bounded Classes}

The following proposition provides an equivalence between intervals of measures and density bounded classes. Recall that $\mu \ll \nu$ denotes that $\mu$ is absolutely continuous with respect to $\nu$.

\begin{proposition}\label{propDRCIoM}
    Given $L, U, \nu \in \Omega$ with $L \le U$ and $U \ll \nu$, there exists some %
    $\nu$-densities $l \le u$ such that 
    \begin{enumerate}[label=(\alph*)]
        \item Every $\mu \in \IInt{L,U}$ has a $\nu$-density $f$ which is in the density bounded class $\IInt{l,u}$; and \label{propDRCIoMProp1}
        \item The measure $\mu(S) = \int_S f d\nu$ given by any density $f \in \IInt{l,u}$ is in the interval of measure $\IInt{L,U}$.\label{propDRCIoMProp2}
    \end{enumerate}
    Moreover, if $\nu = U$ then $u$ is the constant function: $u(t) = 1$ for all $t \in \mathcal T$; and $l$ is a $\nu$-density of $L$, which exists because $L$ is absolutely continuous with respect to $U$.

    Additionally, if $L = a\tau$ and $U = b \tau$ for some $\tau \in \Omega$ with $\tau \ll \nu$ and constants $0 < a \le 1 \le b < \infty$, then $l = a g$ and $u = b g$, where $g$ is a $\nu$-density of $\tau$.

    In the other direction, given some $\nu \in \Omega$ and $\nu$-densities $l \le u$ (which are finite $\nu$-almost everywhere), define $L, U \in \Omega$ by
    \begin{align}
        L(S) &= \int_S l d\nu, \label{eqDRCIoMDefnL}\\
        U(S) &= \int_S u d\nu, \label{eqDRCIoMDefnU}
    \end{align}
    for all $S \in \mathscr F$. Then $L \le U$ and the properties \ref{propDRCIoMProp1} and \ref{propDRCIoMProp2} above hold. 
    
\end{proposition}

Note that the condition that $l$ and $u$ are finite $\nu$-almost everywhere is necessary and sufficient for $L, U$ to be in $\Omega$ (see Lemma~\ref{lemmaDensityFiniteSigma}). 

The first half of this proposition (before ``In the other direction...'') shows that an interval of measure $\IInt{L,U}$ can be considered as a density bounded class $\IInt{l,u}$, and the second half shows the converse -- that a density bounded class $\IInt{l,u}$ can be considered as an interval of measure $\IInt{L,U}$.

\begin{proof}
    Let $L, U \in \Omega$ and suppose $L \le U$. We first consider the case that $\nu = U$. The fact that $L \le U$ implies that $L \ll U$: if $U(S) = 0$ then $L(S) = 0$, for all $S \in \mathscr F$. Thus, by the Radon-Nikodym theorem, $L$ has a density $l$ with respect to $U$. Define the density $u$ by $u(t) = 1$ for all $t \in \mathcal T$. It is straightforward to verify that $u$ is a density of $U$ with respect to $U$ (i.e. $u$ is a $U$-density of $U$). 

    Suppose, for contradiction, that there exists $E \in \mathscr F$ with $U(E) > 0$ and $l(t) > u(t)$ for all $t \in E$. Then 
    \[0 < U(E) = \int_E u dU < \int_E l dU = L(E),\]
    contradicting the assumption $L \le U$. Hence there exists $\tilde l$ such that $\tilde l(t) \le u(t)$ for all $t \in \mathcal T$ and $\tilde l(t) = l(t)$, for $U$-almost all $t \in \mathcal T$. This implies $\tilde l$ is also a $U$-density of $L$ because
    \[\int_S \tilde l dU = \int_S l dU = L(S),\]
    for all $S \in \mathscr F$. From herein, replace $l$ by its modification $\tilde l$. This proves $l \le u$.

    Now we prove \ref{propDRCIoMProp1}. Let $\mu \in \IInt{L,U}$. Since $\mu \le U$, we know that $\mu$ is absolutely continuous with respect to $U$ and hence has a $U$-density $f$ by the Radon-Nikodym theorem. We know that $f(t) \le u(t)$ for $U$-almost all $t \in \mathcal T$. Otherwise there exists $E \in \mathscr F$ with $U(E) > 0$ and $f(t) > u(t)$ for all $t \in E$, which would imply
    \[0 < U(E) = \int_E u dU < \int_E f dU = \mu(E),\]
    contradicting the assumption $\mu \le U$. Thus, $f \le u$, $U$-almost everywhere. This implies there exists $\tilde f$ which differs from $f$ on a $U$-null set such that $\tilde f(t) \le u(t)$ for all $t \in \mathcal T$. Since $f$ and $\tilde f$ differ only on a $U$-null set, $\tilde f$ is also a $U$-density of $\mu$. By exactly the same reasoning as above, we can show that $\tilde f \ge l$, $U$-almost everywhere and hence there exists a modification $\check f$ of $\tilde f$ such that 
    \[l(t) \le \check f(t) \le u(t),\]
    and $\check f\ne \tilde f$ only on a $U$-null set. Thus, $\check f$ is a $U$-density of $\mu$ such that $\check f \in \IInt{l,u}$. This proves \ref{propDRCIoMProp1}.

    Next we prove \ref{propDRCIoMProp2}. Fix some density $f \in \IInt{l,u}$ and define $\mu \in \Omega$ by 
    \[\mu(S) = \int_S f dU.\]
    Then 
    \[L(S) = \int_S l dU \le \int_S f dU = \mu(S),\]
    and 
    \[\mu(S) = \int_S f dU \le \int_S u dU = U(S).\]
    This proves $\mu \in \IInt{L,U}$.

    Now we consider the general case of the first half of the proposition. Let $L, U, \nu \in \Omega$ with $L \le U$ and $U \ll \nu$. Then $L \ll \nu$ as well, and let $l$ and $u$ be $\nu$-densities of $L$ and $U$ respectively. By the same reasoning as above, we can replace $l$ and $u$ with their modifications $\tilde l, \tilde u$. Hence we may assume that $l(t) \le u(t)$ for all $t \in \T$. Now take $\mu \in \IInt{L,U}$. Since $\mu \ll U \ll \nu$, the Radon-Nikodym theorem states that $\mu$ has a $\nu$-density $f$. Analogous to above, we can modify $f$ on a $\nu$-null set to produce a $\nu$-density $\check f$ of $\mu$ such that $\check f \in \IInt{l,u}$. This proves~\ref{propDRCIoMProp1}. The proof of~\ref{propDRCIoMProp2} is analogous to the case when $\nu = U$ with the integrals $\int \cdot\ dU$ replaced by $\int \cdot\ d\nu$.

    Finally, we consider when $L = a\tau$ and $U = b \tau$ for some $\tau \in \Omega$ with $\tau \ll \nu$ and constants $0 < a \le 1 \le b < \infty$. Define $l = a g$ and $u = b g$, where $g$ is some $\nu$-density of $\tau$. Then $l \le u$ and the proof of properties~\ref{propDRCIoMProp1} and~\ref{propDRCIoMProp2} are as before.

    We now prove the second half of the proposition (that which follows ``In the other direction...''). Let $\nu \in \Omega$ and suppose $l$ and $u$ are $\nu$-densities. Define $L, U \in \Omega$ according to equations~\eqref{eqDRCIoMDefnL} and~\eqref{eqDRCIoMDefnU}. Then $L \le U$ since 
    \[L(S) = \int_S l d\nu \le \int_S u d \nu = U(S).\]

    We will show that property \ref{propDRCIoMProp1} is satisfied. Let $\mu \in \IInt{L,U}$. By the same reasoning as in the proof of the first half of the proposition, modify the $\nu$-density $f$ of $\mu$ on a $\nu$-null set to produce $\check f$ satisfying $l(t) \le \check f(t) \le u(t)$. This implies $\check f$ is a $\nu$-density of $\mu$ and that $\check f \in \IInt{l,u}$.

    Finally, property \ref{propDRCIoMProp2} also follows by the same reasoning as in the proof of it in the first half of the proposition.
\end{proof}

\subsection{Characterising Mechanisms with Zero Privacy Loss}

The following proposition formalises the relationship between complete privacy and releasing pure noise, as described on page~\pageref{statementCompletePrivacy}:

\begin{proposition}\label{propCompletePrivacy}
    Let $M$ be a data-release mechanism. Denote the set of connected components of $\X$ by
    \[\ComX = \{ [x] : x \in \X\},\]
    where $[x] = \{ x' \in \X \mid d(x, x') < \infty\}$ (see Definition~\ref{defnConnected}).
    The following statements are equivalent:
    \begin{enumerate}[label=\Roman*]
        \item $M$ satisfies $\epsilon$-DP with $\epsilon = 0$.\label{propCompletePrivacyStatement1}
        \item $M$ is a function of $x$ only through its connected component $[x]$. That is, there exists a function $M' : \ComX \times [0,1] \to \T$ such that $M = M' \circ c$, where \label{propCompletePrivacyStatement2}
        \begin{align*}
            c : \mathcal X \times [0,1] &\to \ComX \times [0,1], \\
            (x, u) &\mapsto ([x], u).
        \end{align*}
        \item The probability $P_x$ induced by $M$ is a function of $x$ only through $[x]$. That is, $P_x = P_{x'}$ whenever $[x] = [x']$, or equivalently, the map $x \mapsto P_x$ factors as $\phi \circ c$, where 
        \begin{align*}
            c : \X &\to \ComX, \\
            x &\mapsto [x],
        \end{align*}
        and $\phi$ is some function which maps $[x]$ to a probability on $(\T, \F)$.\label{propCompletePrivacyStatement3}
    \end{enumerate}
\end{proposition}

\begin{proof}
    The equivalence between \ref{propCompletePrivacyStatement2} and \ref{propCompletePrivacyStatement3} follows by the definition of $P_x$ as the probability induced by $M$ -- see equation~\eqref{eqDefnPx}.

    We now prove that \ref{propCompletePrivacyStatement1} and \ref{propCompletePrivacyStatement3} are equivalent. The mechanism $M$ satisfies $\epsilon$-DP with $\epsilon = 0$ if and only if $\dMult(P_x, P_{x'}) = 0$ for all $x, x'$ with $d(x, x') < \infty$. Because metrics are positive definite (see Definition~\ref{defnMetric}.\ref{defnMetricProp1}), this is equivalent to $P_x = P_{x'}$ for all $x, x'$ with $d(x, x') < \infty$. But $P_x = P_{x'}$ holds for $x' \in [x]$ if and only if $P_x$ is a function of $x$ only through $[x]$.
\end{proof}

The following proposition formalises the result in Remark~\ref{remarkConnect}, which states that publishing $[x]$ alongside $T = M(x, U)$ does not increase the privacy loss but ensures that $\supp (x \mid t, \theta)$ is connected. (Recall that $[x]$ is the connected component $[x] = \{x' \in \X \mid d(x, x') < \infty\}$.)

\begin{proposition}\label{propMakeConnected}
    Let $M : \X \times [0,1] \to \T$ be an $\epsilon$-DP mechanism. Then the mechanism $M'$ defined by
    \begin{align*}
        M' : \X \times [0,1] &\to 2^{\X} \times \T, \\
        (x, u) &= \big( [x], M(x, u) \big),
    \end{align*}
    is also $\epsilon$-DP, and, moreover, the support $\supp(x \mid t, \theta)$ for $M'$ is $d$-connected for all the possible outputs $t \in \T \times 2^{\X}$ of $M'$ and all $\theta \in \Theta$.
\end{proposition}

\begin{proof}
    Proposition~\ref{propCompletePrivacy} implies that the data-release mechanism 
    \begin{align*}
        M_0 : \X \times [0,1] &\to \ComX, \\
        (x, u) &\mapsto [x],
    \end{align*}
    is $\epsilon$-DP with $\epsilon = 0$. Observe that $P_x( M' \in \cdot)$ is the product measure of $P_x(M_0 \in \cdot)$ and $P_x(M \in \cdot)$. Thus $M'$ is also $\epsilon$-DP by Proposition~\ref{propMultDistComp}. This proves the first half of the proposition.
    
    To see the second half of the proposition, fix some $\theta \in \Theta$ and some $t \in \T$. Let $E \subset \X$. If there does not exist some $x \in \X$ with $E = [x]$ then $\supp(x \mid (t, E), \theta) = \emptyset$ which is trivially $d$-connected. On the other hand, if $E = [x]$ for some $x \in \X$ then $\supp(x \mid (t, E), \theta) \subset [x]$. Since $[x]$ is $d$-connected by definition, $\supp(x \mid (t, E), \theta)$ must also be.
\end{proof}

\subsection{Alternative Characterisations of the Privatised Data Probability\texorpdfstring{ $P(T \in \cdot \mid \theta)$}{}}

The following proposition provides a number of characterisations of the privatised data probability $P(T \in \cdot \mid \theta)$, defined in equation~\eqref{eq:marginal-data-probability}. Recall that $\lambda$ is the Lebesgue measure. 

\begin{proposition}\label{propDefnPrivDataProb}
    Let $M$ be a data-release mechanism. Then the privatised data probability is given by:
    \begin{align*}
        P(T \in S \mid \theta) &= \int_{\mathcal X} P_x(S) dP_\theta(x), \\
        &= \int_{\mathcal X \times [0,1]} 1_{M(x, u) \in S} d(P_\theta \times \lambda)(x, u),
    \end{align*}
    for every $S \in \F$, where $P_\theta \times \lambda$ is the product measure of $P_\theta$ and $\lambda$, and $1_{M(x, u) \in S}$ is the indicator function:
    \[1_{M(x, u) \in S} = \begin{cases}
        1 &\mathrm{if\ } M(x, u) \in S, \\
        0 &\mathrm{otherwise.}
    \end{cases}\]
    Moreover, $P(T \in \cdot \mid \theta)$ is a well defined probability on $(\mathcal T, \mathscr F)$, for every $\theta \in \Theta$.
\end{proposition}

\begin{proof}
    Recall that $M$ is $(\mathscr G \otimes \mathcal B[0,1], \mathscr F)$-measurable. Hence $1_{M(x, u) \in S} : \X \times [0,1] \to \mathbb R$ is also measurable for every $S \in \F$. Thus,
    \[\int_{\mathcal X \times [0,1]} 1_{M(x, u) \in S} d(P_\theta \times \lambda)(x, u)\]
    is a well-defined probability on $(\T, \F)$. Then 
    \begin{align*}
        \int_{\mathcal X \times [0,1]} 1_{M(x, u) \in S} d(P_\theta \times \lambda)(x, u) &= \int_{\X} \left( \int_{[0,1]} 1_{M(x,u) \in S} du \right) dP_\theta(x) \\
        &= \int_{\X} \lambda \big( \{ u \in [0,1] : M(x,u) \in S\} \big) dP_\theta(x) \\
        &= \int_{\mathcal X} P_x(S) dP_\theta(x),
    \end{align*}
    where the first line follows by Fubini's theorem, and the second by the definition of the Lebesgue measure $\lambda$ and the third by the definition of $P_x(S)$ in equation~\eqref{eqDefnPx}. (Note that $x \mapsto P_x(S)$ is indeed $(\mathscr G, \F)$-measurable -- see \cite[Lemma~1.7.3]{durrettProbabilityTheoryExamples2019}.)
\end{proof}

\subsection{Alternative Characterisations of the Probability of a Data-Provision Procedure\texorpdfstring{ $M_G$}{}}

The following lemma will be useful in proving that the probability of a data-provision procedure $M_G$ is well defined. 

\begin{lemma}\label{lemmaFunctionsProductMeasurable}
    Let $(X, \Sigma_X), (Y, \Sigma_Y)$ and $(Z, \Sigma_Z)$ be measurable spaces. Suppose that $f : X \to Y$ is $(\Sigma_X, \Sigma_Y)$-measurable and $g : Y \times Z \to W$ is $(\Sigma_Y \otimes \Sigma_Z, \Sigma_W)$-measurable. 

    Then the map
    \begin{align*}
        h : X \times Z &\to W, \\
        (x,z) &\mapsto g(f(x), z),
    \end{align*}
    is $(\Sigma_X \otimes \Sigma_Z, \Sigma_W)$-measurable.
\end{lemma}

\begin{proof}
    Define
    \begin{align*}
        \phi : X \times Z &\to Y \times Z, \\
        (x,z) &\mapsto (f(x), z),
    \end{align*}
    and 
    \[\mathscr E = \{ E_2 \in \Sigma_Y \otimes \Sigma_Z \mid \exists E_1 \in \Sigma_X \otimes \Sigma_Z \text{ s.t. } \phi^{-1}(E_2) = E_1 \}.\]
    It is easy to verify that $\mathscr E$ is a $\sigma$-algebra because $\phi^{-1}(\emptyset) = \emptyset$; $\phi^{-1}(E_2^c) = [\phi^{-1}(E_2)]^c$; and $\phi^{-1}(\cup_i E_2^{(i)}) = \cup_i \phi^{-1} (E_2^{(i)})$. (Here $E_2^c$ is the complement of $E_2$ in $Y \times Z$.) Define the rectangles
    \[\mathscr R = \{ E \times F \mid E \in \Sigma_Y, F \in \Sigma_Z \}.\]
    Because $f$ is measurable, $\mathscr R \subset \mathscr E$ and hence $\Sigma_Y \otimes \Sigma_Z = \sigma(\mathscr R) = \mathscr E$. 

    Now let $E_3 \in \Sigma_W$. Because $g$ is measurable, there exists $E_2 \in \Sigma_Y \otimes \Sigma_Z$ such that $g^{-1}(E_3) = E_2$. Because $\Sigma_Y \otimes \Sigma_Z \subset \mathscr E$, there exists $E_1 \in \Sigma_X \otimes \Sigma_Z$ such that $\phi^{-1}(E_2) = E_1$. Thus, $h^{-1}(E_3) = \phi^{-1}(g^{-1}[E_3]) = E_1$.
\end{proof}

The following proposition provides a number of characterisations of the probability $P(T \in \cdot \mid \theta)$ of the output of a data-provision procedure $M_G$, defined in equation~\eqref{eq:marginal-data-probability2}.

\begin{proposition}\label{propDefnDataProvisionProb}
    Let $M_G$ be a data-provision procedure. The probability on $M_G$'s output $T$ is given by
    \begin{align}
        P(T \in S \mid \theta) &= \lambda \big( \{ u_1, u_2 \in [0,1] : M_G(\theta, u_1, u_2) \in S \} \big) \label{eqPropDefnDataProvisionProb1}\\
        &= \int_{[0,1]^2} 1_{M_G(\theta, u_1, u_2) \in S} d(u_1, u_2) \label{eqPropDefnDataProvisionProb2}\\
        &= \int_{\mathcal X} P_x(S) dP_{\theta}(x),\label{eqPropDefnDataProvisionProb3}
    \end{align}
    for $S \in \F$, where $1_{M_G(\theta, u_1, u_2) \in S}$ is the indicator function:
    \[1_{M_G(\theta, u_1, u_2) \in S} = \begin{cases}
        1 &\mathrm{if\ } M_G(\theta, u_1, u_2) \in S, \\
        0 &\mathrm{otherwise.}
    \end{cases}\]
    Moreover, for each $\theta \in \Theta$, $P(T \in \cdot \mid \theta)$ is a well-defined probability on $(\T, \F)$.
\end{proposition}

\begin{proof}
    We first show that $P(T \in \cdot \mid \theta)$ is a well-defined probability on $(\mathcal T, \mathscr F)$. By the definition of $P(T \in \cdot \mid \theta)$ as a pushforward probability in equation~\eqref{eq:marginal-data-probability2}, it suffices to show that the map
    \begin{align*}
        M_G(\theta, \cdot, \cdot) : [0,1]^2 &\to \mathcal T, \\
        (u_1, u_2) &\mapsto M_G(\theta, u_1, u_2),
    \end{align*}
    is $(\mathcal B[0,1]^2, \mathscr F)$-measurable for every $\theta \in \Theta$. (Here $\mathcal B[0,1]^2$ is the Borel $\sigma$-algebra on $[0,1] \times [0,1]$.) Recall that we assume that $M$ is $(\mathscr G \otimes \mathcal B[0,1], \mathscr F)$-measurable  and that $G(\theta, \cdot)$ is $(\mathcal B[0,1], \mathscr G)$-measurable for all $\theta \in \Theta$. Thus, measurability of $M_G$ follows from Lemma~\ref{lemmaFunctionsProductMeasurable} because $\mathcal B[0,1]^2 = \mathcal B[0,1] \otimes \mathcal B[0,1]$. 
    
    Now we prove equations~\eqref{eqPropDefnDataProvisionProb1}-\eqref{eqPropDefnDataProvisionProb3}. Equation~\eqref{eqPropDefnDataProvisionProb1} is simply the definition of $P(T \in S \mid \theta)$, as given in equation~\eqref{eq:marginal-data-probability2}. Equation~\eqref{eqPropDefnDataProvisionProb2} follows from the definition of the Lebesgue measure $\lambda$ on $[0,1]^2$. Equation~\eqref{eqPropDefnDataProvisionProb3} follows by the calculation:
    \begin{align*}
        \int_{[0,1]^2} 1_{M_G(\theta, u_1, u_2) \in S} d(u_1, u_2) &= \int_{[0,1]} \left( \int_{[0,1]} 1_{M(G(\theta,u_1), u_2) \in S} d u_2 \right) d u_1 \\
        &= \int_{[0,1]} \left( \int_{[0,1]} 1_{M(x, u_2) \in S} d u_2 \right) 1_{G(\theta,u_1) = x} d u_1 \\
        &= \int_{[0,1]} P_x(S) 1_{G(\theta,u_1) = x} d u_1 \\
        &= \int_{\mathcal X} P_x(S) dP_{\theta}(x),
    \end{align*}
    where the first line is Fubini's theorem; the second line is a substitution of $G(\theta, u_1)$ with $x$; the third line follows by the definition of $P_x$ in equation~\eqref{eqDefnPx}:
    \begin{align*}
        P_x(S) &= \lambda \big( \{ u_2 \in [0,1] : M(x,u_2) \in S\} \big) \\
        &= \int_{[0,1]} 1_{M(x,u_2) \in S} du_2,
    \end{align*}
    and the fourth line follows by the definition of $P_\theta$ in equation~\eqref{eqDefnProbPTheta}:
    \begin{align*}
        P_\theta(X \in E) &= \lambda \big( \{u_1 \in [0,1] : G(\theta, u_1) \in E\} \big) \\
        &= \int_{[0,1]} 1_{G(\theta, u_1) \in E} du_1.
    \end{align*}
    (Note that $P_x(S)$ is indeed $(\mathscr G, \F)$-measurable -- see \cite[Lemma~1.7.3]{durrettProbabilityTheoryExamples2019}.) 
\end{proof}

\subsection{Pufferfish Bounds the Prior-to-Posterior Odds Ratio}

The following proposition formalises the result described by equation~\eqref{eqPuffBayes}.

\begin{proposition}\label{propPuffOddsRatio}
    Fix a data-generating process $G$. Let $M$ be a data-release mechanism. Then $M$ satisfies \Puff\ if and only if
    \begin{equation}\label{eqPuffBayes2}
        e^{-\epsilon} \le \left. \frac{P_{\theta^*} (X \in E \mid T = t)}{P_{\theta^*} (X \in E' \mid T = t)} \middle/ \frac{P_{\theta^*} (X \in E)}{P_{\theta^*}(X \in E')} \right. \le e^\epsilon,
    \end{equation}
    for all $\theta^* \in \mathbb D$, all\footref{footnoteWellDefined} $(E, E') \in \Spairs$, and all $t \in \T_*$, where $\T_*^c$ is a null set under $P(T \in \cdot \mid \theta^*, X \in E)$ and under $P(T \in \cdot \mid \theta^*, X \in E')$. 
\end{proposition}

\begin{proof}
    Suppose throughout that $\epsilon < \infty$ (otherwise the proposition is vacuous). 
    
    ``$\Rightarrow$'': Suppose that $M$ satisfies \Puff. Fix some $\theta^* \in \mathbb D$ and some $(E, E') \in \Spairs$ such that $P_{\theta^*}(X \in \cdot \mid X \in E)$ and $P_{\theta^*}(X \in \cdot \mid X \in E')$ are both well-defined. By \ref{thm:pufferfish-iom}.\ref{puffThmStatement4}, there is a common dominating measure $\nu$ of $P(T \in \cdot \mid \theta^*, X \in E)$ and $P(T \in \cdot \mid \theta^*, X \in E')$. Additionally, $P(T \in \cdot \mid \theta^*, X \in E)$ and $P(T \in \cdot \mid \theta^*, X \in E')$ have $\nu$-densities satisfying
    \[e^{-\epsilon} \le \frac{p(t \mid \theta^*, X \in E)}{p(t \mid \theta^*, X \in E)} \le e^{\epsilon},\]
    for all $t \in \T$. Now Bayes rule states that 
    \begin{equation}\label{eqProofPropPuffOddsRatio1}
        P_{\theta^*} (X \in E \mid T = t) \propto P_{\theta^*} (X \in E) p(t \mid \theta^*, X \in E),
    \end{equation}
    for $P(T \in \cdot \mid \theta^*, X \in E)$-almost all $t$. Let $\T_*^{(1)}$ be the null set where~\eqref{eqProofPropPuffOddsRatio1} doesn't hold. Similarly, let $\T_*^{(2)}$ be the $P(T \in \cdot \mid \theta^*, X \in E')$-null set where 
    \[P_{\theta^*} (X \in E' \mid T = t) \propto P_{\theta^*}(X \in E') p(t \mid \theta^*, X \in E'),\]
    does not hold. On $\T_* = (\T_1^* \cup \T_2^*)^c$, the posterior odds is equal to product of the prior odds and the likelihood ratio:
    \[\frac{P_{\theta^*} (X \in E \mid T = t)}{P_{\theta^*} (X \in E' \mid T = t)} = \frac{P_{\theta^*} (X \in E)}{P_{\theta^*}(X \in E')} \frac{p(t \mid \theta^*, X \in E)}{p(t \mid \theta^*, X \in E')}.\]
    Hence~\eqref{eqPuffBayes2} holds for all $t \in \T_*$. 

    ``$\Leftarrow$'': Fix some $\theta^* \in \mathbb D$ and some $(E, E') \in \Spairs$ such that $P_{\theta^*}(X \in \cdot \mid X \in E)$ and $P_{\theta^*}(X \in \cdot \mid X \in E')$ are both well-defined. Suppose that~\eqref{eqPuffBayes2} holds. We will show that \ref{thm:pufferfish-iom}.\ref{puffThmStatement2} must also hold. Firstly, we note that $P(T \in \cdot \mid \theta^*, X \in E)$ is absolutely continuous with respect to $P(T \in \cdot \mid \theta^*)$. Also,
    \[\frac{P_{\theta^*}(X \in E \mid T = t)}{P_{\theta^*}(X \in E)}\]
    is a $P(T \in \cdot \mid \theta^*)$-density for $P(T \in \cdot \mid \theta^*, X \in E)$ because
    \[P(T \in S \mid \theta^*, X \in E) = \int_S \frac{P_{\theta^*}(X \in E \mid T = t)}{P_{\theta^*}(X \in E)} dP(t \mid \theta^*)\]
    by Bayes rule. All the above also applies when $E$ is replaced by $E'$. Hence, for any $S \in \F$, 
    \begin{align*}
        P(T \in S \mid \theta^*, X \in E) &= P(T \in S \cap \T_* \mid \theta^*, X \in E) \\
        &= \int_{S \cap \T_*} \frac{P_{\theta^*}(X \in E \mid T = t)}{P_{\theta^*}(X \in E)} dP(t \mid \theta^*) \\
        &\le e^\epsilon \int_{S \cap \T_*} \frac{P_{\theta^*}(X \in E' \mid T = t)}{P_{\theta^*}(X \in E')} dP(t \mid \theta^*) \\
        &= e^\epsilon P(T \in S \cap \T_* \mid \theta^*, X \in E') \\
        &= e^\epsilon P(T \in S \mid \theta^*, X \in E').
    \end{align*}
    The second half of \ref{thm:pufferfish-iom}.\ref{puffThmStatement2} follows analogously. 
\end{proof}

\section{Proofs Omitted From the Main Text}\label{appendixProofs}

Recall that $\mu \ll \nu$ denotes that $\mu$ is absolutely continuous with respect to $\nu$.

\begin{proof}of Lemma~\ref{lemmaMultIoM}:
    Fix two $\sigma$-finite measures $\mu, \nu \in \Omega$ and a constants $\epsilon > 0$.

    ``$\Rightarrow$'': Suppose that $\nu \in \IInt{e^{-\epsilon} \mu, e^{\epsilon} \mu}$. Then, for all $S \in \mathscr F$,
    \[\nu(S) \ge e^{-\epsilon} \mu(S),\]
    and hence
    \[\ln \frac{\mu(S)}{\nu(S)} \le \epsilon.\]
    Similarly,
    \[-\ln \frac{\mu(S)}{\nu(S)} = \ln \frac{\nu(S)}{\mu(S)} \le \epsilon,\]
    since $\nu(S) \le e^{\epsilon} \mu(S)$. Putting these two results together,
    \[\dMult(\mu, \nu) = \sup_{S \in \mathscr F} \abs{\ln \frac{\mu(S)}{\nu(S)}} \le \epsilon.\]

    ``$\Leftarrow$'': Suppose that $\dMult(\mu, \nu) \le \epsilon$. Then, for all $S \in \mathscr F$,
    \[\frac{\mu(S)}{\nu(S)} \le e^{\epsilon}.\]
    This proves $e^{-\epsilon} \mu \le \nu$. We also have that
    \[- \ln \frac{\mu(S)}{\nu(S)} \le \epsilon,\]
    for all $S \in \mathscr F$, which implies $\nu \le e^{\epsilon} \mu$. Thus, $\nu \in \IInt{e^{-\epsilon} \mu, e^{\epsilon} \mu}$.

    Now we prove the second half of the lemma. Fix two constants $0 < a \le 1 \le b < \infty$. Suppose $\nu \in \IInt{a \mu, b \mu}$. Let $\epsilon = \max (-\ln a, \ln b)$. Since $\IInt{a \mu, b \mu} \subset \IInt{e^{-\epsilon} \mu, e^{\epsilon} \mu}$, we have $\dMult \le \epsilon$ by the result of the first half of the lemma. Now let $\epsilon = \min (-\ln a, \ln b)$ and suppose $\dMult(\mu, \nu) \le \epsilon$. Then $\nu \in \IInt{e^{-\epsilon} \mu, e^{\epsilon} \mu} \subset \IInt{a \mu, b\mu}$ as required.
\end{proof}

\begin{proof}of Theorem~\ref{thm:dp-iom}:
    First we prove that \ref{mainThmStatement1} implies \ref{mainThmStatement2}. Suppose that $M$ satisfies pure $\epsilon$-DP. Fix some $S \in \mathscr F$ and some $x, x' \in \mathcal X$ with $d(x, x') = 1$. By assumption, $\dMult(P_x, P_{x'}) \le \epsilon$. This implies $P_{x'}(S) \le e^\epsilon P_x(S)$ by Lemma~\ref{lemmaMultIoM}, which is exactly \ref{mainThmStatement2}.

    Next we will prove that \ref{mainThmStatement2} implies \ref{mainThmStatement1}. To do this, we need the following lemma ($*$): \ref{mainThmStatement2} implies that $\dMult(P_x, P_{x'}) \le \epsilon$ for all $x, x' \in \mathcal X$ with $d(x, x') = 1$.

    To prove this lemma, suppose \ref{mainThmStatement2} holds and fix some $x, x' \in \mathcal X$ with $d(x, x') = 1$. By assumption $P_{x'} (S) \le e^\epsilon P_x(S)$ for all $S \in \mathscr F$. Yet, by symmetry of $d$ (i.e. because $d(x', x) = d(x, x') = 1$), \ref{mainThmStatement2} also implies that $P_{x} (S) \le e^\epsilon P_{x'} (S)$ for all $S \in \mathscr F$. Then Lemma~\ref{lemmaMultIoM} provides the desired result: $\dMult(P_x, P_{x'}) \le \epsilon$.
    
    Now we return to proving that \ref{mainThmStatement2} implies \ref{mainThmStatement1}. Fix $x, x' \in \mathcal X$. If $d(x, x') = \infty$ then the condition $\dMult(P_x, P_{x'}) \le \epsilon d(x, x')$ is vacuous. Similarly, if $d(x, x') = 0$ then, by the properties of $d$ as a metric, $x = x'$. Thus $P_x = P_{x'}$ and $\dMult(P_x, P_{x'}) = 0$. 
    
    On the other hand, if $d(x, x') = n < \infty$, then by Assumption~\ref{assumptionGraphDistance}, there exists $x = x_0, x_1, \ldots, x_n = x' \in \mathcal X$ with $d(x_i, x_{i+1}) = 1$ for all $0 \le i \le n-1$. Then,
    \begin{equation*}
        \dMult(P_{x}, P_{x'}) \le \sum_{i=0}^{n-1} \dMult(P_{x_i}, P_{x_{i+1}}) 
        \le \sum_{i=0}^{n-1} \epsilon 
        = \epsilon d(x, x'),
    \end{equation*}
     where the first line follows by the triangle inequality (Lemma~\ref{lemmaMultMetric}) of $\dMult$, and the second by $(*)$. This proves that \ref{mainThmStatement2} implies \ref{mainThmStatement1}.

     Now we prove that \ref{mainThmStatement1} is equivalent to \ref{mainThmStatement3}. Fix some $x, x' \in \mathcal X$ with $\delta = d(x, x')$. Then $M$ is $\epsilon$-DP if and only if $\dMult(P_x, P_{x'}) \le \delta \epsilon$ for all such $x, x' \in \mathcal X$. Yet Lemma~\ref{lemmaMultIoM} states that this is equivalent to $P_{x'} \in \IInt{e^{-\delta \epsilon} P_x, e^{\delta \epsilon} P_x}$, which is \ref{mainThmStatement3}.

     We move to proving that \ref{mainThmStatement3} implies \ref{mainThmStatement4}. This follows by Proposition~\ref{propDRCIoM}. Fix $x, x' \in \X$ and $\nu \in \Omega$. Suppose that $x'$ and $x$ are $d$-connected, so that $d(x, x') = \delta < \infty$. Assume that \ref{mainThmStatement3} holds for these $x, x'$. Suppose that $P_x$ has a $\nu$-density $p_x$. By the Radon-Nikodym theorem, this implies $P_x \ll \nu$. %
     Then we can apply Proposition~\ref{propDRCIoM}\ref{propDRCIoMProp1} with $\tau = P_x$, $a = \exp(-\delta \epsilon)$, $b = \exp( \delta \epsilon)$ and $g = p_x$. Then $P_{x'}$ has a $\nu$-density which is in the density bounded class $\IInt{\exp(-\delta \epsilon) p_x, \exp(\delta \epsilon) p_x}$. This is precisely \ref{mainThmStatement4}.

     Finally, we prove that \ref{mainThmStatement4} implies \ref{mainThmStatement3}. Set $\nu = P_x$ so that $P_x$ has a constant $\nu$-density $p_x = 1$. Then \ref{mainThmStatement4} implies that $p_{x'} \in \IInt{l,u}$, where $l = \exp (-\delta \epsilon)$ and $u = \exp(\delta \epsilon)$ are also constant $\nu$-densities. Apply the second half of Proposition~\ref{propDRCIoM}. This states that $P_{x'} \in \IInt{L,U}$ where $L$ and $U$ are defined as 
     \[L(S) = \int_S e^{-\delta \epsilon} dP_x = e^{-\delta \epsilon} P_x(S),\]
     and 
     \[U(S) = \int_S e^{\delta \epsilon} dP_x = e^{\delta \epsilon} P_x(S).\]
     Yet this is exactly \ref{mainThmStatement3}. 
\end{proof}

We now turn to proving Theorem~\ref{thmMarginalData}. We first need to establish some lemmata.

\begin{lemma}\label{lemmaPxT0Zero}
    Consider the same set-up as in Theorem~\ref{thmMarginalData}. If $x \in \supp_0 ( x \mid t_0 )^c \cap \supp (P_\theta)$, then $P_x(\mathcal T_0) = 0$.
\end{lemma}

\begin{proof}
    For $t \in \T_0$, we have
    \begin{align*}
        \supp_0 ( x \mid t_0 )^c \cap \supp (P_\theta) &\subset \left( \supp (P_\theta) \cap \supp_0 (x \mid t_0) \right)^c \cap \supp (P_\theta) \\
        &\subset \left( \supp (P_\theta) \cap \supp_0 (x \mid t) \right)^c \cap \supp (P_\theta) \\
        &\subset \supp_0 (x \mid t)^c \\
        &= \{x \in \X \mid t \in \supp (P_x)^c\}.
    \end{align*}
    Therefore, for $x \in \supp_0 ( x \mid t_0 )^c \cap \supp (P_\theta)$, we have that $\T_0 \in \supp (P_x)^c$. The result then follows by~\eqref{eqZeroOutsideSupport}.
\end{proof}

\begin{lemma}\label{lemmaPxT0Zero2}
    Consider the same set-up as in Theorem~\ref{thmMarginalData}. Then %
    \[\int_{\supp (x \mid t_0, \theta)^c} P_x(\T_0) d P_\theta(x) = 0.\]
\end{lemma}

\begin{proof}
    Compute
    \begin{align*}
        \int_{\supp (x \mid t_0, \theta)^c} P_x(\T_0) d P_\theta(x) %
        &= \int_{\supp (P_\theta)^c} P_x(\T_0) d P_\theta(x) + \int_{\supp_0 (x \mid t_0)^c \cap \supp(P_\theta)} P_x(\T_0) d P_\theta(x) \\
        &= \int_{\supp (P_\theta)^c} P_x(\T_0) d P_\theta(x) \\
        &\le P_\theta(\supp (P_\theta)^c) \\
        &= 0,
    \end{align*}
    where the second line follows by Lemma~\ref{lemmaPxT0Zero} and the fourth by~\eqref{eqZeroOutsideSupport}. 
\end{proof}

\begin{lemma}\label{lemmaMarginalData}
    Consider the same set-up as in Theorem~\ref{thmMarginalData}. Then, a density $p(t \mid \theta)$ for $P(T \in \cdot \mid \theta)$ exists in $\T_0 = \{t \in \T \mid \supp (x \mid t, \theta) \subset \supp (x \mid t_0, \theta)\}$ in the sense that
    \[P(T \in S \cap \T_0 \mid \theta) = \int_{S \cap \T_0} p(t \mid \theta) d\nu(t),\]
    for all $S \in \F$, where $\nu$ is the dominating measure for the density $p(t \mid \theta)$. Moreover, this density satisfies
    \begin{equation}\label{eqThmMarginalData}
        p ( t \mid \theta) \in p_{x_*}(t) \exp \left( \pm \epsilon d_* \right),
    \end{equation}
    for every $t \in \T_0$ and every $x_* \in \supp (x \mid t_0, \theta)$, where $d_{*} = \sup_{x \in \supp (x \mid t_0, \theta)} d(x, x_*)$.
\end{lemma}

\begin{proof}%
    Fix some $x_* \in \supp (x \mid t_0, \theta)$ and some $\nu \in \Omega$ with $P_{x_*} \ll \nu$. Let $p_{x_*}$ be a $\nu$-density of $P_{x_*}$. By Theorem~\ref{thm:dp-iom}.\ref{mainThmStatement4} and the assumption that $\supp(x \mid t_0, \theta)$ is $d$-connected, the probability $P_x$ also has a $\nu$-density $p_x$, for all $x \in \supp (x \mid t_0, \theta)$. For $t \in \T_0$, define
    \[p(t \mid \theta) = \int_{\supp (x \mid t_0, \theta)} p_x(t) d P_\theta(x).\]
    
    We want to prove that $p(t \mid \theta)$ is a $\nu$-density of $P(T \in \cdot \mid \theta)$ in $\T_0$ -- namely, that
    \[P(T \in S \cap \T_0 \mid \theta) = \int_{S \cap \T_0} p(t \mid \theta) d\nu(t),\]
    for all $S \in \mathscr F$. 
    
    We have
    \begin{align*}
    	\int_{S \cap \T_0} p(t \mid \theta) d\nu(t) &= \int_{S \cap \T_0} \left( \int_{\supp (x \mid t_0, \theta)} p_x(t) d P_\theta(x) \right) d\nu(t) \\
    	&= \int_{\supp (x \mid t_0, \theta)} \left( \int_{S \cap \T_0} p_x(t) d\nu(t) \right) d P_\theta(x) \\
    	&= \int_{\supp (x \mid t_0, \theta)} P_x(S \cap \T_0) d P_\theta(x) \\
		&= P(T \in S \cap \T_0 \mid \theta),
    \end{align*}
    where the second line follows by Fubini's theorem and the fourth by Lemma~\ref{lemmaPxT0Zero2}.

    Now we move to proving the upper bound of~\eqref{eqThmMarginalData}:
	\begin{align*}
		p(t \mid \theta) &= \int_{\supp (x \mid t_0, \theta)} p_x(t) dP_\theta(x) \\
		&\le \int_{\supp (x \mid t_0, \theta)} e^{\epsilon d(x, x_*)} p_{x_*}(t) d P_\theta (x) \\
		&\le e^{\epsilon d_{*}} p_{x_*}(t),
	\end{align*}
    where the second line follows from Theorem~\ref{thm:dp-iom}.\ref{mainThmStatement4}. The lower bound follows similarly.
\end{proof}

\begin{proof}of Theorem~\ref{thmMarginalData}: By Theorem~\ref{thm:dp-iom}, we can fix a measure $\nu \in \Omega$ which dominates all $P_x$ for $x \in \supp (x \mid t_0, \theta)$. Let $p_x$ denote a $\nu$-density of $P_x$ for $x \in \supp(x \mid t_0, \theta)$. Take the essential supremum with respect to $\nu$ over the collection of densities $\{\exp \left( - \epsilon d_{*} \right) p_{x_*}(t) : x_* \in \supp (x \mid t_0, \theta)\}$ to produce
    \begin{equation}\label{eqAppendixEssSup}
        l_{\theta, \epsilon}(t) = \esssup_{x_* \in \supp (x \mid t_0, \theta)} \exp \left( - \epsilon d_{*} \right) p_{x_*}(t).
    \end{equation}
    This function $l_{\theta, \epsilon} : \T \to [0,\infty]$ exists and is measurable as $\nu$ is $\sigma$-finite. %
    Thus, $l_{\theta, \epsilon}$ is a $\nu$-density for some measure $L_{\theta, \epsilon}$ on $(\T, \F)$. By Lemma~\ref{lemmaMarginalData}, we can construct a $\nu$-density $p(t \mid \theta)$ for $P(T \in \cdot \cap \T_0 \mid \theta)$ which satisfies
    \[p(t \mid \theta) \ge l_{\theta, \epsilon}(t),\]
    for $\nu$-almost all $t \in \T_0$. This proves that $P(T \in S \mid \theta) \ge L_{\theta, \epsilon}(S)$ for any measurable $S \subset \T_0$.
    
    Since $P(T \in \cdot \cap \T_0 \mid \theta)$ is zero outside of $\T_0$, technically we must modify $l_{\theta, \epsilon}$ to also be zero outside of $\T_0$. It then follows that $L_{\theta, \epsilon} (S) \le P(T \in S \cap \T_0 \mid \theta)$ for all $S \in \F$. This proves the lower bound of~\eqref{eqThmMarginalIoM}.
    
    The argument for the upper measure $U_{\theta, \epsilon}$ is almost analogous. We can construct $u_{\theta, \epsilon}$ on $\T$ as the essential infimum
    \[u_{\theta, \epsilon}(t)=\essinf_{x_* \in \supp (x \mid t_0, \theta)} \exp \left( \epsilon d_{*} \right) p_{x_*}(t).\]
    (Note that it is possible that $u_{\theta, \epsilon}$ is not finite almost everywhere -- even though all of the $p_{x^*}(t)$ are -- so that the measure $U_{\theta, \epsilon}$ is not $\sigma$-finite by Lemma~\ref{lemmaDensityFiniteSigma}.) Then Lemma~\ref{lemmaMarginalData} implies that
    \[p(t \mid \theta) \le u_\epsilon(t),\]
    for $\nu$-almost all $t \in \T_0$. This proves that $P(T \in S \mid \theta) \le U_{\theta, \epsilon}(S)$ for any measurable $S \subset \T_0$. Thus, $P(T \in S \cap \T_0 \mid \theta) \le U_{\theta, \epsilon} (S)$ for any $S \in \F$. This proves the upper bound of~\eqref{eqThmMarginalIoM}.
\end{proof}

\begin{proof}of Theorem~\ref{thmDPHypothesisTesting}: Fix some $x_* \in S_0 \cup S_1$ and some $\nu \in \Omega$ with $P_{x_*} \ll \nu$. By Theorem~\ref{thm:dp-iom}.\ref{mainThmStatement4} and the assumption that $S_0 \cup S_1$ is $d$-connected, the probability $P_x$ has a $\nu$-density $p_x$, for all $x \in S_0 \cup S_1$. For $t \in \T$, define
\[p(t \mid \theta_i) = \int_{S_i} p_x(t) dP_{\theta_i} (x).\]
We can show that $p(t \mid \theta_i)$ is a $\nu$-density of $P(T \in \cdot \mid \theta_0)$: For any $E \in \F$,
\begin{align*}
    \int_E p(t \mid \theta_i) d\nu (t) &= \int_E \int_{S_i} p_x(t) dP_{\theta_i} d\nu (t) \\
    &= \int_{S_i} \int_E p_x(t) d\nu (t) dP_{\theta_i} \\
    &= \int_{S_i} P_x(E) dP_{\theta_i} \\
    &= P(T \in E \mid \theta_i),
\end{align*}
where the second line follows by Fubini's theorem and the fourth by~\eqref{eqZeroOutsideSupport}. 

Let $R$ be the rejection region of a test with size $P(T \in R \mid \theta_0) \le \alpha$. Then
\begin{align*}
    P(T \in R \mid \theta_1) &= \int_R p(t \mid \theta_1) d\nu(t) \\
    &\le \exp ( d_{**} \epsilon) \int_R p(t \mid \theta_0) d\nu(t) \\
    & \le \alpha \exp ( d_{**}\epsilon),
\end{align*}
where the second line follows by the computation:
\begin{align*}
    p(t \mid \theta_1) &= \int_{S_1} p_x(t) dP_{\theta_1}(x) \\
    &= \int_{S_0} \left( \int_{S_1} p_x(t) dP_{\theta_1}(x) \right) dP_{\theta_0}(x') \\
    &\in \int_{S_0} \left( \int_{S_1} \exp \left( \pm \epsilon d_{**} \right) p_{x'}(t) dP_{\theta_1}(x) \right) dP_{\theta_0}(x') \\
    &= \exp \left( \pm \epsilon d_{**} \right) \int_{S_0} \left( \int_{S_1} dP_{\theta_1}(x) \right) p_{x'}(t) dP_{\theta_0}(x') \\
    &= \exp \left( \pm \epsilon d_{**} \right) p(t \mid \theta_0).
\end{align*}
In the above computation, the third line follows from Theorem~\ref{thm:dp-iom}.\ref{mainThmStatement4} and the other lines simply pull constant factors into -- or out of -- some integral over $dP_{\theta_i}$.
\end{proof}

We begin proving Theorem~\ref{thmPriorPredictive} by establishing the following lemma (which is analogous to Lemma~\ref{lemmaMarginalData}).

\begin{lemma}\label{lemmaPriorPredictive}
    Consider the same set-up as in Theorem~\ref{thmPriorPredictive}. Then, a density $p(t)$ for the prior predictive probability $P(T \in \cdot)$ exists in $\T_0$ in the sense that
    \begin{equation}\label{eqLemmaPriorPredictiveDensity}
        P(T \in S \cap \T_0) = \int_{S \cap \T_0} p(t) d\nu(t),
    \end{equation}
    for all $S \in \F$, where $\nu$ is the dominating measure for the density $p(t)$. Moreover, this density satisfies
    \begin{equation}\label{eqThmPriorPredictive}
        p ( t ) \in p_{x_*}(t) \exp \left( \pm \epsilon d_* \right),
    \end{equation}
    for every $t \in \T_0$ and every $x_* \in \supp (x \mid t_0)$, where $d_{*} = \sup_{x \in \supp (x \mid t_0)} d(x, x_*)$.
\end{lemma}

\begin{proof}
    We proceed as for the proof of Lemma~\ref{lemmaMarginalData}. Fix some $x_* \in \supp (x \mid t_0$ and some $\nu \in \Omega$ with $P_{x_*} \ll \nu$. Let $p_{x_*}$ be a $\nu$-density of $P_{x_*}$. By Theorem~\ref{thm:dp-iom}.\ref{mainThmStatement4} and the assumption that $\supp(x \mid t_0)$ is $d$-connected, the probability $P_x$ also has a $\nu$-density $p_x$, for all $x \in \supp (x \mid t_0)$. For $t \in \T_0$, define
    \[p(t) = \int_{\Theta} \left( \int_{\supp (x \mid t_0)} p_x(t) d P_\theta(x) \right) d\pi(\theta).\]
    We can compute
    \begin{align*}
        \int_{S \cap \T_0} p(t) d\nu(t) &= \int_{S \cap \T_0} \int_{\Theta} \int_{\supp (x \mid t_0)} p_x(t) d P_\theta(x) d\pi(\theta) d\nu(t) \\
        &= \int_\Theta \int_{\supp (x \mid t_0)} \int_{S \cap \T_0} p_x(t) d\nu(t) d P_\theta(x) d\pi(\theta) \\
        &= \int_\Theta \int_{\supp (x \mid t_0)} P_x(S \cap \T_0) d P_\theta(x) d\pi(\theta) \\
        &= \int_\Theta P(T \in S \cap \T_0 \mid \theta) d\pi(\theta) \\
        &= P(T \in S \cap \T_0),
    \end{align*}
    where the second line follows by Fubini's theorem and the fourth by Lemma~\ref{lemmaPxT0Zero2}. This proves~\eqref{eqLemmaPriorPredictiveDensity}.

    Now we move to proving the upper bound of~\eqref{eqThmPriorPredictive}:
	\begin{align*}
		p(t) &= \int_\Theta \left( \int_{\supp (x \mid t_0)} p_x(t) dP_\theta(x) \right) d\pi(\theta) \\
		&\le \int_\Theta \left( \int_{\supp (x \mid t_0)} e^{\epsilon d(x, x_*)} p_{x_*}(t) d P_\theta (x) \right) d\pi(\theta) \\
		&\le e^{\epsilon d_{*}} p_{x_*}(t),
	\end{align*}
    where the second line follows from Theorem~\ref{thm:dp-iom}.\ref{mainThmStatement4}. The lower bound of~\eqref{eqThmPriorPredictive} follows similarly.
\end{proof}

\begin{proof}of Theorem~\ref{thmPriorPredictive}: This is exactly analogous to the proof of Theorem~\ref{thmMarginalData}, where references to Lemma~\ref{lemmaMarginalData} are replaced by references to Lemma~\ref{lemmaPriorPredictive}; $p(t \mid \theta)$ is replaced with $p(t)$; and $\supp(x \mid t_0, \theta)$ is replaced with $\supp(x \mid t_0)$.
\end{proof}

\begin{proof}of Theorem~\ref{thmPosterior}: By Theorem~\ref{thm:dp-iom}, we can fix a measure $\nu \in \Omega$ which dominates all $P_x$ for $x \in \supp(x \mid t_0$. Let $p_x$ denote a $\nu$-density of $P_x$ for $x \in \supp(x \mid t_0)$. Define
\[p(t \mid \theta) = \int_{\supp (x \mid t_0)} p_x(t) dP_{\theta} (x).\]
Exactly as in the proof of Lemma~\ref{lemmaMarginalData}, we can show that $p(t \mid \theta)$ is a $\nu$-density of $P(T \in \cdot \mid \theta)$ in $\T_0$. Since $t_0 \in \T_0$, we can use this density to compute the posterior via Bayes rule:
\[\pi(\theta \mid t_0) = \frac{ p(t_0 \mid \theta) \pi(\theta) }{\int_{\Theta} p(t_0 \mid \theta') d\pi(\theta') }.\]
We have that
\begin{align}
    \frac{ p(t_0 \mid \theta) \pi(\theta) }{\int_{\Theta} p(t_0 \mid \theta') d\pi(\theta') } &\le \frac{ p(t_0 \mid \theta) \pi(\theta) }{\int_{\Theta} \exp (-\epsilon d_{**} ) p(t_0 \mid \theta) d\pi(\theta') } \notag\\
    &= \exp ( \epsilon d_{**} ) \pi(\theta),\label{eqProofThmPosterior}
\end{align}
where the first line follows by the calculations
\begin{align*}
    p(t_0 \mid \theta') &= \int_{\supp (x \mid t_0)} p_x(t) dP_{\theta'} (x) \\
    &= \int_{\supp (x \mid t_0)} \left( \int_{\supp (x \mid t_0)} p_x(t) dP_{\theta'} (x) \right) dP_{\theta}(x') \\
    &\le \int_{\supp (x \mid t_0)} \left( \int_{\supp (x \mid t_0)} \exp \left( \epsilon d_{**} \right) p_{x'}(t) dP_{\theta'} (x) \right) dP_{\theta}(x') \\
    &= \exp \left( \epsilon d_{**} \right) \int_{\supp (x \mid t_0)} p_{x'}(t) dP_{\theta}(x') \\
    &= p(t_0 \mid \theta).
\end{align*}
In the above computation, the third line follows from Theorem~\ref{thm:dp-iom}.\ref{mainThmStatement4} and the other lines simply pull constant factors into -- or out of -- some integral over $dP_{\theta}(x')$ or $dP_{\theta'}(x)$.

Using~\eqref{eqProofThmPosterior}, we obtain the upper bound of~\eqref{eq:posterior}:
\begin{align*}
    \pi( \theta \in S \mid t_0) &= \int_S \pi(\theta \mid t_0) d\mu(\theta) \\
    &\le \int_S \exp ( \epsilon d_{**} ) \pi(\theta) d\mu(\theta) \\
    &= \exp ( \epsilon d_{**} ) \pi(\theta \in S),
\end{align*}
where $\mu$ is the dominating measure of the prior density $\pi(\theta)$. The proof of the lower bound of~\eqref{eq:posterior} is analogous.
\end{proof}

\acks{%
We thank the four anonymous reviewers of \cite{bailieDifferentialPrivacyGeneral2023} and the two anonymous reviewers of this paper for their generous suggestions, which have greatly improved this work. JB gratefully acknowledges partial financial support from the Australian-American Fulbright Commission and the Kinghorn Foundation.
}

\section*{Author Contributions}
Both authors conceived the project, and contributed to the research and to the writing of the manuscript.

\section*{Data and Code Availability}
No data was created nor analysed in this work. The code for generating the figures in this article is publicly available here: \url{https://github.com/jameshbailie/BailieGongISIPTA23}.

\bibliography{bailie23.bib}

\end{document}